\documentclass{amsart}[a4paper]
\usepackage{preambleFL}
\usepackage{graphicx} 
\usepackage[all]{xy}
\usepackage{enumerate}

\title{Infinity operad pairs, higher distributivity and spectral Tambara functors}
\author{Aowen Fan}\thanks{Address: fanaowen@amss.ac.cn}
\date{}

\begin{document}
\begin{abstract}
We introduce an $\infty$-categorical model of operad pairs and their associated bimonoidal categories. We study the related concepts such as cocartesian-cartesian operad pairs, algebra pairs and Day convolutions. Using the Day convolution for operad pairs, we can define distributive functors between a symmetric bimonoidal category and a cartesian monoidal category. As an application, we show that for a finite group $G$, the operadic distributive functors from $\Span (\mathsf{Fin}_G)$ to any cartesian monoidal category coincides with the distributive functors defined by the distributive diagram. Therefore, we show that the connective $G$-equivariant normed-ring spectra are equivalently the spectral $G$-Tambara functors.
\end{abstract}
\maketitle

\markright{Infinity operad pairs, higher distributivity and spectral Tambara functors}

\tableofcontents

\section{Introduction}
The theory of operads brings a framework to study different algebraic structures and their modules in a combinatorial way. It was first introduced by \cite{may1972geometry} in a topological context, and was studied thoroughly afterwards by \cite{kriz1995operads}, \cite{Lein}, \cite{MSS}, \cite{Fresse0} and more. Informally, an operad consists of a sequence of objects together with compatible actions, such that the unitary, associative and equivariant laws are satisfied. The modern homotopical story of operads is based on several models of $\infty$-operads. A widely accepted model is given by \cite{lurie2017higher}, while other models include the dendroidal set model in \cite{cisinski2011dendroidal} and Barwick's complete Segal operads in \cite{barwick2018operator}. All these are compared and shown equivalent in \cite{chu2018two} and \cite{hinich2024equivalence}.

On the contrary, the story of operad pairs was fairly later. It was first defined by \cite{may2009construction} using the term an operad action on another operad. Informally, an operad pair consists of two operads and their interplay, mimicking the interaction between addition and multiplication. The operad pairs are useful especially when we want to work with two different operators simultaneously. In \cite{may2009construction}, May uses operad pairs to construct $\Einf$-ring spaces out of a bipermutative category.

This paper is devoted to the study of an $\infty$-categorical formulation of operad pairs. We construct a model of $\infty$-operad pairs where we can talk about the underlying operads and moperads. Our construction of $\infty$-operad pairs is based on Lurie's model of $\infty$-operads. The difficulty is to formulate the data of the distributivity law of operad pairs. In order to solve this difficulty, we first construct a small category $\Fint$, whose position is similar to that of $\Finp$ in mono-operadic world. We then define an $\infty$-operad pair as a category $\opp$ sliced over $\Fint$ such that certain conditions are satisfied. The category $\Fint$ itself does not become an $\infty$-operad. However, it is equipped with two functors $i_0,i_1\colon \Finp\rightrightarrows\Fint$, so we can define  he underlying operads and moperads through pullbacks and denote them by $\mathcal O_0^\otimes$ and $\mathcal O_1^{\otimes}$ respectively.

For mono-operads, certain mapping spaces of Lurie's $\infty$-operad together define a classical operad when the $\infty$-operad has only one color. For our $\infty$-operad pairs, corresponding to the mono-operad case, we can show that they are a reasonable generalization of classical operad in multi-colored case. More explicitly, every active map $\beta_m\colon \langle m\rangle \to \langle 1\rangle \in \Finp$ corresponds to two maps $\beta_m^+,\beta_m^\times \in \Fint$. When the $\infty$-operad pair has only one color $\mathbbl 1$, we can define
\[
\mathscr C (n):= \Map^{\beta_n^+}(\mathbbl 1^{\oplus n},\mathbbl 1) \quad \mathscr J(n):= \Map^{\beta_n^+}(\mathbbl 1^{\otimes n},\mathbbl 1)
\] as the lifts of mapping spaces. This gives two classical operads which are further identified as a classical operad pair.

\begin{mainthm}[\autoref{comparisontheorem}]\label{mainA}
    Let $p\colon\opp \to \mathrm{N}\Fint$ be a unital operad pair with a single color, then $\mathscr C(n)$ and $\mathscr J(n)$ consist of a classical operad pair in the 1-category of spaces.
\end{mainthm}


An $\infty$-operad pair together with an additional cocartesian lifting condition gives rise to the definition of a symmetric bimonoidal category or symmetric rig category. Classically, a symmetric bimonoidal category contains two symmetric monoidal structures and a distributive law connecting them. An important example is \emph{distributive} categories, which are categories such that coproducts distribute among products. In this case, there is a universal symmetric bimonoidal structure over it, which we refer to as the \emph{cocartesian-cartesian bimonoidal structure}.

\begin{mainthm}[\autoref{cct}, \autoref{undercct}]
    For any $\infty$-category $\mathcal C$, there is a functor $q' \colon \mathcal C^{\sqcup,\times}\to \mathrm{N}\Fint$ such that it is a symmetric bimonoidal category if and only if $\mathcal C$ is distributive. In this case, the underlying operad and moperad are respectively cartesian and cocartesian operads defined by Lurie. 
\end{mainthm}

For any $\infty$-operad pair $p\colon \opp \to \mathrm{N}\Fint$, we can talk about the category of $\mathcal C$-algebra over it, where $\mathcal C$ could be either an operad or an operad pair. We will concentrate to the case where the receptacle operad pair is a symmetric bimonoidal category. We denote
\[
\mathrm{CRig} (\mathcal P):= \mathrm{Alg}_{\mathcal C} (\mathcal P)
\] when $\mathcal C$ is the commutative bimonoidal operad pair $\mathsf{Comm}^{\oplus,\otimes}$. This interacts with the case where $\mathcal C$ is the commutative operad in the following way:

\begin{mainthm}[\autoref{calgtheorem}]\label{mainC}
    Let $p\colon \opp \to \mathrm{N}\Fint$ be an $\infty$-operad pair. There exists an $\infty$-operad $p'\colon \mathrm{CAlg}_{\oplus} (\mathcal P)^{\otimes} \to \mathrm{N}\Finp$ such that it is a symmetric monoidal category with underlying category identified with $\mathrm{CAlg}_{\oplus} (\mathcal P)$, whenever $\mathcal P$ is symmetric bimonoidal. Moreover, under this assumption we have
    \[
    \CAlg(\CAlg_{\oplus} (\mathcal P)) \simeq \mathrm{CRig}(\mathcal P)
    \]
\end{mainthm}

One property of great interest of operad pairs is the Day convolution. In \cite{lurie2017higher}, using the norm construction, Lurie defines a Day convolution symmetric monoidal category over the functor category between two symmetric monoidal categories. It is generally not true that two Day convolution monoidal structures together define a symmetric bimonoidal category since the failure of the distributivity. Instead, if we assume the source category is semi-additive, which will guarantee certain cofinality conditions, we can expect these two Day convolution monoidal structure interact with each other, such that one monoidal structure restrict to commutative algebra objects with respect to another.

\begin{mainthm}[\autoref{twotwodayconvolution}]\label{mainD}
    Let $p\colon \opp \to \Fint$ and $q\colon \mathcal Q^{\oplus,\otimes} \to \Fint$ be two symmetric bimonoidal categories such that the Day convolution between the underlying operads and moperads are defined. Suppose that $\mathcal P$ is semi-additive and the underlying operad of $\opp$ is the cartesian symmetric monoidal operad, then the Day convolution symmetric monoidal structure restricts to a symmetric monoidal structure
    \[
    p\colon (\mathrm{Alg}_{\mathcal P}\mathcal Q)^{\otimes} \to \mathrm{N}\Finp
    \] with the underlying category identified with $\mathrm{Alg}_{\mathcal P}\mathcal Q \simeq \Fun^{\oplus\text{-}\mathrm{lax}} (\mathcal P,\mathcal Q)$
\end{mainthm}

The proof of \autoref{mainC} and \autoref{mainD} relies on concrete construction of new operads. These operads are constructed by reshaping the category $\Fint$ such that there exists some proper maps to $\Finp$, allowing us to define the desired operads.

One distinctive feature of Day convolution from a symmetric bimonoidal category is the distributivity. Even though the target category is merely a cartesian monoidal category instead of a bimonoidal category, we can also use Day convolution to define distributive functor category $\Fun^{\mathrm{dist}}(\mathcal P,\mathcal C) \subset \Fun(\mathcal P,\mathcal C)$ which roughly consists of those functors that are symmetric monoidal for the first operation, lax monoidal for the second operation and satisfy a specific distributive law. With these in hand, let us now turn to the story of spans and bispans. 

The $\infty$-category of spans have several names in literature. It was introduced by \cite{barwick2013q} and studied further in \cite{barwick2017spectral}, where the terminology \emph{effective Burnside category} was used. In recent work of \cite{heyer20246}, the author used the term \emph{category of correspondence} and studied the symmetric monoidal structure  and self-enrichment over it. The term \emph{category of span} was used in \cite{haugseng2018iterated} and the following work \cite{elmanto2023distributivity}, \cite{cnossen2024normed} and \cite{cnossen2025homotopical}.

The definition of a suitable $\infty$-category of bispans is more difficult. In \cite{elmanto2023distributivity}, the author defined $L$-\emph{distributive functors} from the $(\infty,2)$-category of spans to an $(\infty,2)$-category. Thus the $(\infty,2)$-category of bispans is born with a universal property:

\[
\mathrm{FUN} (\mathrm{BISPAN}_{F,L} (\mathcal C), \mathcal X) \simeq \mathrm{FUN}_{L\text{-}\mathrm{dist}} (\Span_{F}(\mathcal C) ,\mathcal X)
\] for any span triple $(\mathcal C,\mathcal C_F,\mathcal C_L)$.  The results can be formulated into $(\infty,1)$-categorical settings and we can formally define $(\infty,1)$-category of bispans, which is born with a similar universal property.

This definition of distributivity is somehow intriguing. It involves a complicated definition of distributive diagram and certain Beck-Chevalley conditions, making the distributivity a little hard to handle with. Our next result is to unify the distributivity in this span context and in the operadic context. We will define a symmetric bimonoidal structure over a span category with the help of the cartesian-cocartesian bimonoidal formalism. For any cartesian monoidal category $\mathcal D$, we will compare the (operadic) distributive functors and $L$-distributive functors from $\Span (\mathsf{Fin})$ to $\mathcal C$, showing the definition of $L$-distributive functors is totally operadic in this case.

\begin{mainthm}\label{mainE}
    Let $\mathcal D$ be a presentable $\infty$-category that admits finite products. Then there exists an equivalence
    \[
    \Fun_{L\text{-}\mathrm{dist}}^\times (\Span (\mathsf{Fin}),\mathcal D) \simeq \Fun^{\mathrm{dist}}(\Span (\mathsf {Fin}),\mathcal D) 
    \]
\end{mainthm}

The \autoref{mainE} and its equivariant variants can be used to prove various results about commutative rings in non-equivariant and equivariant settings. Let $\mathcal D$ be any cartesian monoidal category. It is well-known that we have an equivalence
\[
\Fun^\times (\Span(\mathsf{Fin}), \mathcal D) \simeq \mathrm{CMon} (\mathcal D)
\] As a corollary of \autoref{mainE}, we can give a new proof of \cite[Theorem A]{cnossen2025homotopical}.

\begin{mainthm}[\autoref{noneqcor}]
    Let $\mathcal D$ be a cartesian closed presentable category. Then there is an equivalence
    \[
    \mathrm{CRig}(\mathcal D) \simeq \Fun^{\times}(\mathrm{Bispan}(\mathsf{Fin}),\mathcal D)
    \]
\end{mainthm}

Let $G$ be a finite set. In \cite{GM2} shows that the category of genuine $G$-spectra is equivalent to the category of \emph{spectral Mackey functors} $\Fun^{\times}(\Span (\mathsf{Fin}_G), \Sp)$, while a proposition of \cite{nardin2016parametrized} shows that connective $G$-spectra are equivalent to those product-preserving functors from $\Span(\mathsf{Fin}_G)$ to $\mathcal S_{\mathrm{gp}}$. Using \autoref{mainE}, we can show the following theorem, which allows us to view connective $G$-equivariant commutative ring spectra as connective \emph{spectral Tambara functors}.

\begin{mainthm}[\autoref{stf}]
     Let $G$ be a finite set. Then we have an equivalence:
    \[
    \mathrm{NAlg}(\mathrm{Sp}_G^{\geq 0}) \simeq \Fun^{\times}_{\mathrm{gp}} (\mathrm{Bispan} (\mathsf{Fin}_G), \mathcal S)
    \]
\end{mainthm}

Our proof of these two theorems differs from that of \cite{cnossen2025homotopical} and \cite{cnossen2024normed}. In \cite{cnossen2025homotopical}, the authors use a localization functor to the category of bispans of finite sets. While our proof relies on an inspection of distributivity in the different context: the bimonoidal one and the span category one respectively. Moreover, our proof is purely over operadic level, without any specific statement regarding the norm structure on both sides.

The contents of this paper are divided as follows. We will begin with the definition of $\infty$-operad pairs and related concepts in \autoref{oppsection}. In \autoref{comparesection}, we compare our definition with the classical one, proving the \autoref{mainA}. In \autoref{ccsection}, \autoref{algsection} and \autoref{daysection}, we will study the further properties of $\infty$-operad pairs, including the cocartesian-cartesian formalism, algebras and Day convolutions. In this part we will also give the first definition of distributive functor. We will recollect the properties of spans and bispans in \autoref{spansection}, offering another definition of distributivity. In \autoref{distcomparesection}, we will finally show these two distributivities coincide and prove \autoref{mainE}. 

\par
\quad
\par

\noindent\textbf{Acknowledgment} This paper was partially motivated by another project that the author is working on. The author would like to thanks to Zhiyou Wu for bringing up that project. The author also would like to thanks to Foling Zou for helpful conversations.

\par
\quad
\par

\noindent\textbf{Conventions} We fix the following conventions:

\begin{enumerate}[(i)]
    \item We say category for $\infty$-category.
    \item Let $\mathcal C$ and $\mathcal D$ be two symmetric monoidal categories, we let $\Fun(\mathcal C, \mathcal D)$ denote the $(\infty,1)$-category of functors between $\mathcal C$ and $\mathcal D$. Let $\Fun^\times(\mathcal C,\mathcal D)\subset\Fun(\mathcal C,\mathcal D)$ denote the full subcategory generated by those product preserving functors.
    \item We use $\mathsf{Fin}$ to denote the category of non-based finite sets. The objects are denoted by $\mathbf{n}$ for $n\in \mathbb N$.
    \item We use $\Finp$ to denote the category of based finite sets. The objects are denoted by $\langle n \rangle$ for $n \in \mathbb N$. We denote $\langle n \rangle^\circ\subset \langle n \rangle$ for the subset cutting out the base point $0$. We say a map $f\colon \langle m \rangle\in \Finp$ is inert if $f^{-1}\{i\}$ has exactly one element for each $i\in \langle n \rangle^{\circ}$.
    \item Let $G$ be a finite group, we use $\mathsf{Fin}_G$ to denote the category of finite $G$-sets.
    \item We adopt Lurie's model for $\infty$-operads.
    \begin{defn}[\cite{lurie2017higher}]\label{lurieoperad}
    An $\infty$-operad is a functor $\mathcal{O}^\otimes \to \mathrm N\Finp$ that satisfies the followings:
    \begin{enumerate}
        \item For every inert morphism $\alpha \colon \langle n \rangle \to \langle m \rangle$ in $\Finp$ and every $C \in \mathcal O^\otimes_{\langle n \rangle} $, there is a cocartesian lift $f_\alpha \colon C \to C'$ of $\alpha$. This determines a functor $\mathcal O ^\otimes_{\langle n \rangle} \to \mathcal O ^\otimes_{\langle m \rangle}$.

        \item Let $\rho_i \colon \langle m \rangle \to \langle1\rangle$ denote the inert morphism that sends $i$ to $1$. For every $f \colon \langle n \rangle \to \langle m \rangle$ in $\Finp$, the composition with $\rho_i$ induces an equivalence
        \[
        \Map^f (C, C') \simeq \prod_{i=1}^m\Map^{\rho_i f} (C, C_i')
        \] once we choose $C' \to C_i$ lifting $\rho_i$.

        \item For every collection $C_1, \cdots , C_m \in \mathcal O^\otimes _{\langle 1 \rangle }$, there exists $C \in \mathcal O^\otimes _{\langle m \rangle }$ and morphisms $C \to C_i$ lifting $\rho_i$.
    \end{enumerate}
\end{defn}
\end{enumerate}

\section{Infinity operad pairs}\label{oppsection}

An operad pair should consist of two operators and one distributive law between these two operators. In order to define an $\infty$-operad pair, we need to make sure several properties we expect it to have. First, there should be an "initial" operad pair over which we can define algebra pairs. Moreover, these algebras should exactly be identified with bialgebras in symmetric bimonoidal categories in the classical sense. Secondly, there should be two natural forgetful maps from an operad pair to an operad, corresponding to the underlying operad and moperad of an operad pair. Finally, it should be identified with the classical operad pairs defined by \cite{may2009construction} after we restrict to the "one colored" case.

A priori, we may naively think of an operad pair as a mono-operad in the sense of \cite{lurie2017higher} over a universal operad $\Lambda$. The universal operad $\Lambda$ has the same objects as $\Finp$. It serves as a package, containing the data of two operators $\oplus,\otimes$ and a distributive law connecting them. These data in $\Lambda$ are reflected by the active maps $\beta^+,\beta^\times\colon \langle n\rangle \to \langle 1 \rangle$ and distributive maps $d\colon \langle 2ab\rangle \to \langle a+b\rangle$  with certain compatibilities (see \autoref{fin2} for more details). However, there are several difficulties that keep us from realizing this intuition. Firstly, it is difficult to describe such a category $\Lambda$ as an operad. For example, it is hard to find the image in $\Finp$ of the distributive map $d \in \Lambda$. Moreover, the bimonoidal categories are no longer $\Lambda$-monoidal categories. For example, in bimonoidal categories we expect that the distributive law is only defined for $d\colon (A \otimes B ) \oplus (A \otimes C) \to A \otimes (B\oplus C)$, but the cocartesian fibration will lift $d $ starting from any objects that have the form $(A \otimes B ) \oplus (A' \otimes C)$.

For this reason, we will take another approach: We no longer view an operad pair as a mono-operad. Instead we develop a parallelized theory to describe operad pairs.

\begin{notn}\label{Qanddistequivalence}
    We denote $\Seq(\mathbb N):= \coprod_{k\in\mathbb N} \mathbb N^k$ to be the set of finite sequences of natural numbers. An element $Q \in \Seq (\mathbb N)$ is a set of the form $Q = (k, j_1,\cdots ,j_k)$ with $k\geq 0$, $j_1,\cdots,j_k \geq 0$. We denote 
    \[\underline{Q}:=\{(q_1,\cdots,q_k)|1\leq q_i\leq j_i \text{ for every }0\leq i \leq k\}
    \]to be the underlying set of tuples associated with $Q$ and denote $$|Q|:=|\underline{Q}|=j_1\cdots j_k$$. For $Q_1=(k,j_1,\cdots,j_k)$, $Q_2=(l,i_1,\cdots,i_l)$, we denote $$Q_1\star Q_2 := (k+l, j_1,\cdots,j_k,i_1,\cdots,i_l)$$ to be their joint. Finally, suppose we have $Q=(k,j_1,\cdots,j_k)$, $Q'=(k,i_1,\cdots, i_k)$ and a sequence of surjections $f_r\colon \mathbf{j_r} \to \mathbf{i_r}$ for $1 \leq r \leq k$. Then this sequence of surjections induces a map between underlying sets
    \[
    f\colon \underline{Q} \to \underline{Q'} \quad f_1^{-1}(q_1)\times \cdots \times f_k^{-1}(q_k) \mapsto (q_1,\cdots,q_k)
    \] For every $\alpha = (q_1,\cdots,q_k) \in \underline{Q'}$, we define the fiber of $\{f_i\}$ as $$Q_\alpha := (k, |f_1^{-1}(q_1)|,\cdots, |f^{-1}_k (q_k)|) \in \Seq(\mathbb N)$$. 
\end{notn}

In order to consider the equivariant transformations of operad pairs, we will also need the following notation from \cite{Mayinprogress}.

\begin{notn}\label{twonotnforeq}
    Let $j_1,\cdots,j_k$ be natural numbers. Let $\sigma \in \Sigma_k$ and $\sigma_i \in \Sigma_{j_i}$ be permutations. We can define $\sigma\langle j_1,\cdots, j_k\rangle$ and $\sigma_1\otimes\cdots\otimes \sigma_k$ in $\Sigma_{j_1\cdots j_k}$ through following diagrams.

    \[\begin{tikzcd}[nodes={scale=0.8}]
	{j_1\times\cdots\times j_k} & {j_{\sigma(1)}\times\cdots \times j_{\sigma(k)}} && {j_1\times\cdots \times j_k} & {j_1 \times \cdots \times j_k} \\
	{j_1\cdots j_k} & {j_1\cdots j_k} && {j_1\cdots j_k} & {j_1\cdots j_k}
	\arrow["\sigma", from=1-1, to=1-2]
	\arrow["\cong"', from=1-1, to=2-1]
	\arrow["\cong", from=1-2, to=2-2]
	\arrow["{\sigma_1 \times\cdots \times \sigma_k}", from=1-4, to=1-5]
	\arrow["\cong"', from=1-4, to=2-4]
	\arrow["\cong", from=1-5, to=2-5]
	\arrow["{\sigma\langle j_1,\cdots, j_k\rangle}"', from=2-1, to=2-2]
	\arrow["{\sigma_1\otimes\cdots \otimes \sigma_k}"', from=2-4, to=2-5]
\end{tikzcd}\]
\end{notn}

\begin{defn}\label{fin2}
    We define a small category $\Fint$ as follows:
    \begin{enumerate}[(1)]
        \item The objects of $\Fint$ are pairs $(\langle n \rangle,\gamma)$, where we have $\langle n \rangle \in \Finp$ and $\gamma$ is a binary number between $0$ and $2^{n-1}-1$. We can view $\langle 0\rangle= \emptyset$ as the base point. For any object $(\langle n \rangle, \gamma)$, we call $\langle n \rangle$ the \emph{length} of $(\langle n \rangle, \gamma)$. Intuitively, we can think of $\gamma $ as a combination of two different operators of the given length.
        \item There are two monoid structures on the underlying set of $\Fint$ which we denote by
        \[
        \boxplus \colon \big((\langle n \rangle, \gamma_1), (\langle m \rangle , \gamma_2) \big) \mapsto (\langle m+n\rangle, \gamma_1 0 \gamma_2 ) \]
        \[\boxtimes \colon \big((\langle n \rangle, \gamma_1), (\langle m \rangle , \gamma_2) \big) \mapsto (\langle m+n\rangle, \gamma_1 1 \gamma_2 ).
        \] Namely the operators by piecing two binary representation together var $0$ and $1$ respectively. The underlying subset of $\Fint$ without the base point is generated by $\langle 1 \rangle =(\langle 1 \rangle , 0)$ under these two monoid structures.

        \item The morphisms of $\Fint$ are generated by the followings:
        \begin{enumerate}[(a)]
            \item For each $\alpha \colon \langle n \rangle \to \langle m \rangle \in \Finp$, there is an associated morphism $$\alpha^+ \colon (\langle  n\rangle, 0) \to  (\langle  m\rangle, 0)$$.
            \item For each $\alpha \colon \langle n \rangle \to \langle m \rangle \in \Finp$, there is an associated morphism $$\alpha^\times \colon (\langle  n\rangle, 2^{n-1}-1) \to  (\langle  m\rangle, 2^{n-1}-1)$$.
            \item Let $A_1 ,\cdots ,A_k$ be a sequence of objects. Then any $\sigma \in \Sigma_k$ gives rise to two permutation isomorphisms:
            \[
            \sigma^+ \colon A_{\sigma^{-1}(1)}\boxplus \cdots \boxplus A_{\sigma^{-1}(k)} \to A_1 \boxplus \cdots \boxplus A_k
            \]
            \[
            \sigma^\times \colon A_{\sigma^{-1}(1)}\boxtimes \cdots \boxtimes A_{\sigma^{-1}(k)} \to A_1 \boxtimes \cdots \boxtimes A_k
            \]When $A_1 =\cdots = A_k = \langle 1 \rangle$, this coincides with the definition above for $\sigma \colon \langle k\rangle \to \langle k \rangle$ in $\Finp$.
            
            \item For every $Q \in \Seq(\mathbb N)$, there is a map which we refer to as the \emph{distribute map}
            
            \[
            d_Q \colon \bigboxplus_{\alpha \in \underline{Q}} (A ^{1,q^\alpha_1} \boxtimes \cdots \boxtimes A^{k,q^\alpha_k}) \to \bigboxtimes_{1 \leq r\leq k} (A^{r,1}\boxplus \cdots \boxplus A^{r,j_r})
            \] where each $A^{m, q_m}$ is an object for $1 \leq m \leq k$, $1 \leq q_m \leq j_m$ and $\alpha $ runs over all tuples in $\underline{Q}$.
        \end{enumerate} The morphisms of $\Fint$ are generated by the above morphisms under composition, $\boxplus$ and $\boxtimes$. These morphisms must satisfy the following requirements: 
        \begin{enumerate}[(i)]
            \item (Unitary) For $k=1$, we require $d_Q= \mathrm{id}\colon A^1 \boxplus \cdots \boxplus A^j \to A^1 \boxplus \cdots \boxplus A^j$.
            \item For $\alpha_1, \alpha_2 \in \Finp$, we have  $$\alpha_1^+ \boxplus \alpha_2^+ = (\alpha_1+\alpha_2)^+ ,\quad\alpha_1 ^\times \boxtimes \alpha_2^\times = (\alpha_1+\alpha_2)^\times$$ where $+$ denotes the additive monoid structure on $\Finp$. For $f_i\colon A_i \to B_i$, $g_i \colon B_i \to C_i$ with $i=1,2$, we have $$ g_1 f_1\boxplus g_2 f_2 = (g_1 \boxplus g_2) (f_1 \boxplus f_2), \quad g_1 f_1\boxtimes g_2 f_2 = (g_1 \boxtimes g_2) (f_1 \boxtimes f_2) $$ We require $\mathrm{id}:=\mathrm{id}^+ \simeq \mathrm{id}^\times \colon \langle 1 \rangle \to \langle 1 \rangle$. Moreover, let $*\colon \langle 1\rangle \to \langle0 \rangle \in \Finp$, we require $*^+ = *^\times \colon \langle 1 \rangle \to \langle 0 \rangle$ and denote this map by $*$.
            \item (Functoriality) The map $d_Q$ is functorial. That means, given a family of morphisms $\{ f^{m, q_m} \colon A^{m,q_m} \to B^{m, q_m}\}$, it fits in a commutative diagram:
\[\begin{tikzcd}[nodes={scale=0.8}]
	{\bigboxplus_{\alpha \in \underline{Q}} (A ^{1,q^\alpha_1} \boxtimes \cdots \boxtimes A^{k,q^\alpha_k})} & {\bigboxtimes_{1 \leq r\leq k} (A^{r,1}\boxplus \cdots \boxplus A^{r,j_r})} \\
	{\bigboxplus_{\alpha \in \underline{Q}} (B ^{1,q^\alpha_1} \boxtimes \cdots \boxtimes B^{k,q^\alpha_k})} & {\bigboxtimes_{1 \leq r\leq k} (B^{r,1}\boxplus \cdots \boxplus B^{r,j_r})}
	\arrow["{d_Q}", from=1-1, to=1-2]
	\arrow["{\bigboxplus_{\alpha \in \underline{Q}} (f ^{1,q^\alpha_1} \boxtimes \cdots \boxtimes f^{k,q^\alpha_k})}"', from=1-1, to=2-1]
	\arrow["{\boxtimes_{1 \leq r\leq k} (f^{r,1}\boxplus \cdots \boxplus f^{r,j_r})}", from=1-2, to=2-2]
	\arrow["{d_Q}", from=2-1, to=2-2]
\end{tikzcd}\]
\item (Equivariant) The map $d_Q$ respects $\sigma^+$ and $\sigma^\times$. That means, suppose that we have we are given $\sigma \in \Sigma_r$ and $\sigma_r \in \Sigma_{j_r}$ for $1\leq r \leq k$, then let $\sigma_1\otimes\cdots\otimes \sigma_r $ and $\hat{\sigma}:= \sigma \langle j_1 , \cdots , j_k \rangle $ be as in \autoref{twonotnforeq}, we have following commutative diagrams:

\[\begin{tikzcd}[nodes={scale=0.8}]
	{\bigboxplus_{\alpha \in \underline{Q}} (A ^{1,q^\alpha_1} \boxtimes \cdots \boxtimes A^{k,q^\alpha_k})} & {\bigboxtimes_{1 \leq r\leq k} (A^{r,1}\boxplus \cdots \boxplus A^{r,j_r})} \\
	{\bigboxplus_{\alpha \in \underline{Q}} (A ^{1, q^{\alpha'}_1} \boxtimes \cdots \boxtimes A^{k,q^{\alpha'}_k})} & {\bigboxtimes_{1 \leq r\leq k} (A^{r,\sigma_r(1)}\boxplus \cdots \boxplus A^{r,\sigma_r(j_r)})}
	\arrow["{d_Q}", from=1-1, to=1-2]
	\arrow["{\sigma_1 \otimes \cdots \otimes \sigma_r}"', from=1-1, to=2-1]
	\arrow["{\sigma_1 \boxtimes \cdots \boxtimes \sigma_r}", from=1-2, to=2-2]
	\arrow["{d_Q}", from=2-1, to=2-2]
\end{tikzcd}\]

\[\begin{tikzcd}[nodes={scale=0.8}]
	{\bigboxplus_{\alpha \in Q} (A ^{1,q^\alpha_1} \boxtimes \cdots \boxtimes A^{k,q^\alpha_k})} & {\bigboxtimes_{1 \leq r\leq k} (A^{r,1}\boxplus \cdots \boxplus A^{r,j_r})} \\
	{\bigboxplus_{\alpha \in Q} (A ^{\sigma (1), q^{\alpha}_{\sigma(1)}} \boxtimes \cdots \boxtimes A^{\sigma(k),q^{\alpha}_{\sigma(k)}})} & {\bigboxtimes_{1 \leq r\leq k} (A^{\sigma(r),1}\boxplus \cdots \boxplus A^{\sigma(r),j_{\sigma(r)}})}
	\arrow["{d_Q}", from=1-1, to=1-2]
	\arrow["{\tilde \sigma}"', from=1-1, to=2-1]
	\arrow["{\sigma}", from=1-2, to=2-2]
	\arrow["{d_Q}", from=2-1, to=2-2]
\end{tikzcd}\]
\item (The first composition law) Let $Q_1, \cdots,Q_s \in\Seq(\mathbb N)$ with $$Q_i=(k_i, j_{i,1},\cdots,j_{i,k_i}):=(k_i,j_{k_1+\cdots+k_{i-1}+1},\cdots, j_{k_1+\cdots k_i})$$ We let $Q:= Q_1\star\cdots\star Q_s$ denote their joint. Let $Q':= (k, |Q_1|,\cdots ,|Q_k|)$. Then we have a following commutative diagram: \[\begin{tikzcd}[nodes={scale=0.8}, sep=small]
	{\bigboxplus_{\alpha \in \underline{Q}} (A ^{1,q^\alpha_1} \boxtimes \cdots \boxtimes A^{\Sigma k_i,q^\alpha_{\Sigma k_i}})} & {\bigboxtimes_{1 \leq r\leq \Sigma k_i} (A^{r,1}\boxplus \cdots \boxplus A^{r,j_r})} \\
	{\bigboxtimes_{1\leq i \leq s}(\bigboxplus_{\alpha \in \underline{Q_i}} (A ^{j_{i,1},q^\alpha_{j_{i,1}}} \boxtimes \cdots \boxtimes A^{j_{i,k_i},q^\alpha_{j_{i,k_i}}})) }
	\arrow["{d_{Q}}", from=1-1, to=1-2]
	\arrow["{d_{Q'}}"', from=1-1, to=2-1]
	\arrow["{\bigboxtimes d_{Q_i}}"', from=2-1, to=1-2]
\end{tikzcd}\]

\item (The second composition law) Let $Q=(k,j_1,\cdots,j_k)$, $Q'=(k,i_1,\cdots, i_k)$ and $f_r\colon \mathbf{j_r} \to \mathbf{i_r}$ for $1 \leq r \leq k$ be a sequence of surjections. Then we have a following commutative diagram:


\[\begin{tikzcd}[nodes={scale=0.8}]
	{\bigboxplus_{\alpha \in \underline{Q}} (A ^{1,q^\alpha_1} \boxtimes \cdots \boxtimes A^{k,q^\alpha_k})} & {\bigboxtimes_{1 \leq r\leq k} (A^{r,1}\boxplus \cdots \boxplus A^{r,j_r})} \\
	{\bigboxplus_{\alpha' \in \underline{Q'}} (\bigboxplus_{q_1^\alpha \in f^{-1}(q_1^{\alpha'})} A ^{1,q^\alpha_1} \boxtimes \cdots \boxtimes \bigboxplus_{q_k^\alpha \in f^{-1}(q_k^{\alpha'})} A ^{k,q^\alpha_k})}
	\arrow["{d_Q}", from=1-1, to=1-2]
	\arrow["{\bigboxplus_{\alpha'\in\underline{Q'}}d_{Q_{\alpha '}}}"', from=1-1, to=2-1]
	\arrow["{d_{Q'}}"', from=2-1, to=1-2]
\end{tikzcd}\]


\item (Distributivity) Let $Q=(k,j_1,\cdots,j_k)\in \Seq(\mathbb N)$ and all $A^{r,q}=\langle 1 \rangle$. Let $\beta_i \colon \langle i \rangle \to \langle 1 \rangle $ denote the morphisms in $\Finp$ that sends all non-zero elements to $1$. Then we require the following \emph{distributive square} to commute:

            \[\begin{tikzcd}
	{\bigboxplus_{|Q|}\bigboxtimes_k \langle 1 \rangle} & {\bigboxtimes_{1\leq r \leq k} \bigboxplus_{j_r} \langle 1 \rangle} \\
	& {\bigboxtimes_{1\leq r \leq k}\langle 1 \rangle} \\
	{\bigboxplus_{|Q|}\langle 1 \rangle} & {(\langle 1\rangle, 0)}
	\arrow["{d_{Q}}", from=1-1, to=1-2]
	\arrow["{\bigboxplus_{|Q|} \beta_k^\times}"', from=1-1, to=3-1]
	\arrow["{\bigboxtimes_{1\leq r \leq k } \beta^+_{j_r}}", from=1-2, to=2-2]
	\arrow["{\beta_k^\times}", from=2-2, to=3-2]
	\arrow["{\beta_{|Q|}^+}", from=3-1, to=3-2]
\end{tikzcd}\]
\item (Restriction) We use $A^Q_\times$ to denote $\bigboxplus_{\alpha \in \underline{Q}} (A ^{1,q^\alpha_1} \boxtimes \cdots \boxtimes A^{k,q^\alpha_k})$ and use $A^Q_+$ to denote $\bigboxtimes_{1 \leq r\leq k} (A^{r,1}\boxplus \cdots \boxplus A^{r,j_r})$. There are $n(m)= \prod_{r\neq m}j_r$ copies of $A^{m,q_m}$ in $A^Q_\times$ for each $m$ and $q_m$. We use $\rho^{m,q_m,l}_\times\colon A^Q_\times \to A^{m,q_m}$ to denote the projection on these copies for $1\leq l \leq n(m)$. That means, the morphism of the form $*\bullet \cdots \bullet \mathrm{id}\bullet \cdots \bullet *$ where the identity is taken over the $i$-th copy of $A^{m,q_m}$. Similarly, we use $\rho^{m,q_m}_+\colon A^Q_+ \to A^{m,q_m}$ to denote the projection on $A^{m,q_m}$. Then we require the following diagram to commute for any $1 \leq l \leq n(m)$:
\[\begin{tikzcd}
	{A^{Q}_\times} & {A^Q_{+}} \\
	& {A^{m,q_m}}
	\arrow["{d_Q}", from=1-1, to=1-2]
	\arrow["{\rho^{m,q_m,l}_\times}"', from=1-1, to=2-2]
	\arrow["{\rho_+^{m,q_m}}", from=1-2, to=2-2]
\end{tikzcd}\]
        \end{enumerate}
    \end{enumerate}
\end{defn}


\begin{rem}
    Here we use the notation $\Fint$ to stands for ``multi-finite''. Under the operations $\boxplus$ and $\boxtimes$, $\Fint$ becomes a symmetric bimonoidal category in classical sense with both the additive and multiplicative units given by $\langle 0 \rangle$. 
\end{rem}

\begin{exmp}\label{exampleofmorphisms}
    We list some morphisms in $\Fint$ and their properties.
    \begin{enumerate}[(i)]
        \item Recall that a morphism $f$ in $\Finp$ is called \emph{active} if $f^{-1}\{*\} = \{ *\}$. Such examples include $\beta_k$ in the above definition. In $\Fint$, we say that a morphism is active if it can be written as $f= f_1 \bullet \cdots \bullet f_n$ where each $f_i $ is of the form $f'^\times_i$ or $f'^+_i$ for an active $f'_i$ and each bullet represents a $\boxplus$ or a $\boxtimes$.

        \item In $\Finp$, there are some "elementary" morphisms. That is, there are morphisms $\delta^m_{i}$ and $\eta^m_i$. Here $\delta^m_i$ are monotone injections defined as
        \[
        \delta^m_i \colon \langle m \rangle \to \langle m+1 \rangle \quad \{1, \cdots, m\} \mapsto \{ 1, \cdots , \hat{i} , \cdots m+1\}
         \] and $\eta^m_i$ are monotone surjections defined as 

        \[
        \eta^m_i \colon \langle m+1 \rangle \to \langle m \rangle \quad \eta^m_i\{ i\} = \eta^m_i \{i+1\}
        \] It is well-known that the morphisms in $\Finp$ are generated under composition by these elementary morphisms and permutations. As a result, the morphisms in $\Fint$ are generated by $(\delta^m_i)^+, (\delta^m_i)^\times, (\eta^m_i)^+, (\eta^m_i)^\times, \sigma^+,\sigma^\times$ and $d_Q$ under composition, $\boxplus$ and $\boxtimes$.

        \item We've define the projections in [\autoref{fin2}, (viii)]. For $1 \leq i \leq m$, we can define $\rho_i \colon (\langle m \rangle, \gamma  ) \to \langle 1 \rangle$ for any $\gamma$ that makes sense. To define this, we note that any object $(\langle m \rangle, \gamma  ) \in \Fint$ can be written as $\langle 1 \rangle \bullet \cdots \bullet \langle 1 \rangle$, where each $\bullet$ represents a $\boxplus$ or a $\boxtimes$. Then we define $\rho_i = * \bullet\cdots \bullet \mathrm{id} \bullet \cdots \bullet *$ where the identity is taken over the $i$-th term.
    \end{enumerate}
\end{exmp}

\begin{notn}\label{twoembed}
    We use $i_0$ (resp.$i_1$) to denote the functor $\Finp \to \Fint$ sending $\langle n \rangle$ to $(\langle n\rangle,0)$ (resp.$(\langle n\rangle,1)$) and $\alpha \colon \langle n \rangle \to \langle m \rangle $ to $\alpha^+ \colon (\langle n\rangle,0) \to (\langle m\rangle ,0)$ (resp.$\alpha^+ \colon (\langle n\rangle,2^{n-1}-1) \to (\langle m\rangle ,2^{m-1}-1)$).
\end{notn}

It is natural to think of $\Fint$ as a replacement for $\Finp$ in the bimonoidal world. This philosophy inspires the following definition of $\infty$-categorical version of operad pairs. Before that, we need some generalization of terminology.

\begin{defn}
    We say that a morphism $\alpha \in \Fint$ is \emph{inert} if it can be written as $\alpha = \alpha_1 \bullet \cdots \bullet \alpha_n$, where each bullet here represents a composition, a $\boxplus$ or a $\boxtimes$ and each $\alpha_k$ is an $f^+, f^\times, \sigma^+$ or $\sigma^\times$ as in [\autoref{fin2}, (a)-(c)]. where $\sigma$ is a permutation and $f$ is inert in the sense of \cite[Definition 2.1.1.8]{lurie2017higher}. That is,  $f^{-1}(i)$ has exactly one element for each non-zero $i$.
\end{defn}

\begin{defn}\label{operadpairs}
    An $\infty$-\emph{operad pair} is a functor $p\colon\mathcal P^{\oplus, \otimes} \to \mathrm N\Fint$ that satisfies the followings: 
    \begin{enumerate}
        \item For every inert morphism $\alpha\colon (\langle n \rangle , \gamma_1) \to (\langle m \rangle , \gamma_2)$ and every $C \in \opp_{(\langle n \rangle , \gamma_1)}$, there is $p$-cocartesian lift $C \to C'$ of $\alpha$. Especially, we have an induced functor $\alpha_! \colon \opp_{(\langle n \rangle , \gamma_1)} \to \opp_{(\langle m \rangle , \gamma_2)}$
        \item Let $\rho_i$ be defined as in \autoref{exampleofmorphisms}. For every morphism $f$ in $\Fint$, the composition with $\rho_i$ induces an equivalence 
        \[\Map^f (C, C') \simeq \prod_{i=1}^m\Map^{\rho_i f} (C, C_i')
        \] once we choose $C' \to C_i$ lifting $\rho_i$.
        \item For every collection $C_1, \cdots , C_m \in \mathcal P^{\oplus ,\otimes} _{\langle 1 \rangle } := \mathcal P^{\oplus ,\otimes} _{\langle 1 \rangle }$ and every $\gamma \leq 2^{m-1}-1$, there exists $C \in \mathcal \opp_{(\langle m \rangle , \gamma)}$ and morphisms $C \to C_i$ lifting $\rho_i$.
        \item Let notations be as in \autoref{fin2} and \autoref{exampleofmorphisms}. For every $C \in \opp_{A^Q_\times}$ such that $(\rho^{m,q_m,l_1}_\times)_{!}(C) \simeq (\rho^{m,q_m,l_2}_\times)_{!}(C)$ for any fixed $m,q_m$ and $0\leq l_1 \leq l_2 \leq n(m)$, there is a cocartesian lift $\tilde {d}_{Q} \colon C \to C'$ of $d_{Q}\colon A^Q_\times \to A^Q_+$ fitting into a lift of  the diagram in [\autoref{fin2},(viii)].
\[\begin{tikzcd}
	C & {C'} \\
	& {C^{m,q_m}}
	\arrow["{\tilde d_Q}", from=1-1, to=1-2]
	\arrow["{\tilde\rho^{m,q_m,l}_\times}"', from=1-1, to=2-2]
	\arrow["{\tilde\rho_+^{m,q_m}}", from=1-2, to=2-2]
\end{tikzcd}\]
        In particular, we have $(\rho_+^{m,q_m})_! (C')\simeq (\rho^{m,q_m,l}_{\times})_!(C) $ for all $1 \leq l \leq n(m)$.
    \end{enumerate}
\end{defn}

\begin{defn}\label{underlyingopandmop}
    Let $p\colon \opp \to \mathrm{N}\Fint$ be an $\infty$-operad pair and let $i_0$ and $i_1$ be defined as in \autoref{twoembed}. Then we define two categories $\mathcal{O}^\otimes_0$ and $\mathcal O^\otimes_1$ in ${\mathrm{Cat}_\infty}_{/\mathrm{N}\Finp}$ as pullbacks:
    \[\begin{tikzcd}
	{\mathcal O^\otimes_0} & {\mathcal P^{\oplus ,\otimes}} && {\mathcal O^\otimes_1} & {\mathcal P^{\oplus ,\otimes}} \\
	{\mathrm{N}\Finp} & {\mathrm{N}\Fint} && {\mathrm{N}\Finp} & {\mathrm{N}\Fint}
	\arrow[dashed, from=1-1, to=1-2]
	\arrow["{p_0}"', dashed, from=1-1, to=2-1]
	\arrow["p", from=1-2, to=2-2]
	\arrow[dashed, from=1-4, to=1-5]
	\arrow["{p_1}"', dashed, from=1-4, to=2-4]
	\arrow["p", from=1-5, to=2-5]
	\arrow["{i_0}", from=2-1, to=2-2]
	\arrow["\ulcorner"{anchor=center, pos=0.125, rotate=180}, draw=none, from=2-2, to=1-1]
	\arrow["{i_1}", from=2-4, to=2-5]
	\arrow["\ulcorner"{anchor=center, pos=0.125, rotate=180}, draw=none, from=2-5, to=1-4]
\end{tikzcd}\]We refer to them as the \emph{underlying operad} and \emph{underlying moperad} of $\opp$ respectively. 
\end{defn}

This definition makes sense. When restricting to the subcategory $\mathrm{N}\Finp$ via $i_0$ and $i_1$, the conditions (1)-(3) of \autoref{operadpairs} exactly coincide with the operad conditions. Therefore, we have the following:

\begin{prop}
    Let notations be as in \autoref{underlyingopandmop}. The functor $p_0$ and $p_1$ exhibit $\mathcal{O}^\otimes_0$ and $\mathcal O^\otimes_1$ as two operads.
\end{prop}

\begin{rem}
    The term ``moperad'' is borrowed from \cite{Mayinprogress}, which means operads in multiplicative context. However, it turns out that our resulting operad $\mathcal O^\otimes_1$ has nothing to do with May's moperads.
\end{rem}

We shall explain \autoref{operadpairs} more concretely. The first three conditions show that, similarly as in mono-operad case, we have identifications $\opp_{(\langle m \rangle , \gamma)} \simeq (\opp_{(\langle1\rangle, 0 )})^m$ and write objects in $\opp_{(\langle m \rangle , \gamma)}$ as a tuple $(C_1, \cdots , C_m)$. There are two operators, namely $\oplus$ and $\otimes$, corresponding to $\boxplus$ and $\boxtimes$ in the binary representation, such that for $(C_1, \cdots , C_m) \in \opp_{(\langle m \rangle , \gamma_1)}$ and $(C'_1, \cdots , C'_n) \in \opp_{(\langle n \rangle , \gamma_2)}$, we have

\[(C_1, \cdots , C_m) \oplus (C'_1, \cdots , C'_n) \in \opp_{(\langle m+n \rangle , \gamma_1\boxplus \gamma_2)}\] 
\[(C_1, \cdots , C_m) \otimes (C'_1, \cdots , C'_n) \in \opp_{(\langle m+n \rangle , \gamma_1\boxtimes \gamma_2)}\] One can think of $\gamma$ as a combination of two operators. Therefore, we identify objects in $\opp$ as $C_1 \odot C_2 \odot\cdots \odot C_m$, where each $\odot $ either represents a copy of $\oplus$ or $\otimes$. We will sometimes write $\mathcal{P}:= \opp_{\langle 1 \rangle} \simeq (\mathcal{O}^\otimes_0)_{\langle 1 \rangle} \simeq (\mathcal{O}^\otimes_1)_{\langle 1 \rangle}$ and refer to them as the \emph{underlying category} of $\opp$. 

\begin{exmp}\label{exampleofopp}
    We give some examples of operad pairs.
    \begin{enumerate}[(1)]
        \item We define $\Comm^{\oplus,\otimes}$ as $\Comm^{\oplus, \otimes}:=\mathrm{N}\Fint$ and call it the \emph{bicommutative operad pair}. The underlying operad and moperad of $\Comm^{\oplus, \otimes}$ are $\Comm^\otimes$.
        \item We define an operad pair $\Ass^{\oplus,\otimes}$ as follows. The objects of $\Ass^{\oplus,\otimes}$ are the same as that of $\mathrm{N}\Fint$, while the generators of morphisms are almost the same, except that we need to replace (b) by the following: 
        \begin{list}{(b)$'$}{\setlength{\leftmargin}{2em}}
        \item For each $(\alpha,\prec) \colon \langle n \rangle \to \langle m \rangle \in \Ass^\otimes$, there is an associated morphism $(\alpha,\prec)^\times \colon (\langle  n\rangle, 2^{n-1}-1) \to  (\langle  m\rangle, 2^{n-1}-1)$.
        \end{list}
        The morphisms in $\Ass^{\oplus,\otimes}$ are generated by  almost the same relations as in \autoref{fin2}, except we need to replace (vii) by the following since maps are equipped with different orders here:
\begin{list}{(vii)$'$}{\setlength{\leftmargin}{2em}}
 \item Let $Q=(k,j_1,\cdots,j_k)\in \Seq(\mathbb N)$ and all $A^{r,q}=\langle 1 \rangle$. Let $f\colon \langle k\rangle \to \langle 1 \rangle$ be any morphisms in $\Ass^\otimes$ that send every non-zero element to $1$. Then the following square commutes:

 \[\begin{tikzcd}
	{\bigboxplus_{|Q|}\bigboxtimes_k \langle 1 \rangle} & {\bigboxtimes_{1\leq r \leq k} \bigboxplus_{j_r} \langle 1 \rangle} \\
	& {\bigboxtimes_{1\leq r \leq k}\langle 1 \rangle} \\
	{\bigboxplus_{|Q|}\langle 1 \rangle} & {(\langle 1\rangle, 0)}
	\arrow["{d_{Q}}", from=1-1, to=1-2]
	\arrow["{\bigboxplus_{|Q|} f^\times}"', from=1-1, to=3-1]
	\arrow["{\bigboxtimes_{1\leq r \leq k } \beta^+_{j_r}}", from=1-2, to=2-2]
	\arrow["{f^\times}", from=2-2, to=3-2]
	\arrow["{\beta_{|Q|}^+}", from=3-1, to=3-2]
\end{tikzcd}\]
\end{list}
We take $p\colon \Ass^{\oplus, \otimes} \to \mathrm{N}\Fint$ be the projection that sends $f^\times$ to $u(f)^\times$, where $u (f) \in \Finp$ forgets the ordering on $f^{-1}(i)$. This gives $\Ass^{\oplus, \otimes}$ an operad pair structure and we call it the \emph{biassociative operad pair}. Moreover, the underlying operad and moperad of $\Ass^{\oplus, \otimes}$ are $\Comm^\otimes$ and $\Ass^{\otimes}$ respectively.

\item We can also define the \emph{lax operad pair} $\mathsf{Lax}^{\oplus,\otimes}$ by setting $\mathsf{Lax}^{\oplus,\otimes}$ as freely adding a morphism $u \colon (\langle1\rangle,0) \to \langle 1 \rangle$ to $\mathrm{N}\Fint$ and replacing the relation (vii) in \autoref{fin2} by the following:

 \[\begin{tikzcd}
	{\bigboxplus_{|Q|}\bigboxtimes_k \langle 1 \rangle} & {\bigboxtimes_{1\leq r \leq k} \bigboxplus_{j_r} \langle 1 \rangle} \\
	{\bigboxplus_{|Q|}\langle 1 \rangle}& {\bigboxtimes_{1\leq r \leq k}\langle 1 \rangle} \\
	{\langle 1 \rangle)} & {(\langle 1\rangle, 0)}
	\arrow["{d_{Q}}", from=1-1, to=1-2]
	\arrow["{\bigboxplus_{|Q|} \beta_k^\times}"', from=1-1, to=2-1]
	\arrow["{\bigboxtimes_{1\leq r \leq k } \beta^+_{j_r}}", from=1-2, to=2-2]
	\arrow["{\beta_k^\times}", from=2-2, to=3-2]
	\arrow["{\beta_{|Q|}^+}", from=2-1, to=3-1]
    \arrow["{u}", from=3-1, to=3-2]
\end{tikzcd}\] 
Similarly, by inverting the direction of $u$ in the above diagram, we can define the \emph{oplax operad pair} $\mathsf{Oplax}^{\oplus, \otimes}$. We take the structure functor $p\colon \mathsf{Lax}^{\oplus,\otimes} \to \mathrm{N}\Fint$ and $p' \colon \mathsf{Oplax}^{\oplus, \otimes} \to \mathrm{N}\Fint$ to be the projection that sends $u$ to $\mathrm{id}$. \autoref{laxmonoidal} shows that $\mathsf{Lax}^{\oplus,\otimes}$ and $\mathsf{Oplax}^{\oplus, \otimes}$ are operad pairs. The underlying operad and moperad are all $\Comm^\otimes$.
    \end{enumerate}
\end{exmp}



We are going to define the morphisms between operad pairs. Before that, we will need a general definition of inert and distributive morphisms.

\begin{defn}
    Let $p\colon \opp \to \mathrm{N}\Fint$ be an $\infty$-operad pair. We say that a morphism in $\opp$ is \emph{inert} if it is a cocartesian lift of some inert morphism in $\Fint$. We say that a morphism in $\opp$ is \emph{distributive} if it is a cocartesian lift of some $d_{Q}$ and is of the form $\tilde{d}_{Q}\colon C\to C'$ for some $C $ and $C'$, such that $$\rho^{m,q_m,l_1}_\times(C) \simeq \rho^{m,q_m, l_2}_\times(C) = \rho_+^{m,q_m} (C')$$ for any $0\leq l_1 \leq l_2 \leq n(m)$.
\end{defn}

\begin{defn}\label{mapsofpairs}
    Let $p\colon\opp \to \mathrm{N}\Fint$ and $p'\colon (\opp)' \to \mathrm{N}\Fint$ be two operad pairs. A map between them is a functor over $\Fint$ that preserves inert and distributive morphisms. That is, a functor $f\colon \opp \to (\opp)'$ such that the diagram
\[\begin{tikzcd}
	\opp && {(\opp)'} \\
	& {\mathrm{N}\Fint}
	\arrow["f", from=1-1, to=1-3]
	\arrow["p"', from=1-1, to=2-2]
	\arrow["{p'}", from=1-3, to=2-2]
\end{tikzcd}\]
    commutes and takes inert (or distributive) morphisms to inert (or distributive) morphisms.
\end{defn}

We will use $\mathrm{Alg}_{\mathcal P} (\mathcal P') \subset \Fun_{\mathrm N \Fint} ( \opp, (\opp)')$ to denote the full subcategory spanned by those operad pair maps. In mono-operad case, an important class of operad maps is those cocartesian fibrations of operads. These are those cocartesian fibrations that exhibit the source as an operad. We can use the cocartesian fibrations to define monoidal categories over an operad. In operad pair world, however, the definition becomes slightly different. We begin with the definition of cocartesian fibrations of $\infty$-operad pairs.

\begin{defn}
    Let $\opp$ be an $\infty$-operad pair. We say a cocartesian fibration $p\colon \mathcal{C}^{\oplus,\otimes} \to \opp$ is a cocartesian fibration of operad pairs if it exhibits $\mathcal{C}^{\oplus,\otimes} \to \opp$ as an $\infty$-operad pair via $\mathcal{C}^{\oplus,\otimes} \xrightarrow{p} \opp \to \mathrm{N}\Fint$.
\end{defn}


In \cite[Proposition 2.1.1.12]{lurie2017higher}, Lurie shows that the property of being a cocartesian fibration of operad pairs is characterized by the intrinsic property of $p$ itself. For operad pairs, we have a similar result.

\begin{prop}\label{settingmonoidal}
    Let $\opp$ be an $\infty$-operad pair and let $p\colon \mathcal C^{\oplus, \otimes} \to \opp$ be a functor. Consider the following conditions:

    \begin{enumerate}[(a)]
        \item $p\colon \mathcal{C}^{\oplus, \otimes} \to \opp$ is a cocartesian fibration of operad pairs.
        \item $p$ is a cocartesian fibration. For every object $(\langle m \rangle , \gamma) \in \Fint$ and every object $T \simeq T_1 \odot \cdots \odot T_m \in \mathcal{P}^{\oplus,\otimes}_{(\langle m \rangle,\gamma)}$, we have an equivalence
        \[
        \mathcal C^{\oplus, \otimes}_{T} \simeq \prod_{1 \leq i \leq m} \mathcal C ^{\oplus, \otimes}_{T_i}
        \] which is induced by inert maps $T \to T_i$.
        \item Consider the following diagrams consisting of the pullback squares:
        \[\begin{tikzcd}
	{\mathcal C^\otimes_0} & {\mathcal O^\otimes_0} & {\mathrm N \Fint} \\
	{\mathcal{C}^{\oplus, \otimes}} & \opp & {\mathrm N \Finp}
	\arrow["{p_0}", from=1-1, to=1-2]
	\arrow[from=1-1, to=2-1]
	\arrow["{q_0}", from=1-2, to=1-3]
	\arrow[from=1-2, to=2-2]
	\arrow["{i_0}", from=1-3, to=2-3]
	\arrow["p", from=2-1, to=2-2]
	\arrow["\ulcorner"{anchor=center, pos=0.125, rotate=180}, draw=none, from=2-2, to=1-1]
	\arrow["q", from=2-2, to=2-3]
	\arrow["\ulcorner"{anchor=center, pos=0.125, rotate=180}, draw=none, from=2-3, to=1-2]
\end{tikzcd}\]

\[\begin{tikzcd}
	{\mathcal C^\otimes_1} & {\mathcal O^\otimes_1} & {\mathrm N \Fint} \\
	{\mathcal{C}^{\oplus, \otimes}} & \opp & {\mathrm N \Finp}
	\arrow["{p_1}", from=1-1, to=1-2]
	\arrow[from=1-1, to=2-1]
	\arrow["{q_1}", from=1-2, to=1-3]
	\arrow[from=1-2, to=2-2]
	\arrow["{i_1}", from=1-3, to=2-3]
	\arrow["p", from=2-1, to=2-2]
	\arrow["\ulcorner"{anchor=center, pos=0.125, rotate=180}, draw=none, from=2-2, to=1-1]
	\arrow["q", from=2-2, to=2-3]
	\arrow["\ulcorner"{anchor=center, pos=0.125, rotate=180}, draw=none, from=2-3, to=1-2]
\end{tikzcd}\] Then $p_0$ and $p_1$ are cocartesian fibration of underlying operads.

    \item $p_0$ and $p_1$ are cocartesian fibrations. For every object $T \simeq T_1 \oplus \cdots \oplus T_m \in (\mathcal{O}^\otimes_0)_{\langle m\rangle}$ and $S \simeq S_1 \otimes \cdots \otimes S_m \in (\mathcal{O}^\otimes_1)_{\langle m\rangle}$, we have equivalences 
    \[
    (\mathcal{C}^\otimes_0)_{T} \simeq \prod_{1 \leq i \leq m} (\mathcal{ C}^\otimes_0)_{T_i} \quad (\mathcal{C}^\otimes_1)_{S} \simeq \prod_{1 \leq i \leq m} (\mathcal{ C}^\otimes_1)_{S_i}
    \] which are induced by inert maps $T \to T_i$ and $S \to S_i$.
        
    \end{enumerate}
    Then we have (a) $\Leftrightarrow$ (b) $\Rightarrow$ (c) $\Leftrightarrow$ (d).
\end{prop}

\begin{proof}
    The statement (c) $\Leftrightarrow$ (d) is \cite[Proposition 2.1.2.12]{lurie2017higher}. Thes same proposition shows (a) $\Rightarrow$ (b). By \cite[Proposition 2.4.1.3]{lurie2009higher}, $p_0 $ and $p_1$ are cocartesian fibrations. The composition $q_0 \circ p_0$ and $q_1 \circ p_1$ exhibit $\mathcal{ C}^\otimes_0$ and $\mathcal C ^\otimes_1$ as $\infty$-operads by \autoref{underlyingopandmop}. This proves (a) $\Rightarrow$ (c). We now prove (b) $\Rightarrow$ (a). Suppose we have a cocartesian fibration $p \colon \mathcal{C}^{\oplus, \otimes} \to \opp$ that satisfies the condition in (b). We must show that the composition $\mathcal{C}^{\oplus, \otimes} \to \mathrm N \Fint$ satisfies the axioms in \autoref{operadpairs}. The condition (1)-(3) essentially follows from the proof of \cite[Proposition 2.1.2.12]{lurie2017higher}. We need only to check (4).
    
    We let $T \in \opp_{A^Q_\times}$ be of the form
    \[
    T =\bigoplus_{\alpha \in Q} (T ^{1,q^\alpha_1} \otimes \cdots \otimes T^{k,q^\alpha_k})
    \] Let $C \in \mathcal{C}^{\oplus, \otimes}_T$ satisfies that $(\rho^{m,q_m,l_1}_\times)_!(C) \simeq (\rho^{m,q_m,l_2}_\times)_!(C) $ for any $m,q_m$ and $1 \leq l_1 \leq l_2 \leq n(m)$. $p$ sends cocartesian edges to cocartesian edges that exists, so $T$ also satisfies that $T^{m,q_m,l_1}=(\rho^{m,q_m,l_1}_\times)_!(T) \simeq (\rho^{m,q_m,l_2}_\times)_!(T) = T^{m,q_m,l_2}$. Since $\opp \to \mathrm{N} \Fint$ is an $\infty$-operad pair, we can find a distributive morphism $\tilde d_Q\colon T \to T'$ lifting $d_{Q}$, fitting into a composition diagram 
\[\begin{tikzcd}
	T & {T'} \\
	& {T^{m,q_m}}
	\arrow["{\tilde d_Q}", from=1-1, to=1-2]
	\arrow["{\tilde\rho^{m,q_m,l}_\times}"', from=1-1, to=2-2]
	\arrow["{\tilde\rho_+^{m,q_m}}", from=1-2, to=2-2]
\end{tikzcd}\] Since $p$ is a cocartesian fibration, the diagram lift to a composition diagram in $\mathcal C ^{\oplus,\otimes}$

\[\begin{tikzcd}
	C & {C'} \\
	& {C^{m,q_m}}
	\arrow["{\tilde d_Q}", from=1-1, to=1-2]
	\arrow["{\tilde\rho^{m,q_m,l}_\times}"', from=1-1, to=2-2]
	\arrow["{\tilde\rho_+^{m,q_m}}", from=1-2, to=2-2]
\end{tikzcd}\] This finishes the proof of (b) $\Rightarrow$ (a).
\end{proof}

\begin{cor}\label{laxmonoidal}
    The construction $\mathsf{Lax}^{\oplus, \otimes} \to \mathrm{N} \Fint $ and $\mathsf{OpLax}^{\oplus , \otimes} \to \mathrm{N} \Fint$ in \autoref{exampleofopp} become $\infty$-operad pairs.
\end{cor}

\begin{proof}
    By construction, these two functors clearly become cocartesian fibrations and they satisfy (b) in \autoref{settingmonoidal}.
\end{proof}

In \cite{lurie2017higher}, Lurie defines an $\mathcal{O}^\otimes$-monoidal category as a cocartesian fibration of operad pairs $p\colon \mathcal C^\otimes \to \mathcal O^\otimes$. However, for operad pairs this does not make sense. The cocartesian fibration condition is somehow too strong and it will leave redundant data. To see this, suppose that $p\colon \opp \to \Comm^{\oplus, \otimes}$ is a cocartesian fibration of operad pairs, then for any $C \in \opp_{A^Q_\times}$, there is always a cocartesian lift of $d_Q$ starting from $C$. This lift also fits into the same composition diagram. As a consequence, we always have $(\rho^{m,q_m,l_1}_\times)_{!}(C) \simeq (\rho^{m,q_m,l_2}_\times)_{!}(C)$, meaning that $\opp$ has only one color.

To fix the problem, we must reduce the amount of cocartesian edges. For example, we do not want any cocartesian lift of $d_Q$ starting from $A \otimes B \oplus C \otimes D$ for $A$, $B$, $C$ and $D$ that are not equivalent to each other.

\begin{defn}
    We say that a morphism $f \in \Fint$ is \emph{non-distributive} if it can be written as $\alpha = \alpha_1 \bullet \cdots \bullet \alpha_n$, where each bullet here represents a composition, a $\boxplus$ or a $\boxtimes$ and each $\alpha_k$ is an $f^+, f^\times, \sigma^+$ or $\sigma^\times$ as in [\autoref{fin2}, (a)-(c)]. where $\sigma$ is a permutation and $f$ is in $\Finp$. Let $\opp$ be an $\infty$-operad pair. We say a morphism in $\opp$ is non-distributive if it is a lift of some non-distributive $f \in \Fint$.
\end{defn}

\begin{defn}\label{bimonoidal}
    Let $\mathcal{C}^{\oplus, \otimes}$ and $\opp$ be two $\infty$-operad pairs. We say a functor $p\colon \mathcal{C}^{\oplus, \otimes} \to \opp$ exhibits $\mathcal{C}^{\oplus,\otimes}$ as a $\opp$-\emph{bimonoidal category} if 
\begin{enumerate}[(1)]
        \item $p$ is a map of operad pairs.
        \item For any non-distributive $g\colon X \to X' \in \opp$ and any $C \in p^{-1} (X)$, there exists a cocartesian lift of $g$ starting from $C$.
        \item For any distributive $d\colon Y \to Y' \in \opp$ and any $C \in p^{-1}(Y)$ such that $(\rho^{m,q_m,l_1}_\times)_{!}(C) \simeq (\rho^{m,q_m,l_2}_\times)_{!}(C)$, there exists a cocartesian lift $\tilde d$ of $d$, fitting into a cocartesian lift of the composition diagram in $\opp$:
\[\begin{tikzcd}
	C & {C'} \\
	& {C^{m,q_m}}
	\arrow["{\tilde d}", from=1-1, to=1-2]
	\arrow["{\tilde\rho^{m,q_m,l}_\times}"', from=1-1, to=2-2]
	\arrow["{\tilde\rho_+^{m,q_m}}", from=1-2, to=2-2]
\end{tikzcd}\] Note that $\tilde d$ must be the distributive morphism starting from $C$ in the operad pair $\mathcal C^{\oplus, \otimes}$.
\end{enumerate}
\end{defn}

\begin{cor}\label{settingcocartandcart}
    Let $\opp$ be an $\infty$-operad pair and let $p\colon \mathcal C^{\oplus,\otimes} \to \opp$ be a functor that satisfies (2) and (3) in \autoref{bimonoidal}. Then $p$ exhibits $\mathcal C^{\oplus,\otimes}$ as a $\opp$-bimonoidal category if and only if for every $T = T_1 \odot \cdots \odot T_n \in \opp$, we have
    \[
    (\mathcal{C}^{\oplus,\otimes})_{T} \simeq \prod_{1 \leq i \leq m} (\mathcal{ C}^{\oplus,\otimes})_{T_i} \]
\end{cor}

\begin{proof}
    This follows from a same argument as in \autoref{settingmonoidal}.
\end{proof}

We end this subsection by providing some examples of $\opp$-bimonoidal categories for an $\infty$-operad pair $\opp$.

\begin{exmp}\label{monoidalexample}
    Consider $\Comm^{\oplus,\otimes}$ in \autoref{exampleofopp}. We call $\Comm^{\oplus,\otimes}$-bimonoidal categories \emph{symmetric bimonoidal categories} or \emph{symmetric rig categories}. Informally speaking, if we choose cocartesian lifts of $\beta^+\colon (\langle2 \rangle,0) \to \langle 1 \rangle$ and $\beta^\times \colon (\langle 2\rangle,1 ) \to \langle 1 \rangle$, we obtain two symmetric monoidal structures

    \[
    + \colon \mathcal{P}\times \mathcal P \to \mathcal P \quad (A,B) \mapsto \beta^{+}_! (A \oplus B)
    \] and 

    \[
    \times \colon \mathcal{P}\times \mathcal P \to \mathcal P \quad (A,B) \mapsto \beta^{\times}_! (A \otimes B)
    \] Two maps $\beta_0^+, \beta_0^\times \colon \langle 0 \rangle \to (\langle1 \rangle,0) $ lift to two cocartesian edges $\emptyset \to \mathcal I^+$ and $\emptyset \to \mathcal I^\times$, corresponding to the units of two symmetric monoidal structure. Moreover, the diagram in (iii) \autoref{fin2} lifts to a natural distributivity equivalence $$A \times B + A \times C\simeq A \times(B +C)$$. 
\end{exmp}

\begin{exmp}
    Consider $\Ass^{\oplus,\otimes}$ in \autoref{exampleofopp}. We call $\Ass^{\oplus,\otimes}$-bimonoidal categories \emph{bimonoidal categories} or \emph{rig categories}. Similar construction gives a commutative associative addition $+ \colon \mathcal P \times \mathcal P \to \mathcal P$ and an associative multiplication $\times \colon \mathcal P \times \mathcal P \to \mathcal{P}$ together with a distributivity equivalence. 
\end{exmp}

\begin{warn}
    In general, we do not require $\mathcal I^+ \otimes 0 \simeq 0$ or $0 \otimes \mathcal I^+ \simeq 0$ in our definition of (symmetric) rig category. This is a consequence of our banishment of $d_Q$ for $|Q|=0$. However, many examples satisfy this compatibility of units. For instance, if the rig structure comes from a closed symmetric monoidal structure, then this is satisfied automatically.
\end{warn}

\begin{exmp}
    Consider $\mathsf{Lax}^{\oplus, \otimes}$ and $\mathsf{OpLax}^{\oplus, \otimes}$ in \autoref{exampleofopp}. We call $\mathsf{Lax}^{\oplus, \otimes}$ (or $\mathsf{OpLax}^{\oplus, \otimes}$)-bimonoidal categories \emph{lax} (or \emph{oplax}) \emph{bimonoidal} categories. The underlying category of these bimonoidal categories are equipped with a commutative associative addition and a commutative associative multiplication together with a lax distributivity law $A\times B + A \times C \to A \times(B+C)$ (or oplax distributivity law $A \times(B+C) \to A\times B + A \times C $). 
\end{exmp}
\begin{rem}
    It seems reasonable that one can also use the method in this section to define $n$-fold operad pairs and $n$-monoidal categories. Roughly, we need to replace the universal operad pair $\Fint$ by $\mathsf{Fin}^n_*$ where the objects are $\langle m \rangle \in \Finp$ together with an $n$-nary number to distinguish a combination of $n$ different operators. The morphisms should cover the data of inert, active and $n$-fold distribute data. We define an $n$-fold $\infty$-operad pair as a functor $\mathcal{C}^{\otimes^n} \to \mathrm{N}\mathsf{Fin}^n_*$ such that all the data is compatible. For $0\leq k \leq n$ we can pullback along different inclusions $\mathsf{Fin}^k_* \to \mathsf{Fin}^n_*$ to obtain underlying $k$-fold operad pairs. It is expected that all the properties in this section for operad pairs can be generalized to $n$-fold cases for arbituary $n$.
\end{rem}

\section{Comparison with the classical definition}\label{comparesection}
    The $\infty$-operads can be viewed as generalization of classical ones such that classical operads are those operads with a single color. More precisely, let $\mathcal{O}^\otimes$ be an operad, we define $\mathrm{Mul}_{\mathcal O}(\{ X_i\}_{1\leq i \leq n}, Y)$ as the component $\Map^{\beta_n} (X_1 \oplus\cdots \oplus X_n, Y)$ of the total mapping space. For those operads with a single color, one can define $\mathscr{C}(n):= \mathrm{Mul}_{\mathcal O}(\{\mathbbl{1}\}_{1\leq i \leq n}, \mathbbl{1})$, which is a classical operad that receives the operadic product map $\gamma$ from the composition of maps. 

    For operad pairs, the same story happens. We will compare the original definition of \cite{may2009construction} and the single colored operad pair in our definition. We first recall the original definition.

    \begin{defn}[\cite{may2009construction}, \cite{Mayinprogress}]\label{classicaldefn}
        An operad pair $(\mathscr{C}, \mathscr{J})$ in $\mathcal{S}$ consists of two (classical) operads $\mathscr{C}$ and $\mathscr{J}$ and an action 
        \[
        \lambda \colon  \mathscr{J}(k) \times \mathscr{C}(j_1) \times \cdots \times \mathscr{C}(j_k) \to \mathscr{C}(j_1\cdots j_k)
        \] that satisfies the following properties:

\begin{enumerate}[(i)]

\item The operad $\mathscr C$ is unital. Therefore $\lambda$ maps to $*$ if any $j_r=0$.
\item The operad $\mathscr J$ is unital. When $k=0$, $\lambda$ is identified with $\eta \colon *\to \mathscr C(1)$.
\item The following unitary laws are satisfied.
\[\begin{tikzcd}[nodes={scale=0.8}]
	{\mathscr J(k)\times *^k} & {*} && {* \times \mathscr{C}(k)} & {\mathscr C(k)} \\
	{\mathscr{J}(k)\times \mathscr{C}(1)^k} & {\mathscr C(1)} && {\mathscr J(1) \times \mathscr C(k)}
	\arrow["\lambda", from=1-1, to=1-2]
	\arrow["{\mathrm{id}\times\eta_+^k}"', from=1-1, to=2-1]
	\arrow["{\eta_+}", from=1-2, to=2-2]
	\arrow["\simeq", from=1-4, to=1-5]
	\arrow["{\eta_\times\times \mathrm{id}}"', from=1-4, to=2-4]
	\arrow["\lambda", from=2-1, to=2-2]
	\arrow[from=2-4, to=1-5]
\end{tikzcd}\]

\item The associativity law should satisfy. That is, we have the following commutative square:
\[\begin{tikzcd}[nodes={scale=0.8}, sep=small]
	{\mathscr {J}(k)\times(\times_{r=1}^k \mathscr J(j_r))\times (\times_{t=1}^{j_1+\cdots + j_k}\mathscr C(i_t))} & {\mathscr J(\Sigma j_r)\times (\times_{t=1}^{\Sigma j_r} \mathscr C(i_t))} \\
	& {\mathscr C (i_1 \cdots i_{j_1 + \cdots+j_k})} \\
	{\mathscr {J}(k)\times(\times_{r=1}^k \mathscr J(j_r))\times (\times_{r=1}^k \times_{q_r=1}^{j_r}\mathscr C(i_{q_1}))} & {\mathscr J(k) \times (\times_{r=1}^k \mathscr C (\times_{q_r = 1}^{j_r} i_{q_r}))}
	\arrow["{\mathrm{\gamma_\times} \times\mathrm{id}}", from=1-1, to=1-2]
	\arrow["{\mathrm{reshuffle}}"', from=1-1, to=3-1]
	\arrow["\lambda", from=1-2, to=2-2]
	\arrow["{\mathrm{id}\times \lambda^k}", from=3-1, to=3-2]
	\arrow["\lambda"', from=3-2, to=2-2]
\end{tikzcd}\]
\item The action $\lambda$ should be equivariant. Let notations be as in \autoref{twonotnforeq}, we have following commutative diagrams:

\[\begin{tikzcd}[nodes={scale=0.8}]
	{\mathscr {J}(k)\times \mathscr C (j_1) \times \cdots \times \mathscr C(j_k)} & {\mathscr J (k) \times \mathscr C (j_{\sigma(1)}) \times \cdots \times \mathscr C(j_{\sigma(k)})} \\
	{\mathscr C(j_1 \cdots j_k)} & {\mathscr C(j_1 \cdots j_k)}
	\arrow["{\sigma \times \sigma^{-1}}", from=1-1, to=1-2]
	\arrow["\lambda"', from=1-1, to=2-1]
	\arrow["\lambda", from=1-2, to=2-2]
	\arrow["{\sigma(\langle j_{\sigma(1)} ,\cdots, j_{\sigma (k)}\rangle)}"', from=2-1, to=2-2]
\end{tikzcd}\]

\[\begin{tikzcd}[nodes={scale=0.8}]
	{\mathscr {J}(k)\times \mathscr C (j_1) \times \cdots \times \mathscr C(j_k)} && {\mathscr J (k) \times \mathscr C (j_1) \times \cdots \times \mathscr C(j_k)} \\
	{\mathscr C(j_1 \cdots j_k)} && {\mathscr C(j_1 \cdots j_k)}
	\arrow["{\mathrm{id}\times \tau_1\times \cdots \times \tau_k}", from=1-1, to=1-3]
	\arrow["\lambda"', from=1-1, to=2-1]
	\arrow["\lambda", from=1-3, to=2-3]
	\arrow["{\tau_1 \otimes\cdots \otimes \tau_k}"', from=2-1, to=2-3]
\end{tikzcd}\]

\item Let $Q \in \mathrm{Seq}(\mathbb N)$ be as in \autoref{Qanddistequivalence}. Suppose that for each $(m, q_m)$ there is a number $i_{m,q_m}$ consisting a distributive equation:
\[
I:=\coprod_{Q} \prod_{1\leq r\leq k }i_{r,q_r} \simeq \prod_{1\leq r \leq k}\coprod_{1\leq q \leq j_r}i_{r,q}
\]Then we have the following distributivity diagram:
\[\begin{tikzcd}[nodes={scale=0.8}]
	{\mathscr J(k) \times \big(\times_{r=1}^k \mathscr C(j_r) \times (\times_{q=1}^{j_r}\mathscr C (i_{r,q}) )\big)} & {\mathscr J(k) \times \big(\times_{r=1}^k \mathscr C(l_r) \big)} \\
	{\mathscr J(k)^{J+1} \times \big(\times_{r=1}^k \mathscr C(j_r) \times (\times_{q=1}^{j_r}\mathscr C (i_{r,q})^{n(r)} )\big)} \\
	{\mathscr J(k) \times \big(\times_{r=1}^k  \mathscr C(j_r)\big) \times \big(\times_Q( \mathscr J(k)\times\times_{r=1}^k\mathscr C (i_{r,q_r}) )\big)} & {\mathscr C(I)} \\
	{\mathscr J(k) \times \big(\times_{r=1}^k \mathscr C(j_r)\big) \times \big(\times_Q\mathscr C (i_Q) \big)} & {\mathscr C(J) \times \big(\times_{Q} \mathscr C(i_Q) \big)}
	\arrow["{\mathrm{id}\times \gamma^k_+}", from=1-1, to=1-2]
	\arrow["{\Delta\times(\times_r (\mathrm{id} \times\Delta))}"', from=1-1, to=2-1]
	\arrow["\lambda", from=1-2, to=3-2]
	\arrow["\simeq"', from=2-1, to=3-1]
	\arrow["{\mathrm{id}\times \lambda^J}"', from=3-1, to=4-1]
	\arrow["{\lambda\times\mathrm{id}}"', from=4-1, to=4-2]
	\arrow["{\gamma_+}"', from=4-2, to=3-2]
\end{tikzcd}\] Here $l_r=\sum_{q=1}^{j_r}i_{r,q}$, $i_Q=\prod_{r=1}^k i_{r,q_r}$ and $\Delta$ is the diagonal map on spaces $\mathscr J(k)$ and $\mathscr C (i_{r,q_r})$.
\end{enumerate}
\end{defn}

In the above definition. We require both $\mathscr C $ and $\mathscr J$ to be unital. In order to be compatible with this classical definition, we need a counterpart in the $\infty$-categorical setting. Recall that an $\infty$-operad $\mathcal{O}^\otimes$ is called \emph{unital} if for any $X \in \mathcal{O}$, the space $\Map (\emptyset , X) \simeq \mathrm{Mul}(\emptyset,X)$ is contractible. 

\begin{defn}
    Let $p\colon\opp\to \mathrm{N}\Fint$ be an $\infty$-operad pair. Let $\beta_0^+, \beta_0^\times \colon \langle 0 \rangle \to \langle 1 \rangle$ be the maps associated to $\beta_0\colon \langle 0 \rangle \to \langle 1 \rangle \in \Finp$. We use $\mathrm{Mul}^+(\emptyset,X)$ and $\mathrm{Mul}^\times (\emptyset, X)$ to denote $\Map^{\beta_0^+}(\emptyset, X)$ and $\Map^{\beta_0^\times}(\emptyset, X)$ respectively. Then we say that the $\infty$-operad pair $\opp$ is \emph{unital} if 
    
    \begin{enumerate}
        \item  The underlying operad and moperad of $\opp$ are unital. That is, both $\mathrm{Mul}^+(\emptyset,X)$ and $\mathrm{Mul}^\times (\emptyset, X)$ are contractible. By the Segal condition, this is equivalent to say that $\Map^g (\emptyset , X)$ is contractible for any $X \in \opp_{(\langle m \rangle,0)}$ and any $g\colon \langle 0\rangle \to (\langle m \rangle,0)$ in $\Fint$.
        \item The map $\beta_0^+ \colon \langle 0 \rangle \to \langle 1 \rangle$ admits a cocartesian lift $\tilde\beta_0^+\colon\emptyset \to \mathcal{I}$.
    \end{enumerate}
\end{defn}

\begin{rem}
    The requirement (2) can be viewed as a cocartesian lift requirement for degenerated distributive maps $d_Q \in \Fint$, where $Q=(k,j_1,\cdots,j_k)$ such that $j_1=1$ and $j_2=\cdots = j_k=0$. 
\end{rem}


The advantage of working with unital $\infty$-operad pairs is the contractible requirement allows us to have more control over the higher mapping spaces

\begin{lem}\label{unitlemma}
    Suppose that $p\colon \opp \to \mathrm{N}\Fint$ is a unital $\infty$-operad pair and let $\mathcal I$ be an object in $\mathcal{P}\simeq\opp_{\langle 1 \rangle}$. Then for any $k\geq1$, there exists a morphism $\kappa \colon \mathcal I \to \mathcal I^{\otimes k}$ such that for any $f \in \Map^{\beta^\times_k} (\mathcal I^{\otimes k}, \mathcal I)$, $f \kappa \simeq \mathrm{id}_{\mathcal I}$.
\end{lem}

\begin{proof}
    We choose a section of $\beta^\times_k$ in $\Fint$. For example, we can always choose $\underline{\kappa}^\times \colon \langle 1 \rangle \to (\langle k \rangle , 2^{k-1}-1)$ for the injection $\underline{\kappa} \colon \langle 1 \rangle \to \langle k \rangle $ on $1 \in \langle k \rangle $. By the Segal condition and unital requirement, there always exists a morphism $\kappa \colon \mathcal I \to \mathcal I^{\otimes k}$ lifting $ \underline{\kappa}^\times$. Consider the following diagram such that $g$ is the composition of $\kappa $ and $\tilde\beta_0^+$.

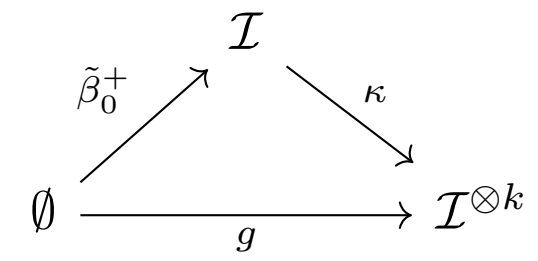
    \[\begin{tikzcd}
	& {\mathcal I} & \\
	\emptyset && {\mathcal I^{\otimes k}}
	\arrow["\kappa", from=1-2, to=2-3]
	\arrow["{\tilde\beta_0^+}", from=2-1, to=1-2]
	\arrow["{g}"', from=2-1, to=2-3]
\end{tikzcd}\] For any $f \colon \mathcal I^{\otimes k} \to \mathcal I$ lifting $\beta^\times_k$, we have the following square:
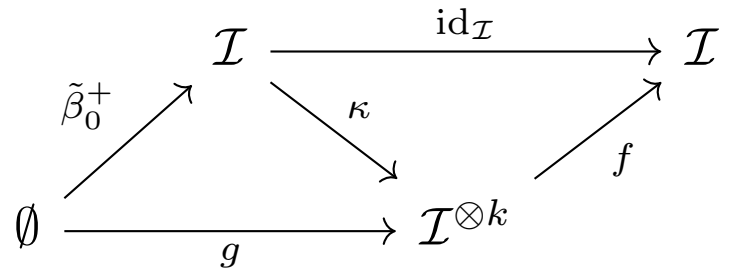
\[\begin{tikzcd}
	& {\mathcal I} && {\mathcal I} \\
	\emptyset && {\mathcal I^{\otimes k}}
	\arrow["{\mathrm{id}_{\mathcal I}}", from=1-2, to=1-4]
	\arrow["\kappa", from=1-2, to=2-3]
	\arrow["{\tilde\beta_0^+}", from=2-1, to=1-2]
	\arrow["{g}"', from=2-1, to=2-3]
	\arrow["f"', from=2-3, to=1-4]
\end{tikzcd}\] This is taken to a commutative diagram in $\mathrm{N}\Fint$. Therefore, the outer square commutes because of the contractibility of the mapping space. As a consequence we must have $f \kappa  \simeq\mathrm{id}_{\mathcal I}$ since $\tilde\beta_0^+$ is a cocartesian edge.
\end{proof}

The above proof in fact tells us more.

\begin{cor}
    Let $p\colon \opp \to \mathrm{N}\Fint$ be a unital $\infty$-operad pair, then unit mapping space $\mathrm{Mul}^{+} (\mathcal I, \mathcal I) \simeq \mathrm{Mul}^{\times} (\mathcal I, \mathcal I)$ is contractible.
\end{cor}

We are going to prove the following:

\begin{thm}\label{comparisontheorem}
    Let $p\colon \opp \to \mathrm N \Fint$ be a unital $\infty$-operad pair such that $\mathcal{P}$ has only one element $\mathbbl{1}$. We define 
    \[
    \mathscr C(n) := \mathrm{Mul}^{+} (\{\mathbbl{1}\}_{n}, \mathbbl{1})=\Map^{\beta_n^+} (\mathbbl{1}^{\oplus n}, \mathbbl{1}) \]
    
    \[\mathscr J(n):= \mathrm{Mul}^{\times} (\{\mathbbl{1}\}_{n}, \mathbbl{1})=\Map^{\beta^\times_n} (\mathbbl{1}^{\otimes n}, \mathbbl{1})
    \]
    Then there exists an action  \[
        \lambda \colon  \mathscr{J}(k) \times \mathscr{C}(j_1) \times \cdots \times \mathscr{C}(j_k) \to \mathscr{C}(j_1\cdots j_k)
        \] making $(\mathscr C , \mathscr J)$ into an operad pair in the sense of \autoref{classicaldefn}.
\end{thm}

The first thing we need to do is to construct the desired $\lambda$. Before that, we need some general notations.

\begin{notn}
    Let $\opp \to \mathrm{N}\Fint$ be an $\infty$-operad pair,  $f_1\colon C_1\to D_1$ and $f_2 \colon C_2 \to D_2$ be two morphisms in $\opp$ lifting $\alpha_1$ and $\alpha_2$ in $\Fint$. Then for $\odot$ that represents a copy of $\otimes $ or $\otimes$, we define $f_1 \odot f_2 \colon C_1 \odot C_2 \to D_1 \odot D_2$ as the morphism in $\Map^{\alpha_1 \bullet \alpha_2} (C_1 \odot C_2 , D_1 \odot D_2)$ corresponding to 

\[(f_1 \tilde\rho_1, f_2 \tilde\rho_2) \in \Map^{\rho_1 \alpha_1 \bullet \alpha_2} (C_1 \odot C_2, D_1) \times \Map^{\rho_2 \alpha_1 \bullet \alpha_2} (C_1 \odot C_2, D_2)\] via the equivalence. Here $\tilde \rho_i \colon C_1 \odot C_2 \to D_i$ denote the cocartesian lift of the inert morphism that represents projection on each term.
\end{notn}

\begin{lem}\label{functoraldist}
    Let$A^Q_\times, A_+^Q,B_\times^Q$ and $B_+^Q$ be notations in [\autoref{fin2},(viii)]. Let  $$\bigoplus_{\alpha \in \underline{Q}} (C ^{1,q^\alpha_1} \otimes \cdots \otimes C^{k,q^\alpha_k}) \in \opp_{A^\times_Q}, \quad \bigoplus_{\alpha \in \underline{Q}} (D ^{1,q^\alpha_1} \otimes \cdots \otimes D^{k,q^\alpha_k}) \in \opp_{B^\times_Q}$$ and suppose that we have a morphism 
    \[
    \bigoplus_{\alpha \in \underline{Q}} (f ^{1,q^\alpha_1} \otimes \cdots \otimes f^{k,q^\alpha_k}) \colon \bigoplus_{\alpha \in \underline{Q}} (C ^{1,q^\alpha_1} \otimes \cdots \otimes C^{k,q^\alpha_k}) \to \bigoplus_{\alpha \in \underline{Q}} (D ^{1,q^\alpha_1} \otimes \cdots \otimes D^{k,q^\alpha_k})
    \] Then the commutative square in [\autoref{fin2},(iii)] in $\Fint$ lifts to a commutative square in $\opp$:
\[\begin{tikzcd}[nodes={scale=0.8}]
	{\bigoplus_{\alpha \in \underline{Q}} (C ^{1,q^\alpha_1} \otimes \cdots \otimes C^{k,q^\alpha_k})} & {\bigotimes_{1 \leq r\leq k} (C^{r,1}\oplus \cdots \oplus C^{r,j_r})} \\
	{\bigoplus_{\alpha \in \underline{Q}} (D ^{1,q^\alpha_1} \otimes \cdots \otimes D^{k,q^\alpha_k})} & {\bigotimes_{1 \leq r\leq k} (D^{r,1}\oplus \cdots \oplus D^{r,j_r})}
	\arrow["{\tilde d_Q}", from=1-1, to=1-2]
	\arrow["{\bigoplus_{\alpha \in \underline{Q}} (f ^{1,q^\alpha_1} \otimes \cdots \otimes f^{k,q^\alpha_k})}"', from=1-1, to=2-1]
	\arrow["{\bigotimes_{1 \leq r\leq k} (f^{r,1}\oplus \cdots \oplus f^{r,j_r})}", from=1-2, to=2-2]
	\arrow["{\tilde d_Q}", from=2-1, to=2-2]
\end{tikzcd}\]
\end{lem}

\begin{proof}
    We need to show that the composition clockwise and counterclockwise coincide in the mapping space of $\opp$. By (ii) of \autoref{operadpairs}, it suffices to check these two maps coincide after apply $\rho_i$ on $\bigotimes_{1 \leq r\leq k} (D^{r,1}\oplus \cdots \oplus D^{r,j_r})$. We consider the following diagram:

    \[\begin{tikzcd}[nodes={scale=0.8}]
	{\bigoplus_{\alpha \in \underline{Q}} (C ^{1,q^\alpha_1} \otimes \cdots \otimes C^{k,q^\alpha_k})} & {\bigotimes_{1 \leq r\leq k} (C^{r,1}\oplus \cdots \oplus C^{r,j_r})} \\
	{\bigoplus_{\alpha \in \underline{Q}} (D ^{1,q^\alpha_1} \otimes \cdots \otimes D^{k,q^\alpha_k})} & {\bigotimes_{1 \leq r\leq k} (D^{r,1}\oplus \cdots \oplus D^{r,j_r})} \\
	{D^{m,q_m}} & {D^{m,q_m}}
	\arrow["{\tilde d_Q}", from=1-1, to=1-2]
	\arrow["{\bigoplus_{\alpha \in \underline{Q}} (f ^{1,q^\alpha_1} \otimes \cdots \otimes f^{k,q^\alpha_k})}"', from=1-1, to=2-1]
	\arrow["{\bigotimes_{1 \leq r\leq k} (f^{r,1}\oplus \cdots \oplus f^{r,j_r})}", from=1-2, to=2-2]
	\arrow["{\tilde d_Q}", from=2-1, to=2-2]
	\arrow["{\tilde\rho^{m,q_m,l}_\times}"', from=2-1, to=3-1]
	\arrow["{\tilde\rho^{m,q_m}_+}", from=2-2, to=3-2]
	\arrow["\simeq", from=3-1, to=3-2]
\end{tikzcd}\] The lower square commutes by compatibility of $\rho $ and $d$ ([\autoref{operadpairs}, (4)]). So in order to show the upper square commutes after applying $\rho^{m,q_m}$, it suffices to show the outer square commutes. The outer square is identified with the following square:

\[\begin{tikzcd}[nodes={scale=0.8}]
	{\bigoplus_{\alpha \in \underline{Q}} (C ^{1,q^\alpha_1} \otimes \cdots \otimes C^{k,q^\alpha_k})} & {\bigotimes_{1 \leq r\leq k} (C^{r,1}\oplus \cdots \oplus C^{r,j_r})} \\
	{C^{m,q_m}} & {C^{m,q_m}} \\
	{D^{m,q_m}} & {D^{m,q_m}}
	\arrow["{\tilde d_Q}", from=1-1, to=1-2]
	\arrow["{\tilde\rho^{m,q_m,l}_\times}"', from=1-1, to=2-1]
	\arrow["{\tilde\rho^{m,q_m}_+}", from=1-2, to=2-2]
	\arrow["\simeq", from=2-1, to=2-2]
	\arrow["{f^{m,q_m}}"', from=2-1, to=3-1]
	\arrow["{f^{m,q_m}}", from=2-2, to=3-2]
	\arrow["\simeq", from=3-1, to=3-2]
\end{tikzcd}\] Again the upper square commutes because of the compatibility of $\rho$ and $d$ and the lower square trivially commutes.
\end{proof}

\begin{rem}\label{specialkappa}
    Let $\opp$ be a unital $\infty$-operad pair and let $\kappa \colon \mathcal I \to \mathcal I^{\otimes k}$ be the morphism lifting $\underline{\kappa}^\times \colon \langle 1 \rangle \to (\langle k \rangle , 2^{k-1}-1)$. We can view $\mathcal{I}^{\oplus |Q|}$ as a trivial product $(\mathcal I \otimes \underbrace{\emptyset \otimes \cdots \otimes \emptyset}_{k-1})^{\oplus|Q|}$. The distributive map on $\mathcal{I}^{\oplus |Q|}$ becomes a projection, which we denote by $\mathrm{pr}$. \autoref{functoraldist} then specifies to give a commutative diagram: 
    \[\begin{tikzcd}
	{\mathcal I^{\oplus |Q|}} & {\mathcal I^{\oplus j_1}} \\
	{(\mathcal I^{\otimes k})^{\oplus |Q|}} & {\mathcal I^{\oplus j_1}\otimes \cdots \otimes \mathcal I^{\oplus j_k}}
	\arrow["{\mathrm{pr}}", from=1-1, to=1-2]
	\arrow["{\kappa^{\oplus |Q|}}"', from=1-1, to=2-1]
	\arrow["{\mathrm{id}\otimes * \otimes \cdots \otimes *}", from=1-2, to=2-2]
	\arrow["{d_Q}"', from=2-1, to=2-2]
\end{tikzcd}\]
\end{rem}

\begin{con}\label{constructionofdist}
    Let notations be as in \autoref{comparisontheorem}. Then we construct $\lambda$ as follows:

\[\begin{tikzcd}[nodes={scale=0.8}]
	{\mathscr J(k) \times\mathscr C (j_1)\times\cdots \times \mathscr C (j_k) } & {\mathscr C(j_1 \cdots j_k)} \\
	{\Map^{\beta^\times_k}(\mathbbl{1}^{\otimes k},\mathbbl 1) \times \Map^{\beta^+_{j_1}}(\mathbbl{1}^{\oplus j_1},\mathbbl 1)\times \cdots \times \Map^{\beta^+_{j_k}}(\mathbbl{1}^{\oplus j_k},\mathbbl 1)} & {\Map^{\beta^+_{|Q|}}(\mathbbl 1^{\oplus |Q|}, \mathbbl 1)} \\
	{\Map^{\beta^\times_k}(\mathbbl{1}^{\otimes k},\mathbbl 1) \times \Map^{\beta^+_{j_1}\boxtimes\cdots \boxtimes\beta^+_{j_k}}(\mathbbl 1^{\oplus j_1}\otimes\cdots \otimes \mathbbl 1^{\oplus j_k},\mathbbl 1 ^{\otimes k})} & {\Map^{\beta^+_{|Q|} \circ \beta^\times_k}((\mathbbl 1^{\otimes k})^{\oplus |Q|},\mathbbl 1 )} \\
	{\Map^{\beta^\times_k}(\mathbbl{1}^{\otimes k},\mathbbl 1) \times \Map^{\beta^+_{j_1}\boxtimes\cdots \boxtimes\beta^+_{j_k}\cdot d_Q}((\mathbbl 1^{\otimes k})^{\oplus |Q|},\mathbbl 1 ^{\otimes k})} & {\Map^{\beta^\times_k \circ\beta^+_{j_1}\boxtimes\cdots \boxtimes\beta^+_{j_k}\circ d_Q}((\mathbbl 1^{\otimes k})^{\oplus |Q|},\mathbbl 1 )}
	\arrow[dashed, from=1-1, to=1-2]
	\arrow["\simeq"', from=1-1, to=2-1]
	\arrow["{\mathrm {id}\times\otimes^{k-1}}"', from=2-1, to=3-1]
	\arrow["\simeq"', from=2-2, to=1-2]
	\arrow["{\mathrm{id}\times(\text{precomposition with } \tilde d_Q)}"', from=3-1, to=4-1]
	\arrow["{\text{composition with }\kappa^{\oplus |Q|}}"', from=3-2, to=2-2]
	\arrow["{\text{composition}}"', from=4-1, to=4-2]
	\arrow["\simeq"', from=4-2, to=3-2]
\end{tikzcd}\] Here the lower right equivalence is by the distributive diagram in \autoref{fin2} and $\kappa$ is defined as in \autoref{unitlemma}.

\end{con}

Now we are going to prove \autoref{comparisontheorem}. With $\lambda$ constructed above, we must verify that this construction satisfies the axioms in \autoref{classicaldefn}.

\begin{proof}[Proof of \autoref{comparisontheorem}]
    Suppose that $j_r=0$ for some $r$, then the mapping space $\Map^{\beta^\times} (\mathbbl1 ^{\oplus |Q|}, \mathbbl 1) \simeq \mathrm{Mul}(\emptyset, \mathbbl 1)$ is a point, so the composition square in \autoref{constructionofdist} maps everything to a point. Suppose that $k=0$, then the composition with $\kappa$ sends $* \in \mathrm{Mul (\emptyset,\mathbbl 1)}$ to $* \in \Map^{\beta_1^\times}(\mathbbl 1 , \mathbbl 1)$. This proves the axioms (i) and (ii). 

    Now we need to show (iii). For the first diagram in (iii), we need to show that when $j_1=\cdots = j_k=1$ for any $f \in \Map^{\beta^\times_k}(\mathbbl{1}^{\otimes k},\mathbbl 1) $, the composition in \autoref{constructionofdist} sends $(f, \mathrm{id},\cdots, \mathrm{id})$ to $\mathrm{id} \in \Map^{\beta^+_1}(\mathbbl 1 , \mathbbl 1)$. This follows from the fact that $d_Q=\mathrm{id}$ for $|Q|=1$ and $f\kappa=\mathrm{id}$ for any $f$ as in \autoref{unitlemma}. The second unit diagram in (iii) follows similarly from the fact that $d_Q=\mathrm{id}$ when $k=1$.

   We need to show (iv). The associative square in \autoref{classicaldefn} now translates into 
\[\begin{tikzcd}[nodes={scale=0.8},sep=small]
	{\Map^{\beta^\times_k}(\mathbbl 1^{\otimes k} , \mathbbl 1) \times \prod_{r=1}^k \Map^{\beta^\times_{j_r}}(\mathbbl 1^{\otimes {j_r}} , \mathbbl 1)\times \prod_{t=1}^{j_+}\Map ^{\beta^+_t}(\mathbbl 1^{\oplus i_t},\mathbbl 1)} & {\Map^{\beta^\times_{j_+}}(\mathbbl 1^{\otimes {j_+}} , \mathbbl 1)\times\prod_{t=1}^{j_+}\Map ^{\beta^+_t}(\mathbbl 1^{\oplus i_t},\mathbbl 1)} \\
	{\Map^{\beta^\times_k}(\mathbbl 1^{\otimes k} , \mathbbl 1)\times \prod_{r=1}^{k}\Map ^{\beta^+_{h_r}}(\mathbbl 1^{\oplus h_r},\mathbbl 1)} & {\Map^{\beta^+_I}(\mathbbl 1^{\oplus I},\mathbbl 1)}
	\arrow[from=1-1, to=1-2]
	\arrow[from=1-1, to=2-1]
	\arrow[from=1-2, to=2-2]
	\arrow[from=2-1, to=2-2]
\end{tikzcd}\] Now choose \[(f, \prod_{r=1}^k g_r , \prod_{t=1}^{j_+} h_t) \in \Map^{\beta^\times_k}(\mathbbl 1^{\otimes k} , \mathbbl 1) \times \prod_{r=1}^k \Map^{\beta^\times_{j_r}}(\mathbbl 1^{\otimes {j_r}} , \mathbbl 1)\times \prod_{t=1}^{j_+}\Map ^{\beta^+_t}(\mathbbl 1^{\oplus i_t},\mathbbl 1)\] Clockwise it is sent to 

\[f \cdot \big( g_1\cdot (h_{1,1} \otimes \cdots \otimes h_{1,j_1}) \tilde d_1 \kappa_1^{\oplus l_1}\otimes \cdots \otimes g_k\cdot (h_{k,1} \otimes \cdots \otimes h_{k,j_k}) \tilde d_k \kappa_k^{\oplus l_k}\big) \tilde d' \kappa' \] while counterclockwise it is sent to

\[
f\cdot (g_1\otimes \cdots \otimes g_k) \cdot (h_{1,1} \otimes\cdots \otimes h_{k,j_k}) \tilde d \kappa^{\oplus I} 
\] Here $l_r=\sum_{q=1}^{j_r}i_{r,q}$, $i_Q=\prod_{r=1}^k i_{r,q_r}$. We use  $\tilde d_r$ to denote the cocartesian lift of $d_{Q_r}$, where $Q_r=(j_r,i_{r,1},\cdots , i_{r,j_r})\in \mathrm{Seq}(\mathbb N)$. We use $\tilde d'$ to denote the cocartesian lift of $d_{Q'}$ where $Q' = (k,l_1,\cdots ,l_k)$. We use $\tilde d$ to denote the cocartesian lift of $d_Q$ where $Q=(l_1+\cdots +l_r, i_{1,1},\cdots,i_{k,j_k})$. It now suffices to show that $\tilde d \kappa^{\oplus I} \simeq \tilde (d_1 \otimes \cdots \otimes \tilde d _k )\cdot (\kappa_1^{\oplus l_1}\otimes \cdots \otimes \kappa_k^{\oplus l_k}) \tilde d' \kappa'$. This follows from the commutative diagram:

\[\begin{tikzcd}[nodes={scale=0.8}]
	{\mathbbl 1^{\oplus l_1\cdots l_k}} & {(\mathbbl 1^{\otimes k})^{\oplus l_1\cdots l_k}} & {\mathbbl 1^{\oplus l_1}\otimes \cdots\otimes \mathbbl 1^{\oplus l_k}} & \\
	& {(\mathbbl 1^{\otimes j_+})^{\oplus l_1 \cdots l_k}} & {(\mathbbl 1^{\otimes j_1})^{\oplus l_1}\otimes \cdots\otimes (\mathbbl 1^{\otimes j_k})^{\oplus l_k}} & {\mathbbl 1^{\oplus i_{1,1}}\otimes \cdots \otimes \mathbbl 1^{\oplus i_{k,j_k}}}
	\arrow["{\kappa'}", from=1-1, to=1-2]
	\arrow["{\kappa^{\oplus l_1\cdots l_k}}"', from=1-1, to=2-2]
	\arrow["{\tilde d'}", from=1-2, to=1-3]
	\arrow["{(\kappa^{\otimes k })^{\oplus l_1\cdots l_k}}", from=1-2, to=2-2]
	\arrow["{\kappa_1^{\oplus l_1}\otimes \cdots \otimes \kappa_k^{\oplus l_k}}", from=1-3, to=2-3]
	\arrow["{\tilde d'}"', from=2-2, to=2-3]
	\arrow["{\tilde d}"', curve={height=18pt}, from=2-2, to=2-4]
	\arrow["{\tilde d_1 \otimes \cdots \otimes \tilde d_k}", from=2-3, to=2-4]
\end{tikzcd}\] where the lower composition is $\tilde d$ since each the maps in the composition are cocartesian lifts of distibutive maps in $\Fint$ and we have composition law in \autoref{fin2}. The middle square commutes by \autoref{functoraldist}.

Now we need to show (v). The first equivariant diagram now translates into 

\[\begin{tikzcd}[nodes={scale=0.8}]
	{\Map^{\beta^\times_k}(\mathbbl 1^{\otimes k},\mathbbl 1) \times \prod_{r=1}^k\Map^{\beta^\times_{j_r}}(\mathbbl 1 ^{\oplus j_r},\mathbbl 1) } & {\Map^{\beta^\times_k}(\mathbbl 1^{\otimes k},\mathbbl 1) \times \prod_{r=1}^k\Map^{\beta^\times_{j_{\sigma(r)}}}(\mathbbl 1 ^{\oplus j_{\sigma (r)}},\mathbbl 1) } \\
	{\Map^{\beta_{j_1}^+ \boxtimes\cdots \boxtimes \beta_{j_k}^+}(\mathbbl 1^{\oplus j_1} \otimes \cdots \otimes \mathbbl 1^{\oplus j_k},\mathbbl 1)} & {\Map^{\beta_{j_{\sigma (1)}}^+ \boxtimes\cdots \boxtimes \beta_{j_{\sigma(k)}}^+}(\mathbbl 1^{\oplus j_{\sigma (1)}} \otimes \cdots \otimes \mathbbl 1^{\oplus j_{\sigma (k)}},\mathbbl 1^{\otimes k})} \\
	{\Map^{\beta^+_{|Q|} \circ \beta^\times_k}((\mathbbl 1^{\otimes k})^{\oplus |Q|},\mathbbl 1 )} & {\Map^{\beta^+_{|Q|} \circ \beta^\times_k}((\mathbbl 1^{\otimes k})^{\oplus |Q|},\mathbbl 1 )} \\
	{\Map^{\beta^+_{|Q|}}(\mathbbl 1^{\oplus |Q|},\mathbbl 1 )} & {\Map^{\beta^+_{|Q|}}(\mathbbl 1^{\oplus |Q|},\mathbbl 1 )}
	\arrow["{\sigma^{-1} \times \sigma}", from=1-1, to=1-2]
	\arrow[from=1-1, to=2-1]
	\arrow[from=1-2, to=2-2]
	\arrow["\sigma", from=2-1, to=2-2]
	\arrow["{\tilde d}"', from=2-1, to=3-1]
	\arrow["{\tilde d}", from=2-2, to=3-2]
	\arrow["{\tilde \sigma}", from=3-1, to=3-2]
	\arrow["\kappa"', from=3-1, to=4-1]
	\arrow["\kappa", from=3-2, to=4-2]
	\arrow["{\tilde \sigma}", from=4-1, to=4-2]
\end{tikzcd}\] The first and the last square commute since $\sigma $ matches each components. The middle square is by the first commutative square in [\autoref{fin2},(iv)] and everything is cocartesian lift. Therefore, the outer square commutes and we've shown the first equivariant square. Applying the second commutative square in [\autoref{fin2},(iv)] similarly gives the second equivariant square.

Finally, we show the distributivity in (vi). We let $f \in \mathscr J (k) = \Map^{\beta^\times_k} (\mathbbl{1}^{\otimes k}, \mathbbl{1})$, $g_r \in \mathscr C (j_r) =\Map^{\beta^+_{j_r}} (\mathbbl{1}^{\oplus j_r}, \mathbbl{1})$ and $h_{r,q} \in \mathscr C(i_{r,q}) = \Map^{\beta^+_{i_{r,q}}} (\mathbbl{1}^{\oplus i_{r,q}}, \mathbbl{1})$. We also use $\tilde d,\tilde d_1, \tilde d_2$ and $\tilde d_\alpha$ for $\alpha \in \underline{Q}$ to denote the cocartesian lifting of distribute maps

\begin{equation*}
    \begin{aligned}
        &\tilde d \colon (\Un^{\otimes k})^{\oplus I} \to (\Un^{\oplus l_1}\otimes \cdots \otimes \Un ^{\oplus l_k}) \\ 
        &\tilde d_1 \colon \Un^{\oplus i_{1,q^1_1}} \otimes \cdots \otimes \Un^{\oplus i_{k,q^1_k}} \oplus \cdots \oplus \Un^{\oplus i_{1,q^Q_1}} \otimes \cdots \otimes \Un^{\oplus i_{k,q^Q_k}} \to \Un^{\oplus l_1}  \otimes \cdots \otimes  \Un^{\oplus l_k} \\
        &\tilde d_2 \colon (\Un^{\otimes k})^{\oplus |Q|} \to \Un^{\oplus j_1} \otimes \cdots \otimes \Un^{\oplus j_k}\\
        &\tilde d_\alpha \colon (\Un^{\otimes k})^{\oplus i_{1,q_1^\alpha} +\cdots + i_{k,q_k^\alpha}}\to  \Un^{\oplus i_{1,q_1^\alpha}} \otimes \cdots \otimes \Un^{\oplus i_{k,q_k^\alpha}}
    \end{aligned}
\end{equation*} Now we translate the distributivity diagram in \autoref{classicaldefn}. Clockwise it sends $(f,\{g_r\},\{h_{r,q}\})$ to 
\[
F_1 = f \cdot \big( g_1\cdot (h_{1,1}\oplus\cdots \oplus h_{1,j_1})\otimes \cdots \otimes g_k\cdot(h_{k,1}\oplus \cdots \oplus h_{k,j_k})\big) \tilde d \kappa 
\] while counterclockwise to 
\[
F_2 = f\cdot (g_1 \otimes \cdots \otimes g_k) \cdot\tilde d_2\cdot \kappa\cdot \bigoplus_{\alpha \in \underline{Q}}  f (h_{1,q_1^\alpha}\otimes \cdots \otimes h_{k,q^\alpha_k})\big)  \cdot \bigoplus_{\alpha \in \underline{Q}} \tilde d_\alpha \cdot \kappa
\] we need to show that they are identified. First, since $\tilde d = \tilde d_1 \cdot\bigoplus_{\alpha \in \underline{Q}} \tilde d_{\alpha}$, $F_1$ is identified with 
\[
F_1 =  f \cdot \big( g_1\cdot (h_{1,1}\oplus\cdots \oplus h_{1,j_1})\otimes \cdots \otimes g_k\cdot(h_{k,1}\oplus \cdots \oplus h_{k,j_k})\big) \tilde d_1 \cdot\bigoplus_{\alpha \in \underline{Q}} \tilde d_{\alpha}\tilde \kappa 
\] Consider the following diagram:

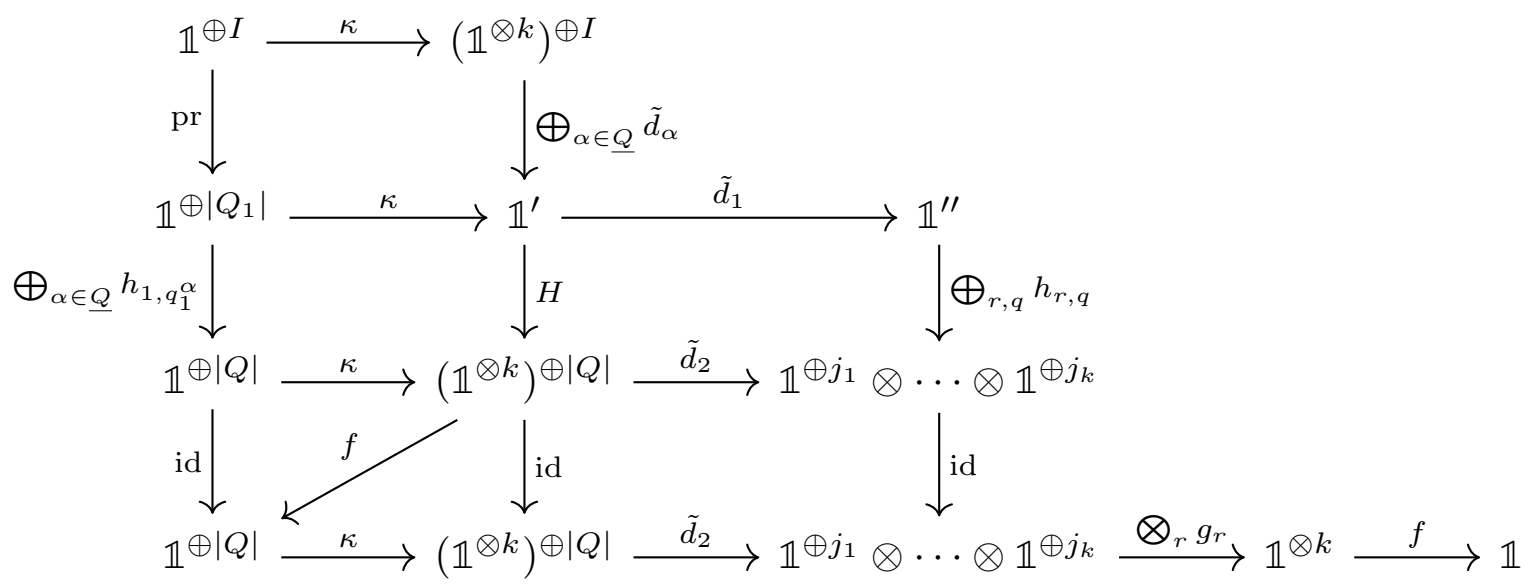
\[\begin{tikzcd}[nodes={scale=0.8}]
	{\mathbbl{1}^{\oplus I}} & {(\mathbbl{1}^{\otimes k})^{\oplus I}} &&& \\
	{\mathbbl{1}^{\oplus |Q_1|}} & {\mathbbl{1}'} & {\mathbbl{1}''} \\
	{\mathbbl{1}^{\oplus |Q|}} & {(\mathbbl{1}^{\otimes k})^{\oplus |Q|}} & {\mathbbl{1}^{\oplus j_1}\otimes \cdots \otimes \mathbbl{1}^{\oplus j_k}} \\
	{\mathbbl{1}^{\oplus |Q|}} & {(\mathbbl{1}^{\otimes k})^{\oplus |Q|}} & {\mathbbl{1}^{\oplus j_1}\otimes \cdots \otimes \mathbbl{1}^{\oplus j_k}} & {\mathbbl{1}^{\otimes k}} & {\mathbbl{1}}
	\arrow["\kappa", from=1-1, to=1-2]
	\arrow["{\mathrm{pr}}"', from=1-1, to=2-1]
	\arrow["{\bigoplus_{\alpha \in \underline{Q}} \tilde d_{\alpha}}", from=1-2, to=2-2]
	\arrow["\kappa", from=2-1, to=2-2]
	\arrow["{\bigoplus_{\alpha \in \underline{Q}}h_{1,q^\alpha_1}}"', from=2-1, to=3-1]
	\arrow["{\tilde d_1}", from=2-2, to=2-3]
	\arrow["H", from=2-2, to=3-2]
	\arrow["{\bigoplus_{r,q}h_{r,q}}", from=2-3, to=3-3]
	\arrow["\kappa", from=3-1, to=3-2]
	\arrow["{\mathrm{id}}"', from=3-1, to=4-1]
	\arrow["{\tilde{d}_2}", from=3-2, to=3-3]
	\arrow["f"', from=3-2, to=4-1]
	\arrow["{\mathrm{id}}", from=3-2, to=4-2]
	\arrow["{\mathrm{id}}", from=3-3, to=4-3]
	\arrow["\kappa", from=4-1, to=4-2]
	\arrow["{\tilde{d}_2}", from=4-2, to=4-3]
	\arrow["{\bigotimes_r g_r}", from=4-3, to=4-4]
	\arrow["f", from=4-4, to=4-5]
\end{tikzcd}\] Here we denote 

\begin{equation*}
    \begin{aligned}
        & \mathbbl 1^{\oplus |Q_1|} = \mathbbl 1^{\oplus i_{1,q_1^1}} \oplus \cdots \oplus \mathbbl 1^{\oplus i_{1,q_1^Q}} \\
        & \mathbbl 1 ' =\Un^{\oplus i_{1,q^1_1}} \otimes \cdots \otimes \Un^{\oplus i_{k,q^1_k}} \oplus \cdots \oplus \Un^{\oplus i_{1,q^{|Q|}_1}} \otimes \cdots \otimes \Un^{\oplus i_{k,q^{|Q|}_k}} \\
        & \mathbbl 1''= \Un^{\oplus l_1}  \otimes \cdots \otimes  \Un^{\oplus l_k} \\
        & H=\bigoplus_{\alpha \in \underline{Q}}   (h_{1,q_1^\alpha}\otimes \cdots \otimes h_{k,q^\alpha_k})
    \end{aligned}
\end{equation*} Every diagram except the lower triangle commutes by \autoref{unitlemma}, \autoref{functoraldist} and \autoref{specialkappa}. The clockwise composition is $F_1$ and the counterclockwise is $F_2$. Thus they are identified.
\end{proof}

\section{Cocartesian-cartesian bimonoidal categories}\label{ccsection}

Let $\mathcal C$ be a classical category. Recall that $\mathcal C$ is called distributive if it admits finite products and finite coproducts, and the following distributivity law is satisfied:
\[
(A \prod B ) \coprod (A \prod C) \cong A \prod (B \coprod C) 
\] In this case, $\mathcal C$ together with $\prod $ and $\coprod$ consist a symmetric bimonoidal category. We call such a symmetric bimonoidal structure the \emph{distributive bimonoidal structure} 

In this subsection, we will formalize this intuition into our $\infty$-categorical settings. We begin with the definition of distributive bimonoidal structure.

\begin{defn}
    Let $\mathcal C$ be an $\infty$-category. We say that $\mathcal C$ is \emph{distributive} if:
    \begin{enumerate}
        \item $\mathcal{C}$ admits finite products and coproducts.
        \item For every object $X$, $Y$ and $Z$, the canonical morphism
        \[
        (X \prod Y ) \coprod (X \prod Z) \to X \prod (Y \coprod Z) 
        \] is an equivalence.
    \end{enumerate}
\end{defn}


Our goal is to show that for any distributive category, we can always extend it to an $\infty$-operad, making it into a symmetric bimonoidal category. To work out this, we follow the spirit of \cite[Chapter 2.4.1]{lurie2017higher}.

\begin{defn}\label{gamma2}
    We define a category $\Gamma_2$ equipped with a symmetric monoidal forgetful functor $\mathrm{fgt}\colon \Gamma_2 \to \Fint$ as follows:
\begin{enumerate}
    \item The objects of $\Gamma_2$ are maps $f\colon (\langle m \rangle , \gamma) \to \langle 1 \rangle$ in $\Fint$ such that $f$ is generated by $\phi^+$ and $\phi^\times$ for some $\phi \in \Finp$ under $\boxplus$, $\boxtimes$ and composition. That is, $f=f_1f_2\cdots f_k$ where each $f_i $ is of the form $f_i=g_{1,i} \bullet\cdots \bullet g_{r_i,i}$, such that $g_{j,i}$ is some $\phi^+$, $\phi^\times$. Note that there are two operations on the underlying set of object, namely $\widetilde {\boxplus}$ and $\widetilde{\boxtimes}$ defined by
    \[
    f_1 \widetilde{\boxplus} f_2 := \beta_2^+ \circ(f_1\boxplus f_2) \quad f_1 \widetilde{\boxtimes} f_2 := \beta_2^+ \circ(f_1\boxtimes f_2)
    \] Note that $\Gamma_2$ is generated by $\phi^+$ and $\phi^\times$ for some $\phi \in \Finp$ under $\widetilde{\boxplus}$ and $\widetilde{\boxtimes}$.
    \item The morphisms of $\Gamma_2$ is generated by the following classes under $\widetilde{\boxplus}$, $\widetilde{\boxtimes}$ and composition:
    \begin{enumerate}
        \item For $\phi_1^+ \colon (\langle m_1 \rangle ,0) \to \langle 1 \rangle$, $\phi_2^+ \colon (\langle m_2 \rangle ,0) \to \langle 1 \rangle$, there is a morphism $\alpha^+ \colon \phi_1^+ \to \phi_2^+$, which would be sent to the morphism with the same name $\alpha^+\colon (\langle m_1 \rangle ,0) \to (\langle m_2 \rangle ,0) \in \Fint$ through the forgetful functor, if in $\Finp$ we have $\alpha \phi_1^{-1} \{1\} \subset \phi_2^{-1}\{1\}$.
        \item For $\phi_1^\times \colon (\langle m_1 \rangle ,2^{m_1-1}-1) \to \langle 1 \rangle$, $\phi_2^\times \colon (\langle m_2 \rangle ,2^{m_2-1}-1) \to \langle 1 \rangle$, there is a morphism $\alpha \colon \phi_1^\times \to \phi_2^\times$, which would be sent to the morphism with the same name $\alpha^\times\colon (\langle m_1 \rangle ,2^{m_1-1}-1) \to (\langle m_2 \rangle ,2^{m_2-1}-1) \in \Fint$ throught the forgetful functor, if in $\Finp$ we have $\alpha^{-1} \phi_2^{-1} \{1\} \subset \phi_1^{-1}\{1\}$.
        \item Suppose we are given $Q \in \mathrm{Seq}(\mathbb N)$, $f_1 \colon A^Q_\times \to \langle 1 \rangle$ and $f_2 \colon A^Q_+ \to \langle 1 \rangle$. There is a morphism $d_Q \colon f_1 \to f_2$, which would be sent to $d_Q \colon A^Q_\times \to A^Q_+ \in \Fint$ through the forgetful functor, if the following diagram commutes in $\mathrm{N}\Fint$:
        \[\begin{tikzcd}
	{A^Q_\times} && {A^Q_+} \\
	& {(\langle1\rangle , 0)}
	\arrow["{d_Q}", from=1-1, to=1-3]
	\arrow["{f_1}"', from=1-1, to=2-2]
	\arrow["{f_2}", from=1-3, to=2-2]
\end{tikzcd}\] 
\end{enumerate}
    \end{enumerate}
    \end{defn}

Roughly, we can think of $\Gamma_2$ as a package, containing the data of methods to combine $m$ different objects through $\boxplus$ and $\boxtimes$.

\begin{rem}\label{smallpullbacklem}
    This definition of $\Gamma_2$ is the counterpart of Lurie's definition of $\Gamma^\times$ (\cite[Notation 2.4.1.2]{lurie2017higher}) in the operad pair settings. It is clear that the objects of Lurie's $\Gamma^\times$ can be written as $(\Finp)_{/\langle 1\rangle}$, where any map $\alpha \colon \langle n \rangle \to \langle 1 \rangle$ distinguishes a subset $S:= \alpha^{-1} \{1\}\subset\langle n \rangle^\circ$. Moreover, one can show that we have a pullback square:
    \[\begin{tikzcd}
	{\Gamma^\times} & {\Gamma_2} \\
	\Finp & \Fint
	\arrow[hook, from=1-1, to=1-2]
	\arrow["{\mathrm{fgt}}"', from=1-1, to=2-1]
	\arrow["{\mathrm{fgt}}", from=1-2, to=2-2]
	\arrow["{i_1}", hook, from=2-1, to=2-2]
	\arrow["\ulcorner"{anchor=center, pos=0.125, rotate=180}, draw=none, from=2-2, to=1-1]
\end{tikzcd}\]
\end{rem}

\begin{exmp}\label{examplesforgamma2}
    \begin{enumerate}
        \item Let $\mathrm{pt}\colon \langle 1 \rangle \to \langle 1 \rangle$ denote the base point map. Then in $\Gamma_2$, we have two objects $\beta_k^+ \colon (\langle k \rangle ,0 ) \to \langle 1 \rangle$ and $\rho^+_i\colon (\langle k \rangle ,0 ) \to \langle 1 \rangle$. If we denote $\mathrm{id}\colon \langle 1 \rangle \to \langle 1 \rangle$ as the identity. Then we have $\beta_k^+ \simeq \mathrm{id}^{\widetilde{\boxplus} k}$ and $\rho_i \simeq \mathrm{pt} \widetilde{\boxplus}\cdots\widetilde{\boxplus} \mathrm{id} \widetilde {\boxplus} \cdots \widetilde{\boxplus} \mathrm{pt}$ where $\mathrm{id}$ is taken over the i-th slot. The identity morphism on $(\langle m \rangle,0)$ then induces a morphism $\mathrm{id}^+_i\colon \rho_i^+ \to \beta_k^+$ in $\Gamma_2$.
        \item Dually, for two objects $\beta_k^\times \colon (\langle k \rangle ,2^{k-1}-1 ) \to \langle 1 \rangle$ and $\rho^\times_i\colon (\langle k \rangle ,2^{k-1}-1 ) \to \langle 1 \rangle$. The identity morphism induces a morphism $\mathrm{id}_i^\times \colon \beta_k^\times \to \rho_i^\times$ in $\Gamma_2$.
    \end{enumerate}
\end{exmp}

We can not expect $\mathrm{fgt} \colon \Gamma_2 \to \Fint$ to be a Grothendieck fibration as in mono-operadic case. However, there are indeed enough cartesian edges.

\begin{lem}\label{cartesianforgamma2}
    Let $\mathrm{fgt}\colon \mathrm{N}\Gamma_2 \to \mathrm{N}\Fint$ be the forgetful functor. 
    \begin{enumerate}
        \item Let $\alpha\in \mathrm{N}\Fint$ be a non-distributive map, that is, $\alpha$ is generated by $f^+, f^\times, \sigma^+$ and $\sigma^\times$. Then for every $\psi\colon B\to \langle 1 \rangle \in \mathrm{N}\Fint$, there is a cartesian edge $\phi\to \psi$ lifting $\alpha$.
        \item Let $Q=(k,j_1,\cdots,j_k)\in\mathrm{Seq}(\mathbb N)$. Let $d_Q\colon A_\times^Q \to A_+^Q$ be the distributive morphism in $\mathrm{N}\Fint$. Suppose that we have $\psi \colon A_+^Q \to \langle 1 \rangle \in \mathrm{N}\Gamma_2$ such that $\phi= \psi \circ d_Q$ is also in $\mathrm{N} \Gamma_2$. Then there is a cartesian edge $\phi \to \psi$ lifting $d_Q$.
    \end{enumerate}
\end{lem}
\begin{proof}
    It suffices to note that we can take $\phi = \psi \circ\alpha$ and $\phi= \psi \circ d_Q$ in each case respectively. Note that in each case $\phi$ locates in $\mathrm{N}\Gamma_2$ and the morphism of the same name gives the desired cartesian lift.
\end{proof}

\begin{con}\label{constructionofcocartcart}
    Let $\mathcal C$ be an $\infty$-category. We construct a simplicial set $\widetilde {\mathcal {C}}^{\sqcup.\times}$ together with a map $q\colon \widetilde {\mathcal {C}}^{\sqcup.\times} \to \mathrm{N}\Fint$ by the following universal property: for every map of simplicial sets $K \to \mathrm{N}\Fint$, we have a bijection
    \[
    \Hom_{\mathrm{N}\Fint} (K, \widetilde {\mathcal {C}}^{\sqcup.\times}) \simeq \Hom_{\sSet} (K \times_{\mathrm{N}\Fint} \mathrm{N}\Gamma_2 , \mathcal C) 
    \] It is clear that for any $(\langle m \rangle , \gamma ) \in \Fint$, $\widetilde {\mathcal {C}}^{\sqcup.\times}_{(\langle m \rangle , \gamma)}$ is identified with $\Fun ((\mathrm {N}\Gamma_2)_{(\langle m \rangle , \gamma)} , \mathcal C)$. We let $\mathcal{C}^{\sqcup, \times}$ denote the full subcategory of $\widetilde {\mathcal {C}}^{\sqcup.\times}$ generated by functors $f$ that satisfy:
    \begin{enumerate}[(i)]
        \item Let $\alpha_1 \colon (\langle m \rangle, \gamma) \to (\langle k \rangle,0 )\in \Fint$. Then by \autoref{examplesforgamma2}, we have morphisms $\mathrm{id}_i^+ \circ \mathrm{id} \colon \rho_i^+ \circ \alpha_1 \to \beta_k^+ \circ \alpha_2$. We require these morphisms together induce $f (\beta_k^+ \circ \alpha_1) \simeq \coprod_{i=1}^k f(\rho_i^\times \circ \alpha_1)$.
        \item Let $\alpha_2 \colon (\langle m \rangle, \gamma) \to (\langle k \rangle,2^{k-1}-1 )\in \Fint$. Similarly we have morphisms $\mathrm{id}_i^\times \circ \mathrm{id} \colon \rho_i^\times \circ \alpha_1 \to \beta_k^\times \circ \alpha_2$. We require these morphisms together induce $f (\beta_k^\times \circ \alpha) \simeq \prod_{i=1}^k f(\rho_i^\times \circ \alpha)$.
        \item Consider the following distributive square:
\[\begin{tikzcd}
	{\bigboxplus_{Q} \bigboxtimes_{1\leq r \leq k} \langle 1 \rangle} & {\bigboxtimes_{1\leq r \leq k} \bigboxplus_{1 \leq q \leq j_r} \langle 1 \rangle} \\
	& {(\langle m \rangle ,2^{m-1}-1)} \\
	{(\langle Q \rangle,0)} & {(\langle 1\rangle, 0)}
	\arrow["{d_Q}", from=1-1, to=1-2]
	\arrow["{\beta^\times_m\boxplus\cdots\boxplus \beta_m^\times}"', from=1-1, to=3-1]
	\arrow["{\beta_{j_1}^+ \boxtimes \cdots \boxtimes \beta_{j_m}^+}", from=1-2, to=2-2]
	\arrow["{\beta_m^\times}", from=2-2, to=3-2]
	\arrow["{\beta_{|Q|}^+}", from=3-1, to=3-2]
\end{tikzcd}\] By (i), (ii) and composition, $f$ takes $\beta_Q^+ \circ (\beta_m^\times \boxplus \cdots \boxplus \beta_m^\times)$ to $\coprod_Q \prod_{1 \leq r \leq k} A_{r,q_r}$ and $\beta_m^\times \circ (\beta_{j_1}^+ \boxtimes \cdots \boxtimes \beta_{j_m}^+)$ to $\prod_{1\leq r \leq k} \coprod_{1\leq q_ \leq j_r} A_{r,q}$, where $A_{r,q}$ is identified with $f(\rho^{r,q})$. Then when viewing as a morphism in $\mathrm{N}\Gamma_2$, $d$ is sent to the canonical comparison morphism.
    \end{enumerate}
\end{con}

\begin{rem}\label{supplimentforinitial}
    Taking $\alpha\colon \langle 1 \rangle \to \langle0\rangle$, the condition (i) implicitly requires that $f$ sends $\alpha$ to the initial object.
\end{rem}

\begin{prop}\label{cct} Let $\mathcal C$ be an $\infty$-category
    \begin{enumerate}
        \item Let $q \colon \widetilde {\mathcal {C}}^{\sqcup.\times} \to \mathrm{N}\Fint$ be the functor constructed above. For every inert or active morphism $\alpha\colon X \to Y\in \mathrm{N}\Fint$ and every $f \in \widetilde{\mathcal C}^{\sqcup.\times}_{X}$, $\alpha$ admits a cocartesian lift $\overline \alpha \colon f \to f'$. 
        \item Moreover, let notations be as in \autoref{fin2} and \autoref{exampleofmorphisms}. Suppose in addition that $\mathcal C$ is distributive. Then for every $f \in \mathcal{C}^{\sqcup,\times}_{A^Q_\times} \subset\widetilde{\mathcal C}^{\sqcup.\times}_{A^Q_\times}$ such that \[(\rho^{m,q_m,l_1}_\times)_{!}(f) \simeq (\rho^{m,q_m,l_2}_\times)_{!}(f)\] for any fixed $m,q_m$ and $0\leq l_1 \leq l_2 \leq n(m)$, there is a cocartesian lift $\tilde {d}_{Q} \colon f \to f'$ of $d_{Q}\colon A^Q_\times \to A^Q_+$ fitting into a lift of  the diagram in [\autoref{fin2},(viii)] such that $f'$ also lies in $\mathcal C^{\sqcup,\times}$.
        \item Let $\overline{\alpha} \colon f \to f'$ be a morphism in $\widetilde {\mathcal {C}}^{\sqcup.\times}$ that lifts $\alpha\colon A_1 \to A_2 \in \mathrm{N}\Fint$ such that $\alpha$ and $f$ satisfies conditions in (1) or (2). Then $\overline {\alpha}$ is cocartesian if and only if for every $g \colon A_2 \to \langle 1 \rangle$, the induced map $f(g\circ \alpha) \to f'(g)$ is an equivalence whenever $g\circ \alpha$ is defined.
        \item The functor $q$ restricts to a projection $q'\colon {\mathcal C}^{\sqcup.\times} \to \mathrm{N}\Fint$ with the same cocartesian edges. 
        \item The projection $q'\colon {\mathcal C}^{\sqcup.\times} \to \mathrm{N}\Fint$ becomes a symmetric bimonoidal category if and only if $\mathcal C$ is distributive.
    \end{enumerate}
\end{prop}

\begin{proof}
    We first prove (3). We can not apply \cite[Proposition 2.1.4.5]{lurie2017higher} here since both the forgetful functor fails to be cartesian and the conclusion is incorrect. Therefore, we choose to give a direct proof. Let $\alpha\colon A \to B \in \Fint$. We say an edge $\overline{\alpha}$ is special if it satisfies the condition of (3). We need to show that an edge is a cocartesian edge if and only if it is special. For any edge $e \in \widetilde {\mathcal {C}}^{\sqcup.\times}$, consider the following lifting problem
\[\begin{tikzcd}
	{\Lambda^n_0} & {\widetilde {\mathcal {C}}^{\sqcup.\times}} \\
	{\Delta^n} & {\mathrm{N}\Fint}
	\arrow["{\theta_0}", from=1-1, to=1-2]
	\arrow[hook, from=1-1, to=2-1]
	\arrow["q", from=1-2, to=2-2]
	\arrow[dashed, from=2-1, to=1-2]
	\arrow["\theta", from=2-1, to=2-2]
\end{tikzcd}\] where $\theta_0 |_{\Delta^{\{0,1\}}} = e$. Using the construction of $\widetilde {\mathcal {C}}^{\sqcup.\times}$, this is equivalent to the lifting problem 
\[\begin{tikzcd}
	{\mathrm{N}\Gamma_2\times_{\mathrm{N}\Fint}\Lambda^n_0} & {\mathcal C} \\
	{\mathrm{N}\Gamma_2\times_{\mathrm{N}\Fint}\Delta^n}
	\arrow[from=1-1, to=1-2]
	\arrow[hook, from=1-1, to=2-1]
	\arrow[dashed, from=2-1, to=1-2]
\end{tikzcd}\] Unwinding the definition, this is to say the cartesian edges in $\mathrm{N}\Gamma_2 \times_{\mathrm{N}\Fint} \Delta^{\{0,1\}}$ lifting $\alpha$ is sent to equivalences. Let $f,f' \in \widetilde {\mathcal {C}}^{\sqcup.\times}$ corresponding to $f \colon \mathrm{N}\Gamma_2 \times_{\mathrm{N}\Fint} \{A\} \to \mathcal C$ and $f'\colon\mathrm{N}\Gamma_2 \times_{\mathrm{N}\Fint} \{B\} \to \mathcal C$. Then this is equivalent to say that we have $f(g \circ \alpha) \simeq f'(g)$, here $\alpha \colon g\circ\alpha \to g$ is the cartesian lift given by \autoref{cartesianforgamma2}. This proves the claim (3). 

By \autoref{cartesianforgamma2}, one can always define $f'(g)= f(g\circ \alpha)$ in this way whenever the condition (1) is satisfied, this proves (1). Suppose additionally that $\mathcal C$ is distributive and (2) is satisfied. Then for $f\colon \mathrm{N}\Gamma_2 \times_{\mathrm{N}\Fint} A_\times^Q\to \mathcal C$ that satisfies the condition, we need to find $f'\colon \mathrm{N}\Gamma_2\times_{\mathrm{N}\Fint} A^Q_+ \to \mathcal C$ such that $f'(\phi)= f(\phi\circ d_Q)$ whenever $\phi\circ d_Q$ is defined. By the requirement (i) and (ii) of \autoref{constructionofcocartcart} and distributivity of $\mathcal C$, for any $\phi$, we can write $f'(\phi)=\coprod_I \prod_{J_i} f'(\rho_{i,j})$ for some indices $I,J$ and some projections $\rho_{i,j}$. Hence, we define $$f'(\phi):= \coprod_I\prod_{J_i} f(\rho_{i,j}\circ d_Q) \simeq \coprod_I\prod_{J_i} f(\rho_{i,j}') \in \mathcal C$$ When $\phi\circ d_Q$ is defined, this coincides with $f(\phi\circ d_Q)$ by condition (iii) in \autoref{constructionofcocartcart} and the assumption for $f$, proving the claim (3). (4) follows immediately since $\mathcal C^{\sqcup,\times}$ contains all these cocartesian edges. 
    
    We now prove (5). If $q'$ is a symmetric bimonoidal category, it is clear that the bimonoidal structure is given by product and coproduct and the distributive law guarantees the distributivity. We need only to show the converse. It is clear from \autoref{supplimentforinitial} that we have $\mathcal C^{\sqcup,\times}_{\langle 1 \rangle} \simeq \mathcal C$. For any $A \in \Fint$ with length $n$, let $\phi \colon \mathcal C^{\sqcup,\times}_{A} \to \mathcal C^{\times n}$ be the map induced by $\rho_i$. By \autoref{settingcocartandcart}, it suffices to show that $\phi$ induces an equivalence. 
    
    Let $(\mathrm N \Gamma_2)_A$ denote the fiber $\{A\}\times _{\mathrm{N}\Fint} \mathrm{N}\Gamma_2$, $(\mathrm N \Gamma_2)_A'\subset (\mathrm N \Gamma_2)_A$ denote the full subcategory consisting of those maps generated by $\alpha^\times$ for $\alpha \in \Finp$, and  let $(\mathrm N \Gamma_2)_A'' \subset (\mathrm N \Gamma_2)_A'$ denote the full subcategory consisting of $\rho_i \colon A \to \langle 1 \rangle$. For any $f \in (\mathrm N \Gamma_2)_A$, by \autoref{gamma2}, $f$ can be written as $f_1 f_2\cdots f_k$. Each $f_i$ is of the form $g_{1,i}\bullet \cdots \bullet g_{r_i,i}$ where each $g_{j,i}$ is $\alpha^+$ or  $\alpha^\times$. We can always choose $f_1$ such that it factors as $f_1=\beta_s^+h_1$ $f_1= \beta_s^\times h_1$ for some active $\beta_s$ since the target is $\langle 1 \rangle$. By definition, for every $F \in \mathcal C ^{\sqcup,\times}_A$, we have 
    \[
    F (f) \simeq \begin{cases}
       \coprod_{1\leq i \leq s} F(\rho_i h_1f_2\cdots f_k) & \text {if } f_1 = \beta_s^+ h_1\\
       \prod_{1\leq i \leq s} F(\rho_i h_1f_2\cdots f_k) & \text {if } f_1 = \beta_s^\times h_1
    \end{cases}
    \] In both cases we can continue with $f'= \rho_i h_1 f_2 \cdots f_k$. Consequently, combining with the distributivity of $\mathcal C$, $F$ is the double Kan-extension by its restricting on $(\mathrm N\Gamma_2)_A''$, meaning $\mathcal C^{\cup,\times}_A \simeq \mathcal C^{\times n}$.

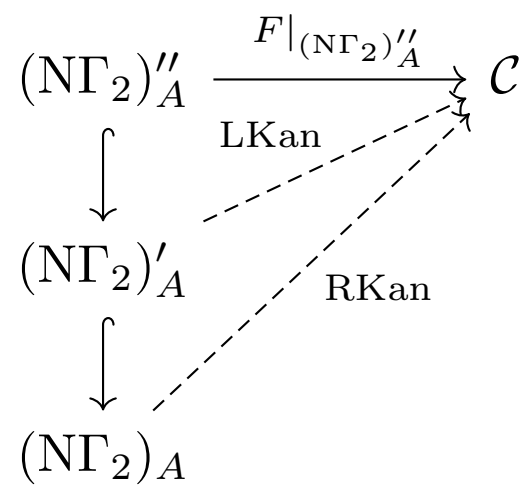
    \[\begin{tikzcd}
	{(\mathrm{N}\Gamma_2)_A''} && {\mathcal C} \\
	{(\mathrm{N}\Gamma_2)_A'} \\
	{(\mathrm{N}\Gamma_2)_A}
	\arrow["{F|_{(\mathrm{N}\Gamma_2)_A''}}", from=1-1, to=1-3]
	\arrow[hook, from=1-1, to=2-1]
	\arrow["{\mathrm{LKan}}", dashed, from=2-1, to=1-3]
	\arrow[hook, from=2-1, to=3-1]
	\arrow["{\mathrm{RKan}}"', dashed, from=3-1, to=1-3]
\end{tikzcd}\] \qedhere
\end{proof}

Let $\mathcal C$ be a distributive category, we call $q' \colon \mathcal C^{\sqcup, \times} \to \mathrm{N}\Fint$ the \emph{cocartesian-cartesian bimonoidal structure}. Recall that in \cite[Proposition 2.4.3.3]{lurie2017higher} and \cite[Proposition 2.4.1.5]{lurie2017higher}, Lurie defines cocartesian and cartesian monoidal operads. It relates to our bimonoidal structure in the following way:

\begin{cor}\label{undercct}
    Suppose that $\mathcal C$ is distributive. Let $q' \colon \mathcal C^{\sqcup, \times} \to \mathrm{N}\Fint$ be the cocartesian-cartesian bimonoidal structure. Then:
    \begin{enumerate}
        \item The underlying operad of $\mathcal C^{\sqcup, \times}$ is the cocartesian operad $\mathcal C^{\sqcup}$.
        \item The underlying moperad of $\mathcal C^{\sqcup, \times}$ is the cartesian operad $\mathcal C^\times$.
    \end{enumerate}
\end{cor}

\begin{proof}
    Let $i_1\colon \mathrm{N}\Finp \to \mathrm{N}\Fint$ be the second inclusion. By \autoref{smallpullbacklem}, pulling back along $i_0$ gives a functor $q'_1 \colon \tilde {\mathcal C}^{\times} \to \mathrm{N}\Finp$, which satisfies
    \[
    \Hom_{\mathrm{N}\Finp} (K, \widetilde {\mathcal {C}}^{\times}) \simeq \Hom_{\sSet} (K \times_{\mathrm{N}\Finp} \mathrm{N}\Gamma^\times , \mathcal C) 
    \] This coincides with the definition of Lurie's $\tilde{\mathcal C}^\times$ and the subcategory $\mathcal C^{\sqcup,\times} $ pulls back to $\mathcal C^\times$ by construction, proving (2). For (1), let $i_0\colon \mathrm{N}\Finp \to \mathrm{N}\Fint$ be the first inclusion and let $q_0'\colon \tilde{\mathcal C}^{\sqcup} \to \mathrm{N} \Finp$ be the pullback of $q'$ along $i_0$. We also have \[
    \Hom_{\mathrm{N}\Finp} (K, \widetilde {\mathcal {C}}^{\sqcup}) \simeq \Hom_{\sSet} (K \times_{\mathrm{N}\Finp} \mathrm{N}\Gamma^+ , \mathcal C) \]Here objects of $\Gamma^+$ are pairs $(\langle m\rangle,S)$ such that $\langle m \rangle \in \Finp$ and $S\subset \langle m \rangle^\circ$. A morphism $f\colon (\langle m \rangle , S) \to (\langle n \rangle ,S')$ is a morphism $f\colon \langle m \rangle \to \langle n \rangle $ such that $f(S) \subset S'$. Let $\mathcal C^{\sqcup} \subset \tilde{\mathcal C}^{\sqcup}$ be the subcategory pulled back from $\mathcal C^{\sqcup,\times}$. Then by construction, $(\mathcal C^{\sqcup})_{\langle m \rangle}$ is generated by functors $f\colon (\mathrm{N}\Gamma^+)_{\langle m \rangle} \to \mathcal C$ such that $f(S)\simeq \coprod_{i\in S} f(\{i\})$. These functors are exactly left Kan-extended from $\mathrm N\Gamma^*\subset \mathrm N \Gamma^+$, where $\Gamma^*$ is defined in \cite[Construction 2.4.3.1]{lurie2017higher}. Therefore one obtains

    \[
    \Hom_{\mathrm{N}\Finp} (K,  {\mathcal {C}}^{\sqcup}) \simeq \Hom_{\sSet} (K \times_{\mathrm{N}\Finp} \mathrm{N}\Gamma^* , \mathcal C) \] which is exactly the definition of Lurie's cocartesian operad. \qedhere
\end{proof}

\begin{rem}\label{dualofcocartcart}
    The \autoref{constructionofcocartcart} has a dualized version. We say a category $\mathcal C$ is \emph{codistributive} if $\mathcal C^{\mathrm{op}}$ is distributive. Then for any codistributive category we can analogously define the \emph{cartesian-cocartesian} bimonoidal structure $\mathcal C^{\times,\sqcup}\to \mathrm{N}\Fint$ by exchanging products and coproducts in \autoref{constructionofcocartcart}. Then the underlying operad and moperad of $\mathcal C^{\times, \sqcup} $ are the cartesian and cocartesian operad respectively.
\end{rem}



\section{Algebras and algebra pairs}\label{algsection}
Let $\opp$ be an $\infty$-operad pair with the underlying operad $\mathcal O^{\otimes}_0$ and the underlying moperad $\mathcal O^{\otimes}_1$. Let $\mathcal C^{\oplus, \otimes} \to \opp$ be a $\opp$-bimonoidal category as defined in \autoref{bimonoidal}. When considering algebras in $\opp$-bimonoidal categories, we have two cases to distinguish: 

\begin{enumerate}[(1)]
  \item Let $\mathcal O^\otimes$ be an operad over the underlying operad or moperad of $\mathcal C^{\oplus, \otimes}$. We can consider an $\mathcal O^{\otimes}$-algebra in $\mathcal C^{\oplus, \otimes}$. That means, we forget the operad pair to an operad and then consider its mono-algbebra.
  \item Let $\mathcal D^{\oplus, \otimes}$ be an operad pair over the $\mathcal C^{\oplus, \otimes}$. We can consider a $\mathcal D^{\oplus, \otimes}$-algebra in $\mathcal C^{\oplus, \otimes}$. That means, we can consider an operad pair acting on another operad pair.
\end{enumerate}
 
These two cases must interact well with each other. Suppose we only consider the ``mono-algebra" of a bimonoidal category, then it is natural to expect that the category consisting of these mono-algebra will inherit a monoidal structure. To begin with, we clarify the definitions. 

\begin{defn}
    Let $\opp$ be an $\infty$-operad pair with the underlying operad $\mathcal O_0^\otimes$ and the underlying moperad $\mathcal O_1^\otimes$. Let $p\colon \mathcal C^{\oplus, \otimes} \to \opp$ be a $\opp$-bimonoidal category as defined in \autoref{bimonoidal}. Let $\mathcal C_0^\otimes$ and $\mathcal C_1^\otimes$ be the underlying operad and moperad of $\mathcal C^{\oplus, \otimes}$. Then:
    \begin{enumerate}[(1)]
        \item Let $q\colon\mathcal O^{\otimes} \to \mathcal O_1^\otimes$ be an $\infty$-operad over $\mathcal O_1^\otimes$. An $\mathcal O^{\otimes}$-\emph{algebra} in $\mathcal C$ is an object in $\mathrm{Alg}_{\mathcal O/\mathcal O_0}(\mathcal C)$ as in \cite[Definition 2.1.3.1]{lurie2017higher}. That is, a map of operads $f$ over $\mathcal O_0^\otimes$:
        \[\begin{tikzcd}
	{\mathcal O^{\otimes}} && {\mathcal C^{\otimes}_0} \\
	& {\mathcal{O}_0^\otimes}
	\arrow["f", from=1-1, to=1-3]
	\arrow["q"', from=1-1, to=2-2]
	\arrow["{p_0}", from=1-3, to=2-2]
\end{tikzcd}\]

\item Let $q\colon\mathcal O^{\otimes} \to \mathcal O_1^\otimes$ be an $\infty$-operad over $\mathcal O_1^\otimes$. An $\mathcal O^{\otimes}$-\emph{malgebra} in $\mathcal C$ is an object in $\mathrm{Alg}_{\mathcal O/\mathcal O_1}(\mathcal C)$.        
    \end{enumerate}
\end{defn}

Here we use the term ``malgebra" to emphasize that it corresponds to the moperad. Since the definition only relies on the underlying operad and moperad structure, there is nothing new in this case and we can study $\mathcal O^\otimes$-algebras and malgebras just as in mono-operad case.

\begin{defn}
    Let $\opp$ be an $\infty$-operad pair. Let $p\colon \mathcal C^{\oplus, \otimes} \to \opp$ be a $\opp$-bimonoidal category. Let $q\colon \mathcal D^{\oplus, \otimes} \to \opp$ be an $\infty$-operad pair over $\opp$. A $\mathcal D^{\oplus, \otimes}$-\emph{algebra pair} is an operad pair map $f\colon \mathcal D^{\oplus, \otimes} \to \mathcal C^{\oplus, \otimes}$ over $\opp$. That is, an operad pair map such that we have the following commutative diagram:
    \[\begin{tikzcd}
	{\mathcal D^{\oplus,\otimes}} && {\mathcal C^{\oplus, \otimes}} \\
	& \opp
	\arrow["f", from=1-1, to=1-3]
	\arrow["q"', from=1-1, to=2-2]
	\arrow["p", from=1-3, to=2-2]
\end{tikzcd}\] We use $\mathrm{PAlg}_{\mathcal D/ \mathcal P} (\mathcal C) \subset \Fun_{\opp}(\mathcal D^{\oplus, \otimes}, \mathcal C^{\oplus, \otimes})$ to denote the full subcategory generated by $\mathcal D^{\oplus, \otimes}$-algebra pairs.
\end{defn}

\begin{exmp}
    \begin{enumerate}
        \item Let $\mathcal D^{\oplus, \otimes}\simeq \opp$ be the commutative operad pair $\Comm^{\oplus,\otimes}$. Then an object in $\mathrm{PAlg}_{\mathcal D/ \mathcal P} (\mathcal C)$ is identified with an object $A \in \mathcal C$ together with two commutative, associative and unitary operations $\mu \colon A\oplus A \to A$ and $m\colon A \otimes A \to A$. These two operations are compatible with the distributive law. In this case we write $\mathrm{PAlg}_{\mathcal D/ \mathcal P} (\mathcal C)= \mathrm{CRig} (\mathcal C)$, which is identified with the category of commutative rig objects in the symmetric rig category $\mathcal C$.
        \item Similarly, if we choose $\mathcal D^{\oplus,\otimes}$ to be the associative operad pair $\Ass^{\oplus,\otimes}$. Then the category of $\mathcal D^{\oplus, \otimes}$-objects are identified with the non-commutative rig objects $\mathrm{Rig}(\mathcal C)$.
    \end{enumerate}
\end{exmp}

We now restrict to commutative algebras. Suppose we have a symmetric bimonoidal category $(\mathcal C, \oplus, \otimes, 0,\mathbbl 1)$, we can simply take commutative algebras with respect to the symmetric monoidal structure $\oplus$ and ask whether $\mathrm{CAlg}_{\oplus}(\mathcal C)$ will inherit a symmetric monoidal structure by $\otimes$. Moreover, we expect that $\mathrm{CAlg} (\mathrm {CAlg}_{\oplus}(\mathcal C))$ can be identified with the category of rig objects in $\mathcal C$. We can ask the same question even when $\mathcal C$ is merely a mono-monoidal category. For example, when $\mathcal C$ is cartesian monoidal, we may expect the category of semirings can be identified with commutative algebra objects in commutative monoids. We will show this in the following section. Our next goal in this section is to show the case where $\mathcal C$ admits a symmetric bimonoidal structure.

\begin{defn}
    Let $Q = (k,j_1,\cdots,j_k) \in \mathrm{Seq}(\mathbb N)$.
    \begin{enumerate}
        \item We let $1^Q_\times \in \Fint$ denote the object $\bigboxplus_{|Q|} (\langle k\rangle,2^{k-1}-1)$.
        \item We let $1^Q_+ \in \Fint$ denote the object 
        $(\langle j_1\rangle,0) \boxtimes\cdots \boxtimes (\langle j_k\rangle,0)$.
    \end{enumerate}
\end{defn}

\begin{lem}\label{wreth}
    Let $Q_1=(k,j_1,\cdots,j_k)$ and $Q_2=(k,i_1,\cdots,i_k)$ be two objects in $\mathrm{Seq (\mathbb N)}$ of the same length. Let $\{f_r\colon \langle j_r \rangle \to \langle i_r\rangle \in \Finp\} $ be a collection of morphisms for $1\leq r \leq k$. There exists a unique map $f_1 \wr \cdots \wr f_k \colon 1^Q_\times \to 1^Q_+$ making the following diagram commute in $\Fint$.
    \[\begin{tikzcd}
	{1^{Q_1}_\times} & {1^{Q_1}_+} \\
	{1^{Q_2}_\times} & {1^{Q_2}_+}
	\arrow["{d_{Q_1}}", from=1-1, to=1-2]
	\arrow["{f_1\wr \cdots\wr f_k}"', dashed, from=1-1, to=2-1]
	\arrow["{f_1\boxplus\cdots\boxplus f_k}", from=1-2, to=2-2]
	\arrow["{d_{Q_2}}"', from=2-1, to=2-2]
\end{tikzcd}\]
\end{lem}

\begin{proof}
    By composition laws in \autoref{fin2}, it suffices to consider the collections $\{f_i\}$ such that $f_2,\cdots, f_k=\mathrm{id}$ and $f_1=\beta_k$. Let $Q'=(j_2,\cdots,j_k)$. We define 
    \[
    \beta_k \wr  \mathrm{id} \wr \cdots \wr \mathrm{id} \colon 1^{Q_1}_\times =1\boxtimes1^{Q'}_\times \boxplus \cdots \boxplus 1\boxtimes1^{Q'}_\times \to 1\boxtimes1^{Q'}_\times 
    \] as the composition:
    \[\begin{tikzcd}
	{1\boxtimes1^{Q'}_\times\boxplus \cdots \boxplus 1\boxtimes1^{Q'}_\times} & {(1\boxplus\cdots \boxplus1)\boxtimes1^{Q'}_{\times}} \\
	& {1\boxtimes1^{Q'}_\times}
	\arrow["d", from=1-1, to=1-2]
	\arrow["{\beta_k\wr\mathrm{id}\wr \cdots\wr\mathrm{id}}"'{pos=0}, from=1-1, to=2-2]
	\arrow["{\beta_k\boxtimes\mathrm{id}}", from=1-2, to=2-2]
\end{tikzcd}\] This is the only choice that makes the diagram in the statement commute by checking relations in \autoref{fin2}.
\end{proof}

\begin{defn}\label{fintv}
    We define a new subcategory $\Fintv \subset \Fint$ as follows:
    \begin{enumerate}
        \item The objects of $\Fintv$ are objects of the form
        \[(\langle m \rangle , 2^{m-1}-1) \boxplus \cdots \boxplus (\langle m \rangle , 2^{m-1}-1)\] These are exactly objects of the form $1^Q_\times$ for some $Q$.
        \item The morphisms from $(\langle m \rangle , 2^{m-1}-1)^{\boxplus n}$ to $(\langle k\rangle , 2^{k-1}-1 )^{\boxplus n} $ are of the form $(f^{\times})^{\boxplus n}$ for $f \colon \langle m \rangle \to \langle k \rangle \in \Finp$.
        \item Let $Q_1=(j_1,\cdots,j_k)$ and $Q_2=(i_1,\cdots,i_k)$ be two tuples of the same length. Let $\{f_r\colon \langle j_r \rangle \to \langle i_r\rangle \in \Finp\} $ be a collection of morphisms for $1\leq r \leq k$. Between objects $1^{Q_1}_\times$ and $1^{Q_2}_\times$, there exists a map $f_1\wr \cdots \wr f_k$. 
    \end{enumerate} The morphisms in $\Fintv$ are generated by above morphisms under compositions.
\end{defn}

We can view $\Fintv$ as a twist of $\Finp \times \Finp$. More precisely, there are two projection maps from $\Fintv$ to $\Finp$.

\begin{defn}
    Let $\Fintv$ be as in \autoref{fintv}.
    \begin{enumerate}
        \item We define $p^\otimes \colon \Fintv \to \Finp$ as follows. It sends $(\langle m \rangle , 2^{m-1}-1)^{\boxplus n}$ to $\langle m\rangle$, the morphisms $(f^+)^{\boxplus n}$ to $f \in \Finp$ and the morphisms in [\autoref{fintv};(3)] to identities.
        \item We define $p^{\oplus } \colon \Fintv \to \Finp$ as follows. It sends $(\langle m \rangle , 2^{m-1}-1)^{\boxplus n}$ to $\langle n\rangle$ and the morphisms in [\autoref{fintv};(2)] to identities. For morphisms in  [\autoref{fintv};(3)], $p^{\oplus}$ sends $f_1\wr \cdots\wr f_k \colon 1_\times^{Q_1} \to 1^{Q_2}_\times$ to $$f_1\times\cdots \times f_k\colon \langle j_1\cdots j_k\rangle \simeq \langle |Q_1|\rangle  \to \langle |Q_2|\rangle \simeq \langle i_1\cdots i_k \rangle$$.
        \end{enumerate}
\end{defn}

In \autoref{shapeoffintv}, we display some arrows in $\Fintv$.  The projection $p^\otimes$ will project this 2-dimensional diagram to the uppermost horizontal line and the projection $p^\oplus$ will project it to the leftmost vertical line.

\begin{figure}
    \centering
\[\begin{tikzcd}[nodes={scale=0.7}]
	{\langle 1 \rangle} & {\langle 1 \rangle \boxtimes \langle 1 \rangle} & {\langle 1 \rangle \boxtimes \langle 1 \rangle  \boxtimes \langle 1 \rangle} & \cdots \\
	{\langle 1 \rangle \boxplus \langle 1 \rangle} & {\langle 1 \rangle \boxtimes \langle 1 \rangle \boxplus \langle 1 \rangle \boxtimes \langle 1 \rangle} & {\langle 1 \rangle \boxtimes \langle 1 \rangle\boxtimes \langle 1 \rangle \boxplus \langle 1 \rangle \boxtimes \langle 1 \rangle \boxtimes \langle 1 \rangle} & \cdots \\
	{\langle 1 \rangle \boxplus \langle 1 \rangle \boxplus \langle 1\rangle} & {\langle 1 \rangle \boxtimes \langle 1 \rangle  \boxplus \langle 1 \rangle \boxtimes \langle 1 \rangle \boxplus \langle 1 \rangle \boxtimes \langle 1 \rangle} & \cdots & \cdots \\
	\cdots & \cdots & \cdots
	\arrow["{f^\times}"', from=1-2, to=1-1]
	\arrow["{g^\times}"', from=1-3, to=1-2]
	\arrow[from=1-4, to=1-3]
	\arrow["{h^+}"', from=2-1, to=1-1]
	\arrow["{h_{1}\wr h_2}"', from=2-2, to=1-2]
	\arrow["{f_1^\times \boxplus f_2^\times}"', from=2-2, to=2-1]
	\arrow["{h_1\wr h_2 \wr h_3}"', from=2-3, to=1-3]
	\arrow["{g_1^\times \boxplus g_2^\times}"', from=2-3, to=2-2]
	\arrow[from=2-4, to=1-4]
	\arrow[from=2-4, to=2-3]
	\arrow["{l^+}"', from=3-1, to=2-1]
	\arrow["{l_1 \wr l_2}"', from=3-2, to=2-2]
	\arrow["{f_1^\times \boxplus f_2^\times \boxplus f_3^\times}"', from=3-2, to=3-1]
	\arrow[from=3-3, to=2-3]
	\arrow[from=3-3, to=3-2]
	\arrow[from=3-4, to=2-4]
	\arrow[from=3-4, to=3-3]
	\arrow[from=4-1, to=3-1]
	\arrow[from=4-2, to=3-2]
	\arrow[from=4-2, to=4-1]
	\arrow[from=4-3, to=3-3]
	\arrow[from=4-3, to=4-2]
\end{tikzcd}\]
    \caption{Some arrows in $\Fintv$.}
    \label{shapeoffintv}
\end{figure}

\begin{con}
    Let $p\colon \opp \to \Fint$ be an $\infty$-operad. We let $\oppv$ denote the pullback:
    
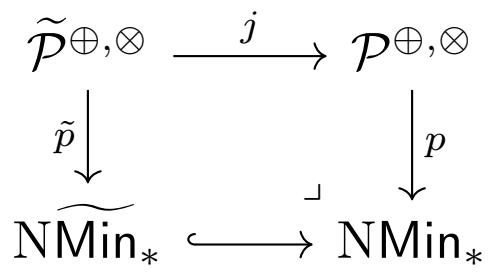
\[\begin{tikzcd}
	{\oppv} & \opp \\
	{\mathrm{N}\Fintv} & {\mathrm{N}\Fint}
	\arrow["j", from=1-1, to=1-2]
	\arrow["{\tilde{p}}"', from=1-1, to=2-1]
	\arrow["p", from=1-2, to=2-2]
	\arrow[hook, from=2-1, to=2-2]
	\arrow["\ulcorner"{anchor=center, pos=0.125, rotate=180}, draw=none, from=2-2, to=1-1]
\end{tikzcd}\] It is clear that the inclusion $\Finp \to \Fintv$ induces an inclusion $\mathcal O_0^\otimes \to \oppv$. Consider the composition $p':= p^{\oplus} \circ \tilde p\colon \oppv \to  \mathrm N \Finp$. We let 
\[\Fun' (\mathrm N \Finp, \oppv) \subset \Fun_{\mathrm N\Finp} (\mathrm N\Finp, \oppv)\] denote the full subcategory generated by the functors $F$ that satisfies the following conditions:
\begin{enumerate}
    \item $F$ sends inert morphisms to inert morphisms, that is, it sends inert maps in $\Finp$ to cocartesian lifts of inert maps in $\oppv$.
    \item For any $\langle k \rangle \in \Finp$, $p^{\otimes}\circ \tilde p \circ F(\langle k \rangle)$ is constant.
\end{enumerate}
We denote $\mathrm{CAlg}_{\oplus} (\mathcal P)^\otimes := \Fun' (\mathrm N \Finp, \oppv)$ and let $p'\colon \mathrm{CAlg}_{\oplus} (\mathcal P)^\otimes  \to \mathrm N \Finp$ denote the functor that sends $F $ to $p^{\otimes}\cdot \tilde p \cdot F(\langle 1 \rangle)$. 
\end{con}

\begin{lem}\label{auxilarylem}
    We say that a morphism in $\oppv$ is inert if its image in $\opp$ is inert. Then
    \begin{enumerate}
        \item $p^{\oplus}\colon \mathrm{N}\Fintv \to \mathrm{N}\Fint$ is a cocartesian fibration.
        \item A morphism $e \in \oppv$ is inert if and only if it is a cocartesian lift of inert maps in $\Finp$ through $p^{\oplus}$.
    \end{enumerate}
\end{lem}
\begin{proof}
    Given any object $(\langle m \rangle , 2^{m-1}-1)^{\boxplus n} \in \mathrm{N}\Fintv$, we can always write it as $1^Q_\times$ for $Q=(m,n,1\cdots ,1) \in \mathrm{Seq}(\mathbb N)$. So for any morphism $\alpha \colon \langle n \rangle \to \langle k \rangle \in \mathrm{N}\Finp$, we can choose the cocartesian lift as $$\overline{\alpha} = \alpha\wr\mathrm{id}\wr \cdots \wr \mathrm{id} \colon (\langle m \rangle , 2^{m-1}-1)^{\boxplus n} \to (\langle m \rangle , 2^{m-1}-1)^{\boxplus k}$$ By equivariant law in \autoref{fin2}, this choice is unique up to equivalence and becomes a cocartesian lift by construction. This proves (1). (2) follows from (1) and \cite[Proposition 2.4.1.3]{lurie2009higher}.
\end{proof}


\begin{thm}\label{calgtheorem} Let $p\colon \opp \to \mathrm N\Finp$ be an $\infty$-operad pair. Let $p'\colon \mathrm{CAlg}_{\oplus}(\mathcal P)^\otimes \to \mathrm N \Finp$ be constructed as above. Then:
\begin{enumerate}
    \item $p'\colon \mathrm{CAlg}_{\oplus}(\mathcal P)^\otimes \to \mathrm N \Finp$ exhibits $\mathrm{CAlg}_{\oplus}(\mathcal P)^\otimes$ as an $\infty$-operad. Moreover, if $\opp$ is in addition a symmetric bimonoidal category in the sense of \autoref{monoidalexample}, then $p'\colon \mathrm{CAlg}_{\oplus}(\mathcal P)^\otimes \to \mathrm N \Finp$ is a symmetric monoidal category with the underlying category identified with $\mathrm{CAlg} _\oplus(\mathcal P)$.
    \item We have a natural equivalence: 
    \[\mathrm {CAlg} (\mathrm {CAlg}_{\oplus}(\mathcal P)) \simeq \mathrm{CRig}(\mathcal P)\]
\end{enumerate}
\end{thm}

\begin{proof}
    We first show that the constructed morphism $p'\colon \mathrm{CAlg}_{\oplus}(\mathcal P)^{\otimes} \to \mathrm{N}\Finp$ is an $\infty$-operad. In order to show this, we note that any object or morphism in $\mathrm{CAlg}_{\oplus}(\mathcal P)^{\otimes}$ is fully determined by its evaluation on $\langle 1 \rangle \in \mathrm{N}\Finp$. More precisely, any $f \in \mathrm{CAlg}_{\oplus}(\mathcal P)^{\otimes}\simeq \Fun' (\mathrm N \Finp, \oppv)$ must send inert morphisms to inert morphisms and takes constant value after applying $p^{\otimes}$. Suppose we have $f,g \in\mathrm{CAlg}_{\oplus}(\mathcal P)^{\otimes} $ with $\alpha \colon f \to g$ such that $f(\langle 1 \rangle)=A \in \oppv, g(\langle 1 \rangle)=B \in \oppv$ and $\alpha_{\langle 1 \rangle} \colon A \to B$, then we must have $f(\langle n \rangle ) = A^{\oplus n}$, $g(\langle n \rangle ) = B^{\oplus n}$ and $\alpha_{\langle n \rangle}= \alpha_{\langle 1 \rangle}^{\oplus n}$. 

    For any inert morphism $t \colon \langle m \rangle \to \langle n \rangle\in \Finp$ and $f \in p'^{-1} (\langle m \rangle )$, let $A^{\oplus n} = f(\langle n\rangle)$. Since $p\colon \opp \to \mathrm N\Fint$ is an $\infty$-operad pair, we can find $B \in p'^{-1} (\langle n \rangle )$ and $\overline t \colon A \to B \in \opp$ the cocartesian lift of $t$. The morphisms $t^{\oplus n} \colon A^{\oplus n} \to B^{\oplus n}$ are also cocartesian, so these $t^{\oplus n}$ together extend to a cocartesian lift $\tilde t\colon f\to g $ in $\mathrm{CAlg}_{\oplus}(\mathcal P)^{\otimes}$. This proves (1) in \autoref{lurieoperad}. (2) and (3) follows immediately from the above argument and the operad pair structure on $\opp$. 

    When $\opp$ is itself a symmetric bimonoidal category, then any edge in $\Finp$ admits a cocartesian lift in $\oppv$. Above argument shows that $\mathrm{CAlg}_{\oplus} (\mathcal{P})^{\otimes} \to \mathrm N\Finp$ is a cocartesian fibration. Moreover, the underlying category of $\mathrm{CAlg}_{\oplus}(\mathcal P)^{\otimes}$ is identified with those $f \in \Fun' (\mathrm N \Finp, \oppv)$ that factors through $\mathcal O^\otimes_0$, which is exactly the category $\mathrm{CAlg}_{\oplus}(\mathcal P)$. These finish the proof of (1).

    We now show the equivalence $\mathrm {CAlg} (\mathrm {CAlg}_{\oplus}(\mathcal P)) \simeq \mathrm{CRig}(\mathcal P)$. Note that \[\mathrm {CAlg} (\mathrm {CAlg}_{\oplus}(\mathcal P))  \subset \Fun_{\mathrm N \Finp} (\mathrm N\Finp, \Fun' (\mathrm N\Finp, \oppv))\] is identified with the full subcategory  \[\Fun^{\mathrm{int}}_{\Fintv}(\mathrm{N}\Fintv , \oppv)\subset \Fun_{\Fintv}(\mathrm{N}\Fintv , \oppv)\] generated by functors that send inert morphisms to inert morphisms. We must show that there is an equivalence:
    \[
    \Fun^{\mathrm{int}}_{\mathrm{N}\Fintv}(\mathrm{N}\Fintv , \oppv) \simeq \Fun'_{\mathrm{N}\Fint}(\mathrm{N}\Fint , \opp)
    \] Here $\Fun'_{\mathrm{N}\Fint}(\mathrm{N}\Fint , \opp) \subset \Fun_{\mathrm{N}\Fint}(\mathrm{N}\Fint , \opp)$ is the full subcategory generated by those operad pair maps. It suffices to show that we can always extend $\tilde f \in \Fun^{\mathrm{int}}_{\mathrm{N}\Fintv}(\mathrm{N}\Fintv , \oppv)$ to a functor $f \in \Fun'_{\mathrm{N}\Fint}(\mathrm{N}\Fint , \opp)$. Suppose that $\tilde f(\langle 1 \rangle) = A$, then for any \[B = \langle 1 \rangle \bullet \cdots \bullet \langle 1 \rangle\in \Fint\] where each $ \bullet$ represent a $\boxplus$ or $\boxtimes$, we set $f(B)=A\odot \cdots \odot A$. The inert and distributive morphisms are sent to cocartesian lift of the corresponding maps in $\Fint$ by requirement. The active maps $\beta_m^+$ and $\beta_m^{\times}$ are already in $\Fintv$. This determines a well-defined functor $f\in \Fun'_{\Fint}(\mathrm{N}\Fint , \opp)$ since the morphisms in $\Fint$ are generated by these active maps,  inert maps and distributive maps. All the relations in \autoref{fin2} are preserved automatically.
\end{proof}

We end this subsection by a quick discussion of bimonoidal functors. Recall that for any $\infty$-operad $\mathcal O^\otimes$ and two $\mathcal O$-monoidal categories $\mathcal C^{\otimes} , \mathcal D^\otimes$, an $\mathcal O$-monoidal functor is an operad map $f\colon \mathcal C^\otimes \to \mathcal D^\otimes$ over $\mathcal O^\otimes$ which sends cocartesian morphisms to cocartesian morphisms. We can use the same principle to generalize this definition to $\infty$-operad pairs.

\begin{defn}
    Let $p\colon \opp \to \mathrm{N}\Fint$ be an $\infty$-operad pair and let $\mathcal C, \mathcal D$ be two $\mathcal P$-bimonoidal categories. A $\mathcal P$-\emph{bimonoidal functor} is a morphism of $\infty$-operad pairs $f \colon \mathcal C^{\oplus,\otimes} \to \mathcal D^{\oplus,\otimes}$ in the sense of \autoref{mapsofpairs} over $\mathcal P^{\oplus,\otimes}$ which sends cocartesian morphisms to cocartesian morphisms.
\end{defn}

\begin{defn}
   Let $\overline {\mathsf{Min}}_* \subset \Fint$ denote the full subcategory of $\Fint$ generated by morphisms in (1), (2) and (3) in \autoref{fin2}. That is, we forget the data of the distributive map $D$. For any $\infty$-operad pair $p\colon \opp\to \Fint$, we use $\overline{\mathcal P}^{\oplus, \otimes}$ to denote the pullback

   \[\begin{tikzcd}
	{\overline{\mathcal P}^{\oplus,\otimes}} & \opp \\
	{\mathrm{N}\overline{\mathsf{Min}}_*} & {\mathrm{N}\Fint}
	\arrow[from=1-1, to=1-2]
	\arrow["{\overline p}"', from=1-1, to=2-1]
	\arrow["p", from=1-2, to=2-2]
	\arrow[hook, from=2-1, to=2-2]
	\arrow["\ulcorner"{anchor=center, pos=0.125, rotate=180}, draw=none, from=2-2, to=1-1]
\end{tikzcd}\]
\end{defn}

\begin{prop}
    Let $f\colon \mathcal C^{\oplus,\otimes} \to \mathcal D^{\oplus, \otimes}$ be a $\mathcal P$-monoidal functor. Then the followings are equivalent:
    \begin{enumerate}[(1)]
        \item The underlying functor $\underline{f} \colon \mathcal C \to \mathcal D$ is an equivalence.
        \item The pullback functor $\overline{f}\colon \overline{\mathcal C}^{\oplus, \otimes} \to \overline{\mathcal D}^{\oplus,\otimes}$ is an equivalence.
        \item For every $X \in \opp$, the induced functor $\mathcal C_X \to \mathcal D_X$ is an equivalence.
    \end{enumerate}
\end{prop}
\begin{proof}
    By \autoref{bimonoidal}, $\overline{\mathcal C}^{\oplus,\otimes} \to \mathrm{N}\overline{\mathsf{Min}}_*$ is a cocartesian fibration and the resulting map $\overline {\mathcal C}^{\oplus, \otimes} \to \overline{\mathcal D}^{\oplus, \otimes}$ is a map of cocartesian fibrations. It sends cocartesian morphisms to cocartesian morphisms. (2) $\Leftrightarrow$ (3) then follows from \cite[Corollary 2.4.4.4]{lurie2009higher}. It is clear that (2) $\Rightarrow$ (1), so it remains to show (1) $\Rightarrow$ (3). If $X \in \mathcal P = \opp_{\langle 1 \rangle}$, it follows immediately. In general, suppose that $X \in \opp_{(\langle m\rangle , \gamma)}$, we have a commutative square induced by $f$:
\[\begin{tikzcd}
	{\mathcal C^{\oplus,\otimes}_{(\langle m\rangle,\gamma)}} & {\mathcal P^{\oplus,\otimes}_{(\langle m\rangle,\gamma)}} \\
	{\mathcal C^{\times m}} & {\mathcal P^{\times m}}
	\arrow["p", from=1-1, to=1-2]
	\arrow["{\prod \rho_i}"', from=1-1, to=2-1]
	\arrow["{\prod \rho_i}", from=1-2, to=2-2]
	\arrow["p"', from=2-1, to=2-2]
\end{tikzcd}\]  where the vertical arrows are equivalences by operad conditions. Passing to the fiber over $X \simeq X_1\odot\cdots \odot X_m$, we obtain an equivalence $\mathcal C_X \simeq \prod_i \mathcal C_{X_i}$. The diagram is functorial with respect to $\mathcal C$. Therefore, the operad map $f$ induces an equivalence $\mathcal C_X \simeq \prod_i \mathcal C_{X_i} \simeq \prod_i \mathcal D_{X_i} \simeq \mathcal D_X$, proving (1) $\Rightarrow$ (3). 
\end{proof}
\section{Day convolutions from an operad pair}\label{daysection}
Let $\mathcal C$ and $\mathcal D$ be two symmetric monoidal categories with additional presentability assumptions. Classically one can define a symmetric monoidal structure over $\Fun (\mathcal C , \mathcal D)$ called the Day convolution monoidal structure, which is given by a left Kan extension. This is generalized to $\infty$-operads by Lurie \cite{lurie2017higher} using a norm construction. However, it becomes more complicated than expected when forming Day convolutions for $\infty$-operad pairs. Basically, we have at least two different settings:

\begin{enumerate}[(1)]
  \item Let $\opp$ and $\mathcal Q^{\oplus,\otimes}$ be two operad pairs. We can consider $\Fun(\mathcal P,\mathcal Q)$ equipped with two different Day-convolution operad structures.
  \item Let $\opp$ be an operad pair and $\mathcal O^\otimes$ be an operad. We can also consider $\Fun(\mathcal P,\mathcal O)$ equipped with two different Day-convolution operad structures.
\end{enumerate}

Unlike the mono-operad case, where we can define an operad with the underlying category identified with the functor category between the underlying categories of two operads. One should not expect the same thing to happen for operad pairs. Classically, suppose that we have two bimonoidal categories $(\mathcal C, \oplus_\mathcal C,\otimes_\mathcal C,0,\mathbbl{1})$ and $(\mathcal D, \oplus_\mathcal D,\otimes_\mathcal D,0,\mathbbl{1})$, the two Day convolution monoidal structures on the functor category do not interact distributively. Instead, what one may expect is that the subcategory generated by $\oplus$-monoidal functor is closed under $\otimes$-Day convolution, which is classically correct when $\mathcal C$ is semi-additive and the underlying operad is cocartesian symmetric monoidal. 

\begin{prop}
    Let $\mathcal A$ be a semi-additive category equipped with a closed symmetric monoidal structure $\otimes_{\mathcal A}$. Let $(\mathcal B,\oplus_\mathcal B, \otimes_\mathcal B, 0, \mathbbl{1})$ be another symmetric rig category. Let $\star_{\otimes}$ denote the Day convolution symmetric monoidal structure over $\Fun(\mathcal A,\mathcal B)$ defined by $\otimes_\mathcal A$ and $\otimes_\mathcal B$. Using $\oplus_A$ to denote the cartesian monoidal structure over $\mathcal A$, suppose that $\oplus_{\mathcal B}$ is closed symmetric monoidal, then the Day convolution of two $\oplus$-monoidal functors is again $\oplus$-monoidal.
\end{prop}
\begin{proof} 
   Let $R = F\star_{\otimes} G (X) \oplus_{\mathcal B} F\star_\otimes G (Y)$. We have\small{\begin{align*}
        R \simeq & \big(\mathrm{colim}_{A\otimes_{\mathcal A} B \to X} F(A) \otimes_{\mathcal B} G(B) \big)\oplus_{\mathcal B} \big( \mathrm{colim}_{C\otimes_{\mathcal A} D \to Y} F(C) \otimes_{\mathcal D} G(B)\big) \\ \simeq & \mathrm{colim}_{\substack{A\otimes_{\mathcal A} B \to X \\C\otimes_{\mathcal A} D \to Y}} \Big( (F(A) \otimes_{\mathcal B} G(B)) \oplus_{\mathcal B} (F(C) \otimes_{\mathcal B} G(D))\Big) \\ \simeq & \mathrm{colim}_{\substack{(A\oplus_\mathcal A C) \otimes_{\mathcal A} B \to X \\ (A\oplus_\mathcal A C) \otimes_{\mathcal A} D \to Y}} \Big( (F(A\oplus_{\mathcal A} C) \otimes G(B)) \oplus_{\mathcal B} (F(A \oplus_{\mathcal A} C) \otimes G(D))\Big) \\ \simeq & \mathrm{colim}_{(A\oplus_\mathcal A C) \otimes_{\mathcal A} (B\oplus_{\mathcal A} D) \to X \oplus_{\mathcal A} Y} \Big( F(A\oplus_{\mathcal A} C) \otimes (G(B) \oplus_{\mathcal B} G(D)) \Big) \\ \simeq &\mathrm{colim}_{(A\oplus_\mathcal A C) \otimes_{\mathcal A} (B\oplus_{\mathcal A} D) \to X \oplus_{\mathcal A} Y} \Big( F(A\oplus_{\mathcal A} C) \otimes G(B \oplus_{\mathcal A} D ) \Big) \\ \simeq & \mathrm{colim}_{K \otimes_\mathcal {A} J  \to X \oplus_{\mathcal A} Y} \Big( F(K) \otimes_{\mathcal B} G(J ) \Big)
    \end{align*}} and the right hand side is identified with $F\star_{\otimes} G (X\oplus Y)$. Here the second equivalence is given by commuting $\oplus_{\mathcal B}$ with colimits, the third and forth equivalences are given by distributivity and the cofinality provided by semiadditivity.
\end{proof} Our goal in this subsection is to understand these relations in $\infty$-categorical settings. To begin with, we need some notations.

Suppose that we are given an operad pair $\opp$, one may expect that the underlying operad $\mathcal O_0^\otimes$, when viewing as a category, is equipped with an operad structure given by the moperad structure. We are going to make this idea precise.


\begin{defn}\label{pfintv}
    We define a category $\pFintv$ equipped with a functor $\pFintv \to \Fintv$ as follows: 
    \begin{enumerate}
        \item The objects of $\pFintv$ are elements of $\mathrm{Seq}(\mathbb N)$, which are tuples $(k,j_1 ,\cdots,j_k) $ for $k \geq 1$ and $j_r \geq 0$ for $1 \leq r \leq k$.
        \item For every $f\colon \langle k \rangle \to \langle m \rangle \in \Finp$, there is a map $$\tilde{f}^\times \colon (k,j_1,\cdots,j_k) \to (m,h_1,\cdots, h_m)$$, where $h_i=\prod_{f(s)=i} j_s$.
        \item For every collection $\{f_i \colon \langle j_i\rangle \to \langle l_i\rangle  \in \Finp\} $ for $1\leq i \leq k$, there exists a map $\boxtimes_{i} \tilde{f}_i^+ \colon (k,j_1, \cdots, j_k) \to (k,i_1 , \cdots, i_k)$. 
        \item The essentially surjective functor $\varepsilon \colon \pFintv \to \Fintv$ sends $Q=(j_1,\cdots , j_r)$ to $1^Q_\times$, sends $\tilde{f}^\times$ to $(f^\times)^{\boxplus n}$, and $\boxtimes_i \tilde{f}_i^+$ to $f_1 \wr \cdots \wr f_k$.
    \end{enumerate}
\end{defn}

We have a composition of functors $i \circ \varepsilon \colon \pFintv \to \Fint$. Moreover, the functors $p^{\oplus}, p^{\otimes}\colon\Fintv \to \Fint$ induce two functors $\tilde p^{\oplus},\tilde p^{\otimes} \colon \pFintv \to \Fint$.

Informally, $\pFintv$ is a strengthening of $\Fintv$ where we can distinguish two different $1^{Q_1}_\times$ and $1^{Q_2}_\times$ even when $|Q_1| = |Q_2|$, as long as $Q_1$ and $Q_2$ comes from different tuples.

Let $p\colon \opp \to \mathrm{N}\Finp$ be an $\infty$-operad pair. We use $\mathrm{p}\opp$ to denote the pullback
 \[\begin{tikzcd}
	{\mathrm{p}{\opp}} & \opp \\
	{\mathrm{N}\pFintv} & {\mathrm{N}\Fint}
	\arrow["j", from=1-1, to=1-2]
	\arrow["{\tilde{p}}"', from=1-1, to=2-1]
	\arrow["p", from=1-2, to=2-2]
	\arrow["i\circ \varepsilon", from=2-1, to=2-2]
	\arrow["\ulcorner"{anchor=center, pos=0.125, rotate=180}, draw=none, from=2-2, to=1-1]
\end{tikzcd}\] For each $Q \in \pFintv$, the objects of the fiber $(\mathrm{p}\opp)_{Q}$ can be identified with that of $(\opp)_{1^Q_\times}$.

\begin{defn}
    Let $p\colon\opp \to \mathrm{N}\Fint$ be an $\infty$-operad pair. We let $\opp_{\Delta}$ denote the full subcategory of $\mathrm{p}\opp$ generated by objects $C$ such that we have equivalences $(\rho^{m,q_m,l_1}_\times)_{!}(C) \simeq (\rho^{m,q_m,l_2}_\times)_{!}(C)$.
\end{defn}

For all objects $C \in (\opp_{\Delta})_{1^Q_\times} \subset \opp$, there is a cocartesian lift $\tilde d_Q \colon C \to C'$. Together with the Segal condition, we have the following:

\begin{lem}\label{identificationofdiagonal}
    The cocartesian lift $\tilde d_Q$ induces an equivalence
    \[
    (\opp_{\Delta})_{Q} \simeq (\opp)_{1^Q_+}
    \]
\end{lem}

\begin{lem}\label{poplusiscocart}
    The projection $\tilde p^{\oplus} \colon \mathrm{N}\pFintv \to \mathrm{N}\Finp $ is a cocartesian fibration of operads.
\end{lem}

\begin{proof}
    Let $Q= (k,j_1,\cdots,j_k)$. $Q \in \mathrm{N}\pFintv$ is sent to $\langle k \rangle \in \mathrm N \Finp$. For any $f\colon \langle k \rangle \to \langle m \rangle \in \Finp$, one can always lift it to $f^\times \colon (j_1,\cdots,j_k) \to (h_1 ,\cdots ,h_m)$, which is a cocartesian lift by construction. By \cite[Proposition 2.1.2.12]{lurie2017higher}, it suffices to show that we have an equivalence $(\mathrm N \pFintv)_{\langle k \rangle} \simeq (\mathrm N \pFintv)_{\langle 1 \rangle }^{\times k}$ induced by $\rho_i$. This follows from the observation that $(\mathrm N \pFintv)_{\langle 1 \rangle } \simeq \mathrm N\Finp$ and the description of morphisms in \autoref{pfintv}.
\end{proof}

\begin{thm}\label{opop}
    Let $p\colon \opp \to \mathrm{N}\Fint$ be a symmetric bimonoidal category. The composition $q=\tilde {p}^{\oplus} \circ \tilde p \colon \opp_{\Delta} \to \mathrm{N}\Finp$ is a cocartesian fibration of operads, with the underlying category identified with the underlying operad $\mathcal O_0^\otimes$ of $\opp$ viewing as a category.
\end{thm}

\begin{proof}
    Consider the following composition diagram:
    \[\begin{tikzcd}
	{\opp_{\Delta}} && {\mathrm{N}\pFintv} \\
	& {\mathrm{N}\Finp}
	\arrow["{\tilde{p}}", from=1-1, to=1-3]
	\arrow["q"', from=1-1, to=2-2]
	\arrow["{\tilde p^{\oplus}}", from=1-3, to=2-2]
\end{tikzcd}\] The functor $\tilde p$ is a cocartesian fibration. Indeed, for every map in $\pFintv$, its image in $\Fint$ always admits a cocartesian lift starting from any object in $\opp_{\Delta}$, and the target of this cocartesian lift also lies in $\opp_{\Delta}$. By \cite[Proposition 2.1.2.12]{lurie2017higher} and \autoref{poplusiscocart}, it suffices to show that for any $Q \simeq Q_1 \otimes\cdots \otimes Q_n$, we have an equivalence $(\opp_\Delta)_{Q} \simeq \prod_i (\opp_\Delta )_{Q_i}$. By \autoref{identificationofdiagonal}, we have an equivalence
\[
(\opp_\Delta)_Q \simeq \opp_{1^Q_+} \simeq \prod_i \opp_{1^{Q_i}_+} \simeq \prod(\opp_\Delta)_{Q_i}
\] which proves the claim.
\end{proof}

\begin{notn}\label{slicedsymmetric}
    Let $p\colon \opp \to \mathrm{N}\Fint$ be a symmetric bimonoidal category. We use $\mathcal P^\otimes_\oplus$ to denote the category $\opp_{\Delta}$ and use $p_{\oplus}\colon \mathcal P^\otimes_\oplus \to \mathrm{N}\Finp$ to denote the resulting operad in \autoref{opop}. We call it the \emph{underlying operad over operad} of $\opp$.
\end{notn}

We are going to study the Day convolution between operad pairs. Suppose we are given two symmetric monoidal categories, parameterized by operads $\mathcal O_1^\otimes$ and $\mathcal O_2^\otimes$, we have an identification:
\[
\Fun^{\mathrm{lax}}(\mathcal O_1,\mathcal O_2) \simeq \mathrm{Alg}_{\mathcal O_1}(\mathcal O_2)
\] The category $\mathrm{Alg}_{\mathcal O_1} (\mathcal O_2)$ is a full subcategory of $\Fun(\mathcal O_1^\otimes, \mathcal O_2^\otimes)$, which is generated by the functors that are over $\mathrm{N}\Finp$ and take inert morphisms to inert morphisms.

\begin{con}
    Let $p\colon \opp \to \mathrm{N}\Fint$ and $q\colon \mathcal {Q}^{\oplus,\otimes} \to \mathrm{N}\Fint$ be two symmetric bimonoidal categories. We use $\mathcal P^{\oplus}$ and $\mathcal Q^\oplus$ to denote their underlying operads. Let $p_{\oplus}\colon \mathcal P^\otimes_{\oplus} \to \mathrm{N} \Finp $ and $q_{\oplus} \colon \mathcal Q^{\otimes}_{\oplus} \to \mathrm{N}\Finp$ be the underlying operads over operad. By \cite[Example 2.2.6.9]{lurie2017higher}, we have a Day convolution symmetric monoidal category $\Fun (\mathcal P^{\oplus},\mathcal Q^{\oplus})^{\otimes} \to \mathrm{N}\Finp$ with underlying category identified with $\Fun (\mathcal P^{\oplus},\mathcal Q^{\oplus})$. We let $\mathrm{Alg}(\mathcal P,\mathcal Q)^{\otimes} \subset \Fun (\mathcal P^{\oplus},\mathcal Q^{\oplus})^{\otimes}$ be the full subcategory generated by the functors that live in $\mathrm{Alg}_{\mathcal P} (\mathcal Q)$ under the Segal condition.
\end{con}

\begin{thm}\label{twotwodayconvolution}
    Let $p\colon \opp \to \mathrm{N}\Finp$ and $q \colon \mathcal {Q}^{\oplus,\otimes} \to \mathrm{N}\Finp$ be two symmetric bimonoidal categories such that the underlying operad of $\mathcal Q^{\oplus, \otimes}$ is a closed symmetric monoidal structure. Suppose that the underlying category of $\opp$ is semi-additive, and the underlying operad of $\opp$ is the cartesian symmetric monoidal operad. Suppose that there exists an uncountable regular cardinal $\kappa$ such that:
    
    \begin{enumerate}
        \item $\mathcal P$ is essentially $\kappa$-small.
        \item $\mathcal Q$ admits $\kappa$-small colimits.
        \item For each $m \geq 1$ and $\phi\colon \langle m \rangle \to \langle 1 \rangle $, $\bigotimes _\phi \colon \prod_m \mathcal P^{\oplus} \to \mathcal P^{\oplus}$ preserves $\kappa$-small colimits. 
    \end{enumerate} Then $\mathrm{Alg}(\mathcal P,\mathcal Q)^\otimes \to \mathrm{N}\Finp$ is a symmetric monoidal category with the underlying category $\mathrm{Alg}_{\mathcal P} \mathcal Q \simeq \Fun^{\mathrm{lax}}(\mathcal P,\mathcal Q)$.
\end{thm}

\begin{proof}
    By construction, the inert morphisms $\rho_i$ induces equivalences \[\mathrm{Alg}(\mathcal P, \mathcal Q)^{\otimes}_{\langle n \rangle} \simeq \mathrm{Alg}_{\mathcal P}(\mathcal Q)^{\times n}\] By \cite[Proposition 2.1.2.12]{lurie2017higher}, it suffices to show that 
    \[r\colon \mathrm{Alg}(\mathcal P, \mathcal Q)^{\otimes} \to \mathrm{N}\Finp\] is a cocartesian fibration. That means, all cocartesian edges starting from objects in $\mathrm{Alg}(\mathcal P, \mathcal Q)^{\otimes}$ must also have targets in $\mathrm{Alg}(\mathcal P, \mathcal Q)^{\otimes}$. By Segal condition, it suffices to check this for cocartesian lifts of the active map $\phi\colon \langle m\rangle \to \langle 1 \rangle$. In order to show this, let us have a look at the cocartesian edges in the Day convolution operad. 

    Let $\{F_i \colon \mathcal {P}^{\oplus} \to \mathcal Q^{\oplus} \}$ be a collection of functors in $\mathrm{Fun}(\mathcal P^{\oplus} , \mathcal Q^{\oplus})$ for $1 \leq i \leq m$. Let $G$ be another functor in $\mathrm{Fun}(\mathcal P^{\oplus} , \mathcal Q^{\oplus})$. By \cite[Corollary 2.2.6.12]{lurie2017higher}, for any $\phi \colon \langle m \rangle \to \langle 1 \rangle \in \Finp$, we have an equivalence:
    \begin{align*}
            \Map^{\phi} (F_1\otimes\cdots\otimes F_m,G) & \simeq \mathrm{Mul}_{\Fun(\mathcal P^{\oplus},\mathcal Q^{\oplus})^{\otimes}}^{\phi} (\{F_i\},G) \\ & \simeq \Map_{\Fun(\prod_m\mathcal P^{\oplus}, \mathcal Q^{\oplus})} (\bigotimes_\phi \circ \prod_iF_i, G \circ \bigotimes_\phi)
    \end{align*} Under the equivalence, the cocartesian lift of $\phi$ corresponds to the left Kan extension:
    \[\begin{tikzcd}
	{\prod_m \mathcal P^{\oplus}} & {\prod_m \mathcal Q^{\oplus}} & {\mathcal Q^{\oplus}} \\
	{\mathcal P^{\oplus}}
	\arrow[""{name=0, anchor=center, inner sep=0}, "{\prod_iF_i}", from=1-1, to=1-2]
	\arrow["{\bigotimes_\phi}"', from=1-1, to=2-1]
	\arrow["{\bigotimes_{\phi}}", from=1-2, to=1-3]
	\arrow[""{name=1, anchor=center, inner sep=0}, "G"', dashed, from=2-1, to=1-3]
	\arrow["\alpha"', between={0.2}{0.8}, Rightarrow, from=0, to=1]
\end{tikzcd}\] 

In order to show that $\mathrm{Alg}(\mathcal P, \mathcal Q)$ is closed under taking cocartesian lifts, it suffices to show that if $F_i$ takes inert morphisms to inert morphisms for every $1 \leq i \leq m$. and let $G$ be the left Kan extension as shown above, then $G$ must be a functor over $\mathrm{N}\Finp$ and send inert morphisms to inert morphisms. For the ``over $\mathrm{N}\Finp$ '' part, note that for any $X \in \mathcal P^{\oplus}_{\langle n \rangle}$, any map $\bigotimes_{\phi} \prod_{1\leq i \leq m} A_i \to X $ lifting $f\colon \langle k \rangle \to \langle n \rangle$ must factors through $f_! \big( \bigotimes_{\phi} \prod_{1\leq i \leq m} A_i\big) \in \mathcal P^{\oplus}_{\langle n \rangle}$. By the cofinality provided by the semi-additivity, we have
\begin{align*}
    G(X) & \simeq \mathrm{colim}_{\bigotimes_\phi \prod C_i \to X} \hat{F} (\prod C_i) \\
    & \simeq \mathrm{colim}_{\substack{\bigotimes_\phi \prod C_i \to X \\ \bigotimes_{\phi}\prod C_i \in\mathcal P^{\oplus}_{\langle n \rangle}}} \hat F (\prod C_i)
\end{align*} This colimit is taken within each fiber, so $G(-)$ is a functor over $\mathrm{N}\Finp$ since each $F_i$ is a functor over $\mathrm{N}\Finp$.

It can be shown that $\otimes_\phi$ takes inert morphisms to inert morphisms. Indeed, for any inert morphism $f \colon D \to D'$, the image $\bigotimes_\phi (f)$ is defined through the cocartesian lifting of the diagram:

\[\begin{tikzcd}
	C & {\bigotimes_\phi (C)} \\
	{C'} & {\bigotimes_\phi (C')}
	\arrow["{\tilde \phi_1}", from=1-1, to=1-2]
	\arrow["f"', from=1-1, to=2-1]
	\arrow["{\bigotimes_\phi(f)}", dashed, from=1-2, to=2-2]
	\arrow["{\tilde \phi_2}"', from=2-1, to=2-2]
\end{tikzcd}\] where $\tilde \phi_1$ and $\tilde \phi_2$ are cocartesian lift of $\phi$. By (the dual of) \cite[Proposition 2.4.1.3]{lurie2009higher}, this implies that $\bigotimes_\phi (f)$ is also inert. In order to show that $G$ takes inert morphisms to inert morphisms, let $\hat{F}$ denote the composition $\bigotimes_{\phi}\circ \prod_i F_i$. Let $\alpha\colon \langle n \rangle \to \langle l \rangle$ be an inert map in $\mathrm{N}\Finp$. By straightening, it suffices to show the following diagram commutes:

\[\begin{tikzcd}
	{\mathcal P^{\oplus}_{\langle n \rangle}} & {\mathcal P^{\oplus}_{\langle l \rangle}} \\
	{\mathcal Q^{\oplus}_{\langle n \rangle}} & {\mathcal Q^{\oplus}_{\langle l \rangle}}
	\arrow["{\alpha_!}", from=1-1, to=1-2]
	\arrow["{(\mathrm{Lan}_{\bigotimes_\phi}\hat F)_{\langle n \rangle}}"', from=1-1, to=2-1]
	\arrow["{(\mathrm{Lan}_{\bigotimes_\phi}\hat F)_{\langle l \rangle}}", from=1-2, to=2-2]
	\arrow["{\alpha_!}"', from=2-1, to=2-2]
\end{tikzcd}\] Unraveling everything, for every $C \in \mathcal P^{\oplus}_{\langle n \rangle}$, we have:

\begin{align*}
    \alpha_! (\mathrm{Lan}_{\bigotimes_\phi}\hat F)_{\langle n \rangle} (C) & \simeq\alpha_!\mathrm{colim}_{\substack{\bigotimes_\phi \prod C_i \to C \\ \bigotimes_{\phi}\prod C_i \in\mathcal P^{\oplus}_{\langle n \rangle}}} \hat F (\prod C_i) \\
    &\simeq\mathrm{colim}_{\substack{\bigotimes_\phi \prod C_i \to C \\ \bigotimes_{\phi}\prod C_i \in\mathcal P^{\oplus}_{\langle n \rangle}}} \alpha_!\hat F (\prod C_i)\\
    & \simeq \mathrm{colim}_{\substack{\bigotimes_{\phi} \prod C_i \to \alpha_!C \\ \bigotimes_{\phi}\prod C_i \in\mathcal P^{\oplus}_{\langle l \rangle}}} \hat F(\prod C_i) \\
    & \simeq (\mathrm{Lan}_{\bigotimes_\phi}\hat F)_{\langle l \rangle} \alpha_! (C)
\end{align*} where the first two equivalences follow from commuting $\alpha_!$ with colimits, and the last two equivalences are by the cofinality of the diagram. \qedhere

\end{proof} The proof of the theorem can also be applied to the subcategory $\Fun^{\oplus}(\mathcal P. \mathcal Q) \subset \mathrm{Alg}(\mathcal P , \mathcal Q)$ consisting of $\oplus$-monoidal functors. Therefore, we obtain the following:

\begin{cor}
    Let $\mathcal P^{\oplus,\otimes}$ and $\mathcal Q^{\oplus,\otimes}$ satisfy the condition as in \autoref{twotwodayconvolution}. Suppose additionally that the operator $\oplus$ in $\mathcal Q$ commutes with colimits. Then $\mathrm{Alg}(\mathcal P,\mathcal Q)^{\otimes} \to \mathrm{N}\Finp$ restricts to a symmetric monoidal subcategory \[\Fun^{\oplus}(\mathcal P. \mathcal Q)^\otimes \to \mathrm{N}\Finp\]
\end{cor}
\begin{proof}
    This follows from a same proof as in \autoref{twotwodayconvolution} word for word.
\end{proof}
Our next goal is to study the Day convolution between a symmetric bimonoidal category $\opp$ and a symmetric monoidal category $\mathcal C^{\otimes}$. This becomes of vital importance in the next section, where we begin to study functors from a span category to a cartesian monoidal category and no additional operad pair structure is considered on the target. 

Again, since the source category is equipped with two different monoidal structures, there are two subcategories of the functor category, namely $\Fun^{\oplus\text{-}\mathrm{lax}} (\mathcal P,\mathcal C)$ and $\Fun^{\otimes\text{-}\mathrm{lax}}(\mathcal P,\mathcal C)$. In the case when the target category $\mathcal C$ is cartesian, we can also talk about how a functor is moreover distributive. Informally, a distributive functor should be some $F\colon \mathcal P \to \mathcal C$ such that:

\begin{enumerate}
    \item $F$ is $\oplus$-lax monoidal, that is, we have a natural comparison morphism:
    \[
    F(A)\times F(B) \to F(A\oplus B)
    \]
    \item $F$ is $\otimes$-lax monoidal, that is, we have a natural comparison morphism:
    \[
    F(A) \times F(B) \to F(A\otimes B)
    \]
    \item Let $Q=(k,j_1 ,\cdots,j_k)\in\mathrm{Seq}(\mathbb N)$ and suppose we are given a distributive equivalence 
        \[
        \prod_{|Q|} \prod_{1\leq r \leq k} A_{r,q_r} \simeq \prod_{1 \leq r \leq k} \prod_{1 \leq q \leq j_r} A_{r,q}
        \] The following distributive square commutes, where $\Delta$ is the diagonal morphism, $D$ is the distributive equivalence in $\mathcal P$ and the vertical arrows are comparisons.
\[\begin{tikzcd}
	{\prod_{|Q|} \prod_{1\leq r \leq k} F ( A_{r,q_r})} & {F(\bigoplus_{|Q|}\bigotimes_{1\leq r\leq k} A_{r,q_r})} \\
	{\prod_{1 \leq r \leq k} \prod_{1 \leq q \leq j_r} F(A_{r,q})} & {F(\bigotimes_k\bigoplus_{1 \leq q \leq j_r} A_{r,q})}
	\arrow[from=1-1, to=1-2]
	\arrow["D", from=1-2, to=2-2]
	\arrow["\Delta", from=2-1, to=1-1]
	\arrow[from=2-1, to=2-2]
\end{tikzcd}\]
\end{enumerate} In general, those functors that are $\otimes$-lax monoidal and $\oplus$-lax monoidal simultaneously are not necessarily become distributive automatically. However, they can be identified with certain commutative algebra objects.

\begin{con}
    Let $p\colon \opp \to \mathrm{N}\Fint$ be a symmetric bimonoidal category and let $s\colon \mathcal C^\oplus \to \mathrm{N}\Finp $ be a symmetric monoidal category. If $\mathcal C$ admits finite product, then so is $\mathcal C^\oplus$. We let $s_{\oplus}\colon(\mathcal C^{\oplus})^{\times} \to \mathrm{N}\Finp$ denote the cartesian fibration defined in \cite[Section 2.4.1]{lurie2017higher} and let $p_{\oplus}\colon \mathcal P_{\oplus}^{\otimes} \to \mathrm{N} \Finp$ denote the symmetric monoidal category as in \autoref{slicedsymmetric}. We then have a Day convolution symmetric monoidal category
    \[
    r\colon \Fun(\mathcal P^{\oplus}, \mathcal C^{\oplus})^{\otimes} \to \mathrm{N}\Finp 
    \] 
\end{con}

\begin{cor}\label{twooneoperad}
     Let $p\colon \opp \to \mathrm{N}\Fint$ be a symmetric bimonoidal category and let $s\colon \mathcal C^\oplus \to \mathrm{N}\Finp $ be a closed symmetric monoidal category. Suppose that $\mathcal P$ satisfies the conditions in \autoref{twotwodayconvolution} and $\mathcal C$ is presentable and admits finite products. Let $\mathrm{Alg}(\mathcal P,\mathcal C) \subset \Fun(\mathcal P^{\oplus},\mathcal C^{\oplus})$ (respectively $\Fun^{\oplus}(\mathcal P,\mathcal C)$) denote the full subcategory generated by operad maps (respectively, symmetric monoidal functors). Then $r$ defined above restrict to a symmetric monoidal category with underlying category identified with $\mathrm{Alg}_\mathcal P \mathcal C \simeq \Fun^{\mathrm{lax}}(\mathcal P,\mathcal C)$ (respectively, $\Fun^{\oplus}(\mathcal P,\mathcal C)$).
\end{cor}

\begin{proof}
    We can apply the proof of \autoref{twotwodayconvolution} word for word and replace the operad $\mathcal Q^{ \otimes}_{\oplus}$ by $(\mathcal C^{\oplus})^\times$.
\end{proof}

We now give our first definition of distributive functors.

\begin{defn}\label{operadpairdistributive}
    Let $p\colon \opp \to \mathrm{N}\Fint$ and $s\colon \mathcal C^\oplus \to \mathrm{N}\Finp $ be defined as above. 
    \begin{enumerate}
        \item We say a functor $F \colon \mathcal P \to \mathcal C$ is \emph{lax-distributive} if it is a commutative algebra object in $\mathrm{Alg}_{\mathcal P}\mathcal C \simeq \Fun^{\mathrm{lax}}(\mathcal P,\mathcal C)$. That is, we define
        \[
        \Fun^{\mathrm{ldist}}(\mathcal P,\mathcal C) := \mathrm{CAlg}(\Fun^{\mathrm{lax}}(\mathcal P,\mathcal C))
        \]
        \item We say $F$ is \emph{distributive} if it is a commutative algebra object in $\Fun^{\oplus}(\mathcal P,\mathcal C)$. That is, we define
        \[
        \Fun^{\mathrm{dist}}(\mathcal P,\mathcal C):= \mathrm{CAlg}(\Fun^{\oplus}(\mathcal P,\mathcal C))
        \]
    \end{enumerate}
\end{defn}

\begin{rem}\label{explanationfordist}
    When $\mathcal C^{\oplus} \simeq \mathcal C^{\times} $ is the cartesian symmetric monoidal structure over $\mathcal C$, the definition coincides with the classical ones. To see this, let us unpack the definition.

    Let $+$ and $\cdot$ denote two symmetric monoidal operators in $\mathcal P$ corresponding to $\oplus$ and $\otimes$. Let $\tilde \times$ denote the cartesian operadic operator in $(\mathcal C^{\oplus})^{\times}$. Let $F$ be a (lax) distributive functor from $\mathcal P$ to $\mathcal C$. We have the following commutative diagram:
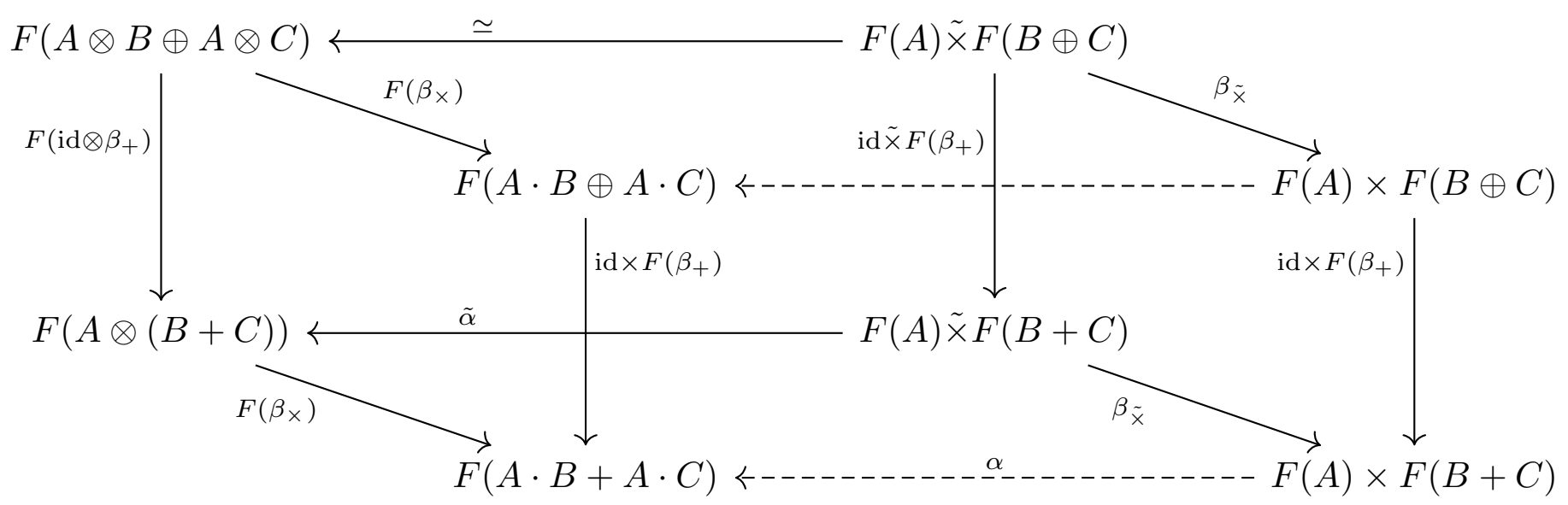
\[\begin{tikzcd}[nodes={scale=0.8}]
	{F(A\otimes B\oplus A\otimes C)} && {F(A)\tilde{\times}F(B\oplus C)} & \\
	& {F(A\cdot B \oplus A \cdot C)} && {F(A)\times F(B\oplus C)} \\
	{F(A \otimes(B+C))} && {F(A)\tilde \times F(B+C)} \\
	& {F(A\cdot B + A\cdot C)} && {F(A)\times F(B+C)}
	\arrow["{F(\beta_\times)}", from=1-1, to=2-2]
	\arrow["{F(\mathrm{id}\otimes\beta_+)}"'{pos=0.3}, from=1-1, to=3-1]
	\arrow["\simeq"'{pos=0.7}, from=1-3, to=1-1]
	\arrow["{\beta_{\tilde \times}}", from=1-3, to=2-4]
	\arrow["{\mathrm{id}\tilde{\times}F(\beta_+)}"'{pos=0.3}, from=1-3, to=3-3]
	\arrow["{\mathrm {id} \times F(\beta_+)}"{pos=0.2}, from=2-2, to=4-2]
	\arrow[dashed, from=2-4, to=2-2]
	\arrow["{\mathrm{id} \times F (\beta_+)}"'{pos=0.2}, from=2-4, to=4-4]
	\arrow["{F(\beta_{\times})}"'{pos=0.3}, from=3-1, to=4-2]
	\arrow["{\tilde\alpha}"'{pos=0.7}, from=3-3, to=3-1]
	\arrow["{\beta_{\tilde \times}}"'{pos=0.3}, from=3-3, to=4-4]
	\arrow["\alpha"', dashed, from=4-4, to=4-2]
\end{tikzcd}\]The back surface sits in $(\mathcal C^{\oplus})^\times_{\langle 2 \rangle}$ and the front surface sits in $\mathcal C^{\oplus} \simeq (\mathcal C^{\oplus})^\times_{\langle 1 \rangle}$. For $X = A \otimes B \oplus A \otimes C \in (\mathcal P_\oplus^\otimes)_{\langle 2 \rangle}$, we have $(\rho_1)_! X \simeq A$ and $(\rho_2)_! X \simeq B\oplus C$. Therefore the commutative algebra structure induces the top equivalence in the backward surface. The dotted arrow $\alpha$ is induced by the cocartesian edge $\beta_{\tilde \times}$. Especially, we have another comparison map

 \[ F(A) \times F(B) \to F(A\cdot B)
 \] Together with the $\oplus$-monoidality gives 

 \[
 \alpha_1 \colon F(A )\times F(B)\times F(C) \to F(A\cdot B + A \cdot C)
 \]
 \[
 \alpha_2 \colon F(A)\times F(B)\times F(A)\times F(C) \to F(A \cdot B +A\cdot C)
 \] Finally, starting from $F(A \otimes B \oplus A \otimes C)$, we have two morphisms:
\[\begin{tikzcd}
	{F(A\otimes B \oplus A \otimes C)} & {F(A\otimes B)} & {F(A)}
	\arrow["{F(\rho_1^+)}", from=1-1, to=1-2]
	\arrow["{F(\rho_1^\otimes)}", from=1-2, to=1-3]
\end{tikzcd}\]
\[\begin{tikzcd}
	{F(A\otimes B \oplus A \otimes C)} & {F(A\otimes C)} & {F(A)}
	\arrow["{F(\rho_2^+)}", from=1-1, to=1-2]
	\arrow["{F(\rho_1^\otimes)}", from=1-2, to=1-3]
\end{tikzcd}\] These two coincide by relations in $\opp$. After precomposing with the equivalence $F(A) \tilde \times F(B\oplus C) \simeq F(A\otimes B \oplus A \otimes B)$, both of them becomes projection to $F(A)$. This implies on category level we have a factor through:
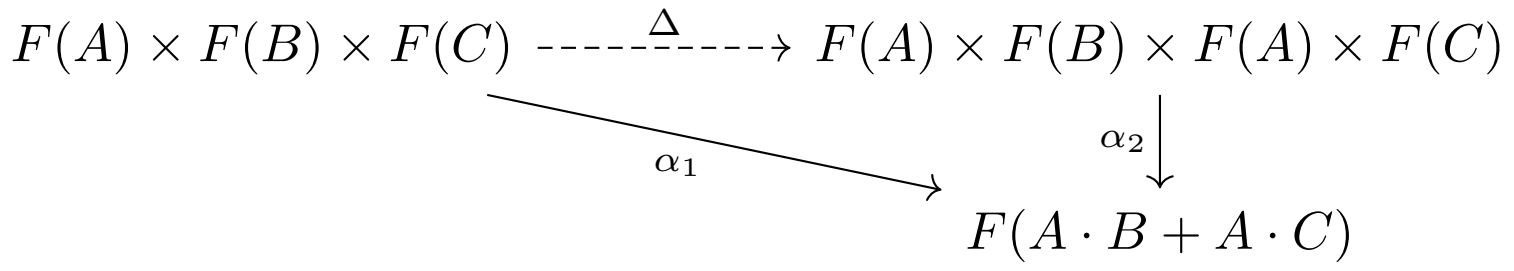
\[\begin{tikzcd}
	{F(A)\times F(B)\times F(C)} && {F(A)\times F(B) \times F(A) \times F(C)} \\
	&& {F(A\cdot B +A\cdot C)}
	\arrow["\Delta", dashed, from=1-1, to=1-3]
	\arrow["{\alpha_1}"', from=1-1, to=2-3]
	\arrow["{\alpha_2}"', from=1-3, to=2-3]
\end{tikzcd}\]
\end{rem}

\section{The universal property of spans and bispans}\label{spansection}
The most important ingredients to define spectral Mackey functors and spectral Tambara functors are categories of spans and bispans. In this section we review the construction of the span and bispan category and how they are related together. Moreover, we will construct a symmetric bimonoidal structure over the category of spans in any category. This will plays an important role in the next section.

\begin{defn}[\cite{haugseng2018iterated}, \cite{heyer20246}]
    A \emph{span pair} $(\mathcal C ,\mathcal C_F)$ consists of an $\infty$-category $\mathcal C$ together with a wide subcategory $\mathcal C_F \subset \mathcal C$ such that:
    \begin{enumerate}[(a)]
        \item $\mathcal C_F$ is closed under pullbacks along morphisms in $\mathcal C$.
        \item $\mathcal C_F$ admits and the inclusion $\mathcal C_F \subset \mathcal C$ preserves pullback diagrams.
    \end{enumerate}
\end{defn}

We can formally organize all span pairs into a category, which is denoted by $\mathsf{Pair}$. One can show that the category of span pairs is complete, see for example \cite[Section 2.1]{heyer20246}. 
\begin{exmp}
    Let $\mathbbl {\Sigma}^n$ denote the poset of pairs $(i,j)$ such that $0 \leq i \leq j \leq n$. The partial ordering is given by $(i,j)\leq (k,l)$ if $i\leq k$ and $j\leq l$. Let $\mathbbl {\Sigma}^n_2$ be the sub-poset with the same elements, while the partial ordering is given by $(i,j)\leq (k,l)$ if $i\leq k$ and $j= l$. Then $(\mathbbl {\Sigma}^n, \mathbbl{\Sigma}^n_2)$ becomes a span pair. Moreover, the construction gives a cosimplicial object
    \[
    (\mathbbl {\Sigma}^\bullet, \mathbbl{\Sigma}^{\bullet}_2) \colon \Delta \to \mathsf{Pair}
    \]
\end{exmp}

\begin{defn}[\cite{heyer20246}]
    We define a functor
    \[
    \Span \colon \mathsf{Pair} \to \Fun(\Delta^{\mathrm{op}}, \Ani)
    \] through the composition
\[\begin{tikzcd}
	{\Delta^{\mathrm{op}}\times \mathsf{Pair}} && {\mathsf{Pair}^{\mathrm{op}}\times\mathsf{Pair}} & \Ani
	\arrow["{(\mathbbl{\Sigma}^\bullet,\mathbbl{\Sigma}^\bullet_2)\times \mathrm{id}}", from=1-1, to=1-3]
	\arrow["{\mathrm{Hom}}", from=1-3, to=1-4]
\end{tikzcd}\] By \cite{heyer20246}, this functor factors through the full subcategory of Segal Animas, thus defines a functor
\[
\Span\colon \mathsf{Pair} \to \mathrm{Cat}_{\infty}
\]
    
\end{defn}
Let $(\mathcal C,\mathcal C_F)$ be a span pair. The $(\infty,1)$-category $\Span(\mathcal C,\mathcal C_F)$ can be described as follows: the objects of $\Span(\mathcal C,\mathcal C_F)$ are the same as $\mathcal C$. For two objects $X$ and $Y$, a morphism from $X$ to $Y$ is given by spans
\[\begin{tikzcd}
	& U & \\
	X && Y
	\arrow["f"', from=1-2, to=2-1]
	\arrow["g", from=1-2, to=2-3]
\end{tikzcd}\] where $g$ is in $\mathcal C_F$. The composition is given by pullbacks
\[\begin{tikzcd}
	&& W && \\
	& U && V \\
	X && Y && Z
	\arrow["{i'}"', from=1-3, to=2-2]
	\arrow["{j'}", from=1-3, to=2-4]
	\arrow["f"', from=2-2, to=3-1]
	\arrow["g", from=2-2, to=3-3]
	\arrow["i"', from=2-4, to=3-3]
	\arrow["j", from=2-4, to=3-5]
	\arrow["\ulcorner"{anchor=center, pos=0.125, rotate=135}, draw=none, from=3-3, to=1-3]
\end{tikzcd}\]

The category of span enjoys a beautiful universal property. It can be viewed as the``universal source" of those ``$F$-left adjointable" functors.

\begin{defn}
    Let $(\mathcal C,\mathcal C_F)$ be a span pair and let $\mathcal D$ be another $\infty$-category. A functor $G \colon \mathcal C^{\mathrm{op}} \to \mathcal D$ is \emph{left $F$- adjointable} if:
    \begin{enumerate}
        \item For every $g\colon X \to Y$ in $\mathcal C_F$, there is a functorial morphism $g_{!}\colon G(X) \to G(Y)$ serving as the ``left adjoint" of $F(g)$.
        \item For any $f\colon A \to B$ in $\mathcal C$, we denote $f^*:=G(f)$. Then the following Beck-Chevalley condition is satisfied. For any cartesian square:
        \[\begin{tikzcd}
	{X'} & X \\
	{Y'} & Y
	\arrow["{g'}", from=1-1, to=1-2]
	\arrow["{f'}"', from=1-1, to=2-1]
	\arrow["f", from=1-2, to=2-2]
	\arrow["g"', from=2-1, to=2-2]
	\arrow["\ulcorner"{anchor=center, pos=0.125, rotate=180}, draw=none, from=2-2, to=1-1]
\end{tikzcd}\] we always have $g'_!f'^* \simeq f^*g_!$.
    \end{enumerate}
\end{defn}

\begin{warn}
    Throughout the paper we are working with $(\infty,1)$-categories. Therefore, our definition is in fact the "degenerated" version of the $(\infty,2)$-definition of \cite{gaitsgory2019study} and \cite{elmanto2023distributivity}. In the $(\infty, 2)$-categorical settings, not only should we require $g_!$ to be a left adjoint of $G(g)$, but also the equivalence $g'_!f'^* \simeq f^*g_!$ should be given by Beck-Chevalley.
\end{warn}

The left $F$-adjointable functors form a category $\Fun_{F\text{-}\mathrm{ladj}} (\mathcal C^{\mathrm{op}}, \mathcal D)$. There is a natural forgetful functor $u\colon \Fun_{F\text{-}\mathrm{ladj}} (\mathcal C^{\mathrm{op}}, \mathcal D) \to \Fun(\mathcal C^{\mathrm{op}},\mathcal D)$ by forgetting the left adjointable data. The following result is a degenerated version of \cite[Theorem 3.2.2]{gaitsgory2019study} \cite[Theorem 4.2.6]{macpherson2022bivariant} in $(\infty,1)$-categories.

\begin{prop}\label{spand}
    Let $(\mathcal C,\mathcal C_F)$ be a span pair. There exists an equivalence \[\Fun(\Span (\mathcal C,\mathcal C_F),\mathcal D) \simeq \Fun_{F\text{-}\mathrm{ladj}}(\mathcal C^{\mathrm{op}},\mathcal D)\] fitting into the following commutative diagram
    \[\begin{tikzcd}
	{\Fun(\Span (\mathcal C,\mathcal C_F),\mathcal D)} && {\Fun_{F\text{-}\mathrm{ladj}}(\mathcal C^{\mathrm{op}},\mathcal D)} \\
	& {\Fun (\mathcal C^{\mathrm{op}} , \mathcal D)}
	\arrow["\simeq", from=1-1, to=1-3]
	\arrow["i^*"', from=1-1, to=2-2]
	\arrow["u", from=1-3, to=2-2]
\end{tikzcd}\] where $i^*$ is induces by the inclusion $i\colon \mathcal C^{\mathrm{op}} \to \Span (\mathcal C,\mathcal C_F)$ and $u$ is the forgetful functor.
\end{prop}

\begin{proof}
    Let $\Phi\colon \Span(\mathcal C,\mathcal C_F) \to \mathcal D$ be a functor. The left adjoint of $\Phi|_{\mathcal C^{\mathrm{op}}}$ is exactly given by $\Phi|_{\mathcal C_F}$. The Beck-Chevalley condition is exactly the composition formula in $\Span (\mathcal C,\mathcal C_F)$. 
\end{proof}
We now consider the category of bispans. 

\begin{defn}[\cite{elmanto2023distributivity}]\label{distdiagram}
    Let $X\xrightarrow{l}Y \xrightarrow{f} Z$ be morphisms in an $\infty$-category $\mathcal C$. A \emph{distributive diagram} for $l$ and $f$ is a commutative diagram

    \[\begin{tikzcd}
	& {W\times_ZY} & W \\
	X & Y & Z
	\arrow["{\tilde f}", from=1-2, to=1-3]
	\arrow["\epsilon"', from=1-2, to=2-1]
	\arrow["{\tilde g}", from=1-2, to=2-2]
	\arrow["g"', from=1-3, to=2-3]
	\arrow["l", from=2-1, to=2-2]
	\arrow["f", from=2-2, to=2-3]
	\arrow["\ulcorner"{anchor=center, pos=0.125, rotate=180}, draw=none, from=2-3, to=1-2]
\end{tikzcd}\] where the square is cartesian, with the property that for any morphism $\phi\colon U \to Z$, the composition
\[
\Map_{/Z} (\phi, g) \to \Map_{/Y}(f^*\phi,\tilde g) \xrightarrow{\epsilon_*} \Map_{/Y}(f^*\phi, l)
\] is an equivalence. Note that the distributive diagram is unique whenever it exists.
\end{defn}

\begin{exmp}\label{noneqexmp}
    Let $\mathcal C=\mathsf{Fin}$ be the category of finite sets. Then for any morphisms $I \xrightarrow{u}J\xrightarrow{v} K$, the distributive diagram exists and is of the form
    \[\begin{tikzcd}
	& {\coprod_{j\in J}\prod_{j'\in J_{v(j)}}I_{j'}} & {\coprod_{k\in K}\prod_{j\in J_k}I_j} \\
	I & J & K
	\arrow["{\tilde v}", from=1-2, to=1-3]
	\arrow["\epsilon"', from=1-2, to=2-1]
	\arrow["{\tilde h}", from=1-2, to=2-2]
	\arrow["h", from=1-3, to=2-3]
	\arrow["u", from=2-1, to=2-2]
	\arrow["v", from=2-2, to=2-3]
	\arrow["\ulcorner"{anchor=center, pos=0.125, rotate=180}, draw=none, from=2-3, to=1-2]
\end{tikzcd}\] Here $h$ and $\tilde h$ sends each coproduct summand to $\mathbf{1} \in \mathsf{Fin}$. Let $L={\prod_{j'\in J_{v(j)}} I_{j'}}$. The morphism $\epsilon$ is defined via the distributive law and codiagonal:
\[
\epsilon\colon \coprod_{j\in J}\prod_{j'\in J_{v(j)}}I_{j'} \simeq \coprod_{j\in J}\coprod_{L} \mathbf{1} \xrightarrow{\Delta} \coprod_{i\in J} \coprod_{I_j} \mathbf{1} \simeq I
\] 
\end{exmp}

\begin{exmp}
    Let $G$ be a finite set and let $\mathsf{Fin}_G$ be the category of finite $G$-sets. Then each $X \in \mathsf{Fin}_G$ can be written as $X \cong X_1\sqcup \cdots \sqcup X_m$ where each $X_i$ is a $G$-orbit and any morphisms between finite $G$-sets can be written as a coproduct of morphisms over these summands. For any morphism $A \xrightarrow{u} B \xrightarrow{v} C$, written them as the orbit decomposition $\coprod_{i\in I} A_i \xrightarrow{u} \coprod_{j\in J} B_j \xrightarrow{v} \coprod_{k\in K} C_k$. Then the distributive diagram exists and is of the form 
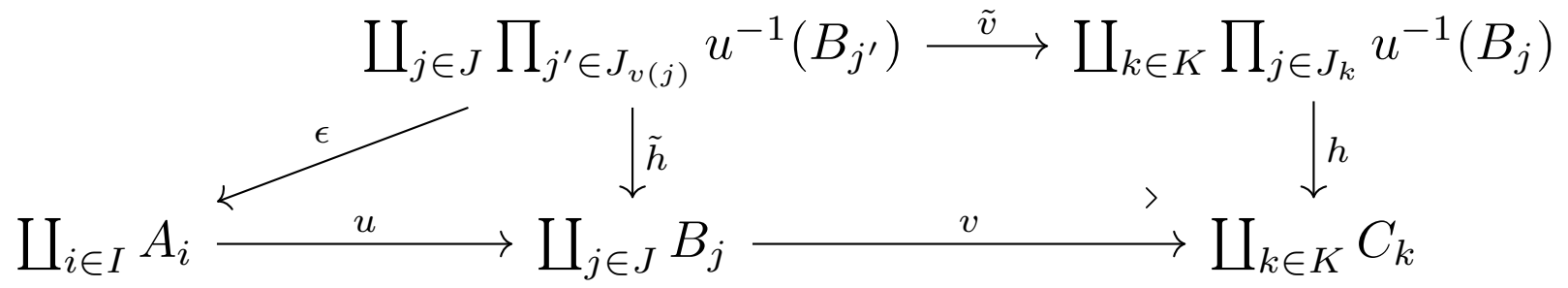
\[\begin{tikzcd}
	& {\coprod_{j\in J}\prod_{j'\in J_{v(j)}}u^{-1}(B_{j'})} & {\coprod_{k\in K}\prod_{j\in J_k}u^{-1}(B_j)} \\
	{\coprod_{i\in I}A_i} & {\coprod_{j\in J}B_j} & {\coprod_{k\in K}C_k}
	\arrow["{\tilde v}", from=1-2, to=1-3]
	\arrow["\epsilon"', from=1-2, to=2-1]
	\arrow["{\tilde h}", from=1-2, to=2-2]
	\arrow["h", from=1-3, to=2-3]
	\arrow["u", from=2-1, to=2-2]
	\arrow["v", from=2-2, to=2-3]
	\arrow["\ulcorner"{anchor=center, pos=0.125, rotate=225}, draw=none, from=2-3, to=1-2]
\end{tikzcd}\] Here $h$, $\tilde{h}$ and $\epsilon$ are defined analogously to \autoref{noneqexmp}.
\end{exmp}

\begin{defn}[\cite{elmanto2023distributivity}]
    A \emph{bispan triple} consists of an $\infty$-category $\mathcal C$ together with two subcategories $\mathcal C_F, \mathcal C_L$ such that the followings are satisfied:
    \begin{enumerate}
        \item $(\mathcal C, \mathcal C_F)$ is a span pair.
        \item $(\mathcal C, \mathcal C_L)$ is a span pair.
        \item For every $l\colon X \to Y \in \mathcal C_L$ and every $f \colon Y\to Z \in \mathcal C_F$, there exists a distributive diagram as in \autoref{distdiagram} such that $g \in \mathcal C_L$.
    \end{enumerate}
\end{defn}

The category of bispans also borns with a beautiful universal property. To see this, let us review the definition of span-distributive functors.

\begin{notn}
    Let $(\mathcal C,\mathcal C_F)$ be a span pair, let $\Phi\colon \Span (\mathcal C,\mathcal C_F) \to \mathcal D$ be a functor. For $f\colon X \to Y \in \mathcal C$, we denote
    \[
    f^*:= \Phi(Y \xleftarrow{f} X \xrightarrow{\mathrm{id}} X)\colon F(Y) \to F(X)
    \] and denote
    \[
    f_{\otimes}:= \Phi(X \xleftarrow{\mathrm{id}} X \xrightarrow{f} Y)\colon F(X) \to F(Y)
    \]
\end{notn}

\begin{defn}\label{spandistributive}
    Let $(\mathcal C,\mathcal C_F,\mathcal C_L)$ be a span triple. A functor $\Phi\colon \Span (\mathcal C,\mathcal C_F) \to \mathcal D$ is called $L$-\emph{span-distributive} if
    \begin{enumerate}
        \item $\Phi|_{\mathcal C^{\mathrm{op}}}$ is left $L$-adjointable. We denote the left adjoint of $f\colon Y \to Z$ by $f_{\oplus} \colon \Phi(Y) \to \Phi (Z)$. 
        \item For any distributive diagram with $l \in L$ and $f \in F$:
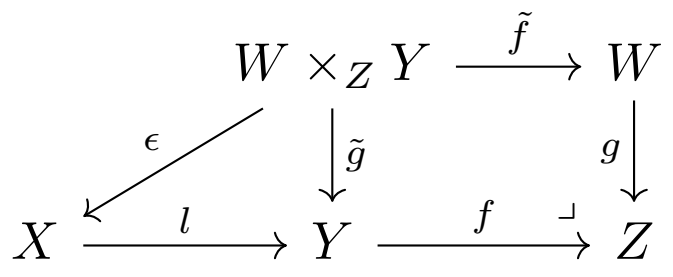
         \[\begin{tikzcd}
	& {W\times_ZY} & W \\
	X & Y & Z
	\arrow["{\tilde f}", from=1-2, to=1-3]
	\arrow["\epsilon"', from=1-2, to=2-1]
	\arrow["{\tilde g}", from=1-2, to=2-2]
	\arrow["g"', from=1-3, to=2-3]
	\arrow["l", from=2-1, to=2-2]
	\arrow["f", from=2-2, to=2-3]
	\arrow["\ulcorner"{anchor=center, pos=0.125, rotate=180}, draw=none, from=2-3, to=1-2]
\end{tikzcd}\] We have a commutative diagram in $\mathcal D$:
\[\begin{tikzcd}
	& {\Phi(W\times_Z Y)} & {\Phi(W)} \\
	{\Phi(X)} & {\Phi(Y)} & {\Phi(Z)}
	\arrow["{\tilde f_{\otimes}}", from=1-2, to=1-3]
	\arrow["{g_\oplus}"', from=1-3, to=2-3]
	\arrow["{\epsilon^*}", from=2-1, to=1-2]
	\arrow["{l_{\oplus}}", from=2-1, to=2-2]
	\arrow["{f_{\otimes}}", from=2-2, to=2-3]
\end{tikzcd}\]
    \end{enumerate}
\end{defn}

\begin{rem}
    We use the term ``bispan distributive '' instead of distributive to emphasis the context and to distinguish from the distributivity we defined using $\infty$-operad pairs in \autoref{operadpairdistributive}. We will see that under good assumptions, these two definitions surprisingly coincide. 
\end{rem}

\begin{thm}[\cite{elmanto2023distributivity}]\label{spanbispan}
    Suppose $(\mathcal C_,\mathcal C_F)$ is a span pair.
    \begin{enumerate}
        \item $(\mathcal C,\mathcal C_F,\mathcal C_L)$ is a span triple if and only if $(\Span (\mathcal C,\mathcal C_F)^{\mathrm{op}},\mathcal C_L)$ is a span pair.
        \item If (1) holds, a functor $\Phi\colon \Span (\mathcal C,\mathcal C_F) \to \mathcal D$ is $L$-bispan distributive if and only if it is $L$-left adjointable.
    \end{enumerate}
\end{thm}
\begin{defn}
    Let $(\mathcal C,\mathcal C_F,\mathcal C_L) $ be a span triple. Then we define
    \[
    \mathrm{Bispan} (\mathcal C,\mathcal C_F,\mathcal C_L) := \Span (\Span(\mathcal C,\mathcal C_F)^{\mathrm{op}}, \mathcal C_L)
    \]
\end{defn}

Informally, the objects of $\mathrm{Bispan} (\mathcal C,\mathcal C_F,\mathcal C_L)$ are the same as the objects of $\mathcal C$. A morphism from $X$ to $Y$ are bispans 
\[\begin{tikzcd}
	& A & B & \\
	X &&& Y
	\arrow["f", from=1-2, to=1-3]
	\arrow["h"', from=1-2, to=2-1]
	\arrow["l", from=1-3, to=2-4]
\end{tikzcd}\] and the composition is formed by using pullback and distributive diagrams.

Let $\Fun_{L\text{-}\mathrm{dist}} (\Span (\mathcal C, \mathcal C_F)^{\mathrm{op}},\mathcal D)$ denote the category of distributive functors. As a consequence of \autoref{spand} and \autoref{spanbispan}, we have the following universal property of bispans:

\begin{cor}\label{universalofbispan}
    Let $(\mathcal C, \mathcal C_F,\mathcal C_L)$ be a bispan triple. There exists an equivalence 
    \[
    \Fun (\mathrm{Bispan} (\mathcal C,\mathcal C_F,\mathcal C_L), \mathcal D) \simeq \Fun_{L\text{-}\mathrm{dist}} (\Span (\mathcal C, \mathcal C_F)^{\mathrm{op}},\mathcal D)
    \] fitting into the following commutative diagram:
    \[\begin{tikzcd}[nodes={scale=0.8}]
	{\Fun(\mathrm{Bispan} (\mathcal C,\mathcal C_F,\mathcal C_L),\mathcal D)} && {\Fun_{L\text{-}\mathrm{dist}}(\Span (\mathcal C,\mathcal C_F)^{\mathrm{op}},\mathcal D)} \\
	& {\Fun (\Span (\mathcal C,\mathcal C_F)^{\mathrm{op}} , \mathcal D)}
	\arrow["\simeq", from=1-1, to=1-3]
	\arrow["j^*"', from=1-1, to=2-2]
	\arrow["u", from=1-3, to=2-2]
\end{tikzcd}\] where $j^*$ is induces by the inclusion $j\colon  \Span (\mathcal C,\mathcal C_F)^{\mathrm{op}} \to \mathrm{Bispan} (\mathcal C ,\mathcal C_F ,\mathcal C_L)$ and $u$ is the forgetful functor.
\end{cor}

Let $(\mathcal C,\mathcal C_F)$ be an operad pair. Suppose $\mathcal C$ admits finite products, \cite{heyer20246} shows that we can endow $\Span (\mathcal C)$ with a symmetric monoidal structure which is obtained from the cartesian monoidal structure of $\mathcal C$. To conclude this section, let us show that when $\mathcal C$ is in addition distributive, we can endow $\Span (\mathcal C,\mathcal C_F)$ with a bimonoidal structure, coming from the cocartesian-cartesian bimonoidal structure of $\mathcal C$.

\begin{con}\label{bimonoidalonspan}
    Let $(\mathcal C, \mathcal C_F)$ be a span pair such that $\mathcal C$ and $\mathcal C_F$ are distributive. Therefore $\mathcal C^{\mathrm{op}}$ and $\mathcal C_F^{\mathrm{op}}$ are codistributive. We let $(\mathcal C^{\mathrm{op},\times,\sqcup})$ be the cartesian-cocartesian bimonoidal category as in \autoref{dualofcocartcart}. We let $\mathcal C_F^-$ denote the category
    \[
    \mathcal C_F^- := (\mathcal C_F^{\mathrm{op},\times,\sqcup})^{\mathrm{op}}\times_{(\mathrm{N}\Fint)^{\mathrm{op}}}(\mathrm{N}\Fint)^{\simeq}
    \] By a same argument as in \cite[lemma 2.3.1]{heyer20246}, $((\mathcal C^{\mathrm{op},\times,\sqcup})^{\mathrm{op}},\mathcal C_F^-)$ is a span pair. We denote $\Span(\mathcal C,\mathcal C_F^-)^{\oplus,\otimes} := \Span ((\mathcal C^{\mathrm{op},\times,\sqcup})^{\mathrm{op}},\mathcal C_F^-)$. Moreover, we have a map $((\mathcal C^{\mathrm{op},\times,\sqcup})^{\mathrm{op}},\mathcal C_F^-) \to ((\mathrm{N}\Fint)^{\mathrm{op}}, (\mathrm{N}\Fint)^{\mathrm{\simeq}})$. This induces a functor
    \[
    p \colon \Span(\mathcal C,\mathcal C_F^-)^{\oplus,\otimes} \to \Span((\mathrm{N}\Fint)^{\mathrm{op}}, (\mathrm{N}\Fint)^{\mathrm{\simeq}}) \simeq \mathrm{N}\Fint
    \]
\end{con}

\begin{thm}
    Let $(\mathcal C,\mathcal C_F)$ be an operad pair. Suppose that $\mathcal C$ and $\mathcal C_F$ are distributive. Then:
\begin{enumerate}
    \item $p\colon \Span(\mathcal C,\mathcal C_F^-)^{\oplus,\otimes} \to \mathrm{N}\Fint$ is a symmetric bimonoidal category.
    \item The underlying operad of $\Span(\mathcal C,\mathcal C_F^-)^{\oplus,\otimes}$ is the cocartesian operad (and of course the cartesian operad) of $\Span (\mathcal C,\mathcal C_F)$ and the underlying moperad coincides with the operad 
    \[\Span(\mathcal C,\mathcal C_F)^\otimes := \Span ((\mathcal C^{\mathrm{op}})^{\sqcup,\mathrm{op}},(\mathcal C_F^{\mathrm{op}})^{\sqcup,\mathrm{op}}\times_{\mathrm{N}\Finp^{\mathrm{op}}} \mathrm{N}\Finp^{\simeq} )\] constructed in \cite[proposition 2.3.3]{heyer20246}.
\end{enumerate}
\end{thm}
\begin{proof}
    For the first claim, we need to verify the axioms of operad pairs and the cocartesian lifting of active morphisms. It suffices to show the distributive part since the remaining follows from \cite[Proposition 2.3.3]{heyer20246}. For any distributive map $d_Q\colon 1^Q_\times \to 1^Q_+ \in \Fint$ and any $C\in \Span(\mathcal C,\mathcal C_F^-)^{\oplus,\otimes}$, since $\mathcal C^{\mathrm{op},\times,\sqcup}$ is an operad, we can find a cocartesian lift $D \colon C\to C' \in \mathcal C^{\mathrm{op},\times,\sqcup}$. Then $C \xleftarrow{D'} C' \xrightarrow{\mathrm{id}} C'$ is the desired cocartesian lift of $d$ to $\Span(\mathcal C,\mathcal C_F^-)^{\oplus,\otimes}$. The second claim follows from the fact that the cocartesian monoidal structure over $\Span (\mathcal C,\mathcal C_F)$ is the same as that of $\mathcal C$ and the \autoref{dualofcocartcart}.
\end{proof}
\begin{rem}
    Let $\mathsf{Fin}$ be the category of finite sets. Then by \autoref{twooneoperad}, the above symmetric bimonoidal structure over $\Span (\mathsf{Fin})$ induces a symmetric monoidal structure over $\Fun^{\oplus}(\Span (\mathsf{Fin}), \mathcal C)$ for any cartesian closed monoidal category $\mathcal C$. Under the equivalence $\Fun^{\oplus}(\Span(\mathsf{Fin}), \mathcal C) \simeq \mathrm{CMon} (\mathcal C)$, this symmetric monoidal structure coincide that of \cite{gepner2016universality} by uniqueness.
\end{rem}
\section{Spectral Tambara functors}\label{distcomparesection}

In \autoref{operadpairdistributive}, we define distributive functors using Day convolution between a symmetric bimonoidal category and a cartesian symmetric monoidal category. In the context of spans and bispans, we also have a definition of distributivity (\autoref{spandistributive}). Now we are going to identify these two definitions.

\begin{thm}\label{maintheorem}
    Let $\mathcal D$ be a presentable $\infty$-category and let $\mathcal D^{\times}$ be its cartesian symmetric monoidal operad. Let $(\mathcal C,\mathcal C_F,\mathcal C_L)=(\mathsf{Fin},\mathsf{Fin},\mathsf{Fin})$ be the span triple. Then there is an equivalence
    \[
    \Fun_{L\text{-}\mathrm{dist}}^\times (\Span (\mathsf{Fin}),\mathcal D) \simeq \Fun^{\mathrm{dist}}(\Span (\mathsf {Fin}),\mathcal D) 
    \] where we equip $\Span (\mathsf{Fin})$ with the symmetric bimonoidal structure constructed in \autoref{bimonoidalonspan}. 
\end{thm}

\begin{proof}
    We construct the desired equivalence as follows. Suppose we have $F\in \Fun^{\mathrm{dist}}(\Span (\mathsf {Fin}),\mathcal D)$, by \autoref{explanationfordist}, $F$ is lax monoidal. The algebraic structures give functorial transformations
\[
\alpha_i\colon\prod_i F(A_i) \to F(\prod_i A_i)
\] and it satisfies the distributivity law. For any maps $\epsilon,l\in \mathsf{Fin}$, we denote $$\epsilon^*:= F|_{\mathsf{Fin}^{\mathrm{op}}} (\epsilon)$$ and $$l_{\oplus} := F|_{\mathsf{Fin}}(l)$$ For any map $v\colon \mathbf{m} \to \mathbf{n} \in \mathsf{Fin}$, we define $v_\otimes \colon F(\mathbf {m}) \to F(\mathbf{n})$ as the composition

\[
 F(\mathbf{m}) \simeq \prod_m F(\mathbf{1}) \simeq \prod_{n}\prod_{v^{-1}(n)} F(\mathbf{1}) \xrightarrow{\alpha_{v^{-1}(n)}} \prod_n F(\prod_{v^{-1}(n)} \mathbf{1}) \simeq \prod_n F(\mathbf 1) \simeq F(\mathbf{n})
\] The first and the last equivalence is from the assumption that $F$ is cartesian symmetric monoidal and $\alpha$ is the lax-monoidal comparison. Moreover, since $\alpha$ is functorial with respect to morphisms in $\Span (\mathsf{Fin})$, suppose we are given a cartesian square
\[\begin{tikzcd}
	{A'} & {B'} \\
	A & B
	\arrow["{\epsilon'}", from=1-1, to=1-2]
	\arrow["{v'}"', from=1-1, to=2-1]
	\arrow["v", from=1-2, to=2-2]
	\arrow["\epsilon"', from=2-1, to=2-2]
	\arrow["\ulcorner"{anchor=center, pos=0.125, rotate=180}, draw=none, from=2-2, to=1-1]
\end{tikzcd}\] in $\mathsf{Fin}$, we have an equivalence $v'_{\otimes}\epsilon'^*\simeq\epsilon^*v_{\otimes}$. Therefore, $(-)^*$ and $(-)_{\otimes}$ together define a left-adjointable functor. As a consequence, this corresponds to a functor $\Phi F \in \Fun (\Span(\mathsf{Fin}), \mathcal D)$ by \autoref{spand} and it is clear such a correspondence is functorial. 

We need to show that $\Phi F$ is $\mathrm{All}$-span-distributive. This is equivalent to say that any distributive diagram as in \autoref{noneqexmp}:

 \[\begin{tikzcd}
	& {\coprod_{j\in J}\prod_{j'\in J_{v(j)}}I_{j'}} & {\coprod_{k\in K}\prod_{j\in J_k}I_j} \\
	I & J & K
	\arrow["{\tilde v}", from=1-2, to=1-3]
	\arrow["\epsilon"', from=1-2, to=2-1]
	\arrow["{\tilde h}", from=1-2, to=2-2]
	\arrow["h", from=1-3, to=2-3]
	\arrow["u", from=2-1, to=2-2]
	\arrow["v", from=2-2, to=2-3]
	\arrow["\ulcorner"{anchor=center, pos=0.125, rotate=180}, draw=none, from=2-3, to=1-2]
\end{tikzcd}\] is sent to a commutative diagram in $\mathcal D$:

\[\begin{tikzcd}
	& {F(\coprod_{j\in J}\prod_{j'\in J_{v(j)}}I_{j'})} & {F(\coprod_{k\in K}\prod_{j\in J_k}I_j)} \\
	{F(I)} & {F(J)} & {F(K)}
	\arrow["{\tilde v_{\otimes}}", from=1-2, to=1-3]
	\arrow["{h_{\oplus}}", from=1-3, to=2-3]
	\arrow["{\epsilon^*}", from=2-1, to=1-2]
	\arrow["{u_{\oplus}}", from=2-1, to=2-2]
	\arrow["{v_{\otimes}}", from=2-2, to=2-3]
\end{tikzcd}\] Without lost of generality, we can assume that $K=\mathbf {1}$. Since $\epsilon$ is defined by coproduct diagonal, $F$ sends coproduct in $\mathsf{Fin}$ to products in $\mathcal D$, $\epsilon^*$ must be a diagonal morphism $\Delta$. Unpack the definition, we consider the following diagram:
\[\begin{tikzcd}
	{\prod_IF(\mathbf{1})} & {\prod_J F(\prod_{j\in J} I_j)} & \\
	{\prod_J F(\mathbf{1})} & {\prod_{J}F(I_j)} & {F(\prod_{j\in J}I_j)} \\
	& {\prod_J F(\mathbf{1})} & {F(\mathbf{1})}
	\arrow["\Delta", from=1-1, to=1-2]
	\arrow["{u_{\oplus}}"', from=1-1, to=2-1]
	\arrow["\simeq", from=1-1, to=2-2]
	\arrow["\alpha", from=1-2, to=2-3]
	\arrow["\simeq"', from=2-1, to=3-2]
	\arrow["\alpha", from=2-2, to=2-3]
	\arrow["{u_{\oplus}}"', from=2-2, to=3-2]
	\arrow["{h_{\oplus}}", from=2-3, to=3-3]
	\arrow["\alpha"', from=3-2, to=3-3]
\end{tikzcd}\] On the outer hexagon, the clockwise composition is $h_\oplus \tilde v_{\otimes} \epsilon^*$ and the counterclockwise composition is $v_{\otimes} u_{\oplus}$. It suffices to show that the three small squares commute. The upper square commutes by the distributivity of $F$ (see \autoref{explanationfordist}). The left diagram commutes since $F$ is cartesian monoidal and the right diagram commutes by functoriality of $\alpha$. This finishes the proof that $\Phi F$ is $\mathrm{All}$-span distributive and gives one direction
\[
\Phi\colon \Fun^{\mathrm{dist}}(\Span (\mathsf {Fin}),\mathcal D) \to \Fun_{L\text{-}\mathrm{dist}}^\times (\Span (\mathsf{Fin}),\mathcal D)
\]

Now suppose conversely we are given an $\mathrm{All}$-span-distributive functor $$G \in \Fun_{L\text{-}\mathrm{dist}}^\times (\Span (\mathsf{Fin}),\mathcal D)$$ we want to show that $G$ defines a commutative algebra object in $\Fun^{\oplus}(\mathcal P,\mathcal C)$. Let $\Fun^{\oplus}(\Span (\mathsf{Fin}),\mathcal D)^{\otimes} \subset \mathrm{Alg}(\Span (\mathsf{Fin}),\mathcal D)^{\otimes} $ denote the operad defined as in \autoref{twooneoperad}. Let $\star$ denote the operadic operations. Then for every $\langle n \rangle \in \Finp$, we need to construct the desired morphism
\[
\tilde \beta\colon G \star\cdots\star G \to G \in \Fun^{\oplus}(\Span (\mathsf{Fin}),\mathcal D)^{\otimes}
\] lifting $\beta \colon \langle n \rangle \to \langle 1 \rangle$. 

Recall that by construction, $\Fun^{\oplus}(\Span (\mathsf{Fin}),\mathcal D)^{\otimes} \subset \Fun (\Span (\mathsf{Fin})^{\oplus},\mathcal D^\times)^{\otimes}$ where the latter is the Day convolution symmetric monoidal operad. By \cite[Corollary 2.2.6.12]{lurie2017higher}, such a morphism $\tilde \beta$ lies in
\[
\mathrm{Mul}^{\beta}(\{G\}_n,G) \simeq \Map_{\Fun((\Span (\mathsf{Fin})^{\oplus})^{\times n},\mathcal D^\times)} (\bigotimes \circ \prod_n G , G\circ \bigotimes)
\] Therefore, $\tilde \beta$ is a natural transformation
\[\begin{tikzcd}
	& {\mathcal D^{\times} \times \cdots\times\mathcal D^{\times}} & \\
	{\Span(\mathsf{Fin})^{\oplus} \times\cdots\times\Span(\mathsf{Fin})^{\oplus}} & {\Span(\mathsf{Fin})^{\oplus}} & {\mathcal D^{\times}}
	\arrow[between={0.3}{0.8}, Rightarrow, from=1-2, to=2-2]
	\arrow["\bigotimes", from=1-2, to=2-3]
	\arrow["{\prod_{n}G}", from=2-1, to=1-2]
	\arrow["\bigotimes"', from=2-1, to=2-2]
	\arrow["G"', from=2-2, to=2-3]
\end{tikzcd}\] Since all functors lie in $\Fun^{\oplus}(\Span (\mathsf{Fin}),\mathcal D)$, $\Span(\mathsf{Fin})^{\oplus}$ is the cartesian operad and $F$ preserves cartesian product, by \cite[Corollary 2.4.1.8]{lurie2017higher}, such a natural transformation is equivalent to a transformation 
\[\begin{tikzcd}
	& {\mathcal D \times \cdots\times\mathcal D} & \\
	{\Span(\mathsf{Fin}) \times\cdots\times\Span(\mathsf{Fin})} & {\Span(\mathsf{Fin})} & {\mathcal D}
	\arrow[between={0.3}{0.8}, Rightarrow, from=1-2, to=2-2]
	\arrow["\prod", from=1-2, to=2-3]
	\arrow["{\prod_{n}G}", from=2-1, to=1-2]
	\arrow["\bigotimes"', from=2-1, to=2-2]
	\arrow["G"', from=2-2, to=2-3]
\end{tikzcd}\] This is exactly a natural transformation $G(A_1) \times \cdots \times G(A_n) \to G(A_1 \times \cdots \times A_n)$.

We now construct the desired natural transformation $\tilde \beta$. Without lost of generality we can do this for $\langle n\rangle = \langle 2 \rangle$. Recall that $G\in \Fun_{L\text{-}\mathrm{dist}}^\times (\Span (\mathsf{Fin}),\mathcal D)$. Therefore for every $\psi\colon A \to B \in \mathsf{Fin}$, it defines a morphism $\psi_\otimes \colon G(A) \to G(B)$. We define
\[
\tilde \beta \colon  G(A) \times G(B) \to G(A\times B) 
\] as the composition 
\[
G(A)\times G(B) \simeq G(A\sqcup B) \xrightarrow{(\Delta_1\sqcup \Delta_2)^*} G(A\times B \sqcup A\times B) \xrightarrow{v_\otimes} G(A\times B)
\] Here $v\colon A\times B\sqcup A\times B \to A\times B$ is the codiagonal, $\Delta_1$ and $\Delta_2$ are given by
\[
\Delta_1 \colon \coprod_A\coprod_B \mathbf{1} \xrightarrow{\coprod_A \Delta} \coprod_A \mathbf{1} \quad \Delta_2 \colon \coprod_B\coprod_A \mathbf{1} \xrightarrow{\coprod_B \Delta} \coprod_B \mathbf{1}
\] We must show that this construction is functorial with respect to morphisms in $\Span (\mathsf{Fin})$. That is, for any map $f\colon A \to B \in \Span (\mathsf{Fin})$, the outer square of the following diagram commute:
\[\begin{tikzcd}
	{G(A) \times G(C)} && {G(A\times C\sqcup A\times C)} & {G(A\times C)} \\
	{G(B)\times G(C)} && {G(B\times C \sqcup B\times C)} & {G(B\times C)}
	\arrow["{(\Delta_1\sqcup\Delta_2)^*}", from=1-1, to=1-3]
	\arrow["{G(f)\times\mathrm{id}}"', from=1-1, to=2-1]
	\arrow["{v^2_\otimes}", from=1-3, to=1-4]
	\arrow["{G(f\times\mathrm{id}\sqcup \mathrm{id} \times f)}"', from=1-3, to=2-3]
	\arrow["{G(f\times \mathrm{id})}", from=1-4, to=2-4]
	\arrow["{(\Delta_1\sqcup\Delta_2)^*}"', from=2-1, to=2-3]
	\arrow["{v^1_\otimes}"', from=2-3, to=2-4]
\end{tikzcd}\] It suffices to show for the case where $G(f) = f^*$ or $G(f)=f_\oplus$. When $G(f)=f^*$, it follows immediately from the Beck-Chevalley condition. It remains to show the case where $G(f)=f_\oplus$. We consider the following diagram

\[\begin{tikzcd}[nodes={scale=0.8}]
	&&&& {G(A\times C)} \\
	{G(A) \times G(C)} && {G(A\times C\sqcup A\times C)} & {G(J')} & {G(K')} \\
	{G(B)\times G(C)} && {G(B\times C \sqcup B\times C)} && {G(B\times C)}
	\arrow["{\Delta^*}", from=1-5, to=2-5]
	\arrow["{(f\times \mathrm{id})_{\oplus}}", curve={height=-30pt}, from=1-5, to=3-5]
	\arrow["{(\Delta_1\sqcup\Delta_2)^*}", from=2-1, to=2-3]
	\arrow["{f_{\oplus}\times\mathrm{id}}"', from=2-1, to=3-1]
	\arrow["{v^2_\otimes}", from=2-3, to=1-5]
	\arrow["{\Delta^*}"', from=2-3, to=2-4]
	\arrow["{(f\times\mathrm{id}\sqcup \mathrm{id} \times f)_{\oplus}}"', from=2-3, to=3-3]
	\arrow["{\tilde v_\otimes}", from=2-4, to=2-5]
	\arrow["{h_\oplus}", from=2-5, to=3-5]
	\arrow["{(\Delta_1\sqcup\Delta_2)^*}"', from=3-1, to=3-3]
	\arrow["{v^1_\otimes}"', from=3-3, to=3-5]
\end{tikzcd}\] The right square is the distributive diagram for \[I= A\times C \sqcup A\times C, \quad  J=B\times C \sqcup B\times C , \quad K=B\times C\] Thus the right square commutes by span-distributivity. The left square and the upper triangle commutes by Beck-Chevalley conditions. This proves the functoriality of the construction of $\tilde \beta$. As a consequence, this gives a functor 

\[
\Psi \colon \Fun_{L\text{-}\mathrm{dist}}^\times (\Span (\mathsf{Fin}),\mathcal D) \to \Fun^{\mathrm{dist}}(\Span (\mathsf {Fin}),\mathcal D) 
\] Finally, comparing the definition of $v_\otimes$ and $\tilde \beta$ on both sides shows that $\Phi \Psi \simeq \mathrm{id}$ and $\Psi \Phi \simeq \mathrm{id}$. \qedhere

\end{proof}

\begin{cor}[\cite{cnossen2025homotopical}]\label{noneqcor}
    Let $\mathcal D$ be a cartesian closed presentable category. Then there is an equivalence
    \[
    \mathrm{CRig}(\mathcal D) \simeq \Fun^{\times}(\mathrm{Bispan}(\mathsf{Fin}),\mathcal D)
    \]
\end{cor}

\begin{proof}
    By definition, we have \[\Fun^{\mathrm{dist}}(\Span (\mathsf {Fin}),\mathcal D) \simeq \mathrm{CAlg} (\Fun^{\oplus}(\Span (\mathsf{Fin}), \mathcal D)) \simeq \mathrm{CAlg} \CMon (\mathcal D) \simeq \mathrm{CRig} (\mathcal D)\] while by \autoref{universalofbispan} and the theorem, we have
    \[
 \Fun^{\times}(\mathrm{Bispan}(\mathsf{Fin}),\mathcal D)\simeq\Fun_{L\text{-}\mathrm{dist}} (\Span (\mathsf{Fin}),\mathcal D) \simeq \Fun^{\mathrm{dist}}(\Span (\mathsf {Fin}),\mathcal D)
    \qedhere \]
\end{proof}

\begin{cor}
     Let $\mathcal D$ be a presentable $\infty$-category and let $\mathcal D^{\times}$ be its cartesian symmetric monoidal operad. Let $G$ be a finite group. Let $(\mathcal C,\mathcal C_F,\mathcal C_L)=(\mathsf{Fin}_G,\mathsf{Fin}_G,\mathsf{Fin}_G)$ be the span triple. Then there is an equivalence
    \[
    \Fun_{L\text{-}\mathrm{dist}}^\times (\Span (\mathsf{Fin}_G),\mathcal D) \simeq \Fun^{\mathrm{dist}}(\Span (\mathsf {Fin}_G),\mathcal D)
    \]
\end{cor}

\begin{proof}
    The key property of $\mathsf{Fin}$ we used when proving \autoref{maintheorem} is that $\mathsf{Fin}$ is generated by non-unitary cartesian monoids and admits distributive diagram, which is also shared by $\mathsf{Fin}_G$.
\end{proof}

\begin{cor}[\cite{cnossen2024normed}]\label{stf}
    Let $G$ be a finite set. Then we have an equivalence:
    \[
    \mathrm{NAlg}(\mathrm{Sp}_G^{\geq 0}) \simeq \Fun^{\times}_{\mathrm{gp}} (\mathrm{Bispan} (\mathsf{Fin}_G), \mathcal S)
    \]
\end{cor}
\begin{proof}
    We can replace $\mathsf{Fin}$ by $\mathsf{Fin}_G$ and repeat the proof of \autoref{noneqcor}.
\end{proof}

\printbibliography
\end{document}